\documentclass{article}
\usepackage[utf8]{inputenc}
\usepackage{amsmath, amssymb, amsthm, color,dsfont, mathtools, yfonts, amsfonts, eqparbox}
\usepackage[all]{xy}
\usepackage{calligra, mathrsfs}
\usepackage[nodisplayskipstretch]{setspace}
\usepackage{thmtools}
\usepackage{thm-restate}
\usepackage{verbatim}
\usepackage{relsize}
\usepackage{bm}
\usepackage{tasks}
\usepackage{exsheets}
\usepackage{faktor}
\usepackage{breqn}
\usepackage{tikz-cd}
\usepackage[mathlines]{lineno}
\usepackage[bookmarks]{hyperref}

\usepackage{graphicx, floatrow, wrapfig}
\usepackage{xstring}

\newcommand{\myincludesvg}[3]{\includegraphics[#1]{svg-inkscape/#3_svg-raw.png}}
\newcommand{\myincludesvggroup}[4]{\includegraphics[#1]{svg-inkscape/#3-#4_svg-raw.png}}

\newcommand{\diagramhh}[4]{
\raisebox{-#3}{\myincludesvggroup{scale=0.25}{inkscape=nolatex,inkscapeformat=png,inkscapename=#1-#2,inkscapeopt=-i #2 -j}{#1}{#2}}\hspace{#4}}
\newcommand{\smalldiagram}[1]{\diagramhh{skeinrelations}{#1}{10pt}{1pt}}

\usepackage{fullpage}
\allowdisplaybreaks

\usepackage{url}
\usepackage[backend=bibtex ,style=alphabetic, doi=false,isbn=false,url=false, maxcitenames = 99, maxbibnames = 99]{biblatex}
\usepackage[noabbrev, capitalise]{cleveref}

\newtheorem{theorem}{Theorem}[section]
\newtheorem{lemma}[theorem]{Lemma}
\newtheorem{proposition}[theorem]{Proposition}
\newtheorem{corollary}[theorem]{Corollary}
\newtheorem{example}[theorem]{Example}

\newtheorem{definition}[theorem]{Definition}
\newtheorem*{theorem*}{Theorem}
\newtheorem*{proposition*}{Proposition}
\newtheorem*{corollary*}{Corollary}
\theoremstyle{definition}
\theoremstyle{remark}
\newtheorem{remark}[theorem]{Remark}
\theoremstyle{remark}

\newcommand{\catname}[1]{{\normalfont\textbf{#1}}}

\newcommand{\Cocont}{\catname{Cocont}}
\newcommand{\Rep}{\catname{Rep}}
\newcommand{\Repfd}{\catname{Rep}^{\mathrm{fd}}}

\newcommand{\Mfld}[1]{\catname{Mfld}^{\mathrm{or}, \sqcup}_{#1}}
\newcommand{\Surf}{{\catname{Surf}\,}^{\mathrm{or}, \sqcup}}
\newcommand{\Disc}[1]{\catname{Disc}^{\mathrm{or}, \sqcup}_{#1}}

\newcommand{\kMod}{k\catname{Mod}}

\newcommand{\LFP}[1]{\catname{LFP}_{#1}}

\newcommand{\Rex}{\mathscr{R}\!\mathit{ex}}
\newcommand{\uChg}{\underline{\operatorname{Ch}}_{G}}
\newcommand{\Chg}{\operatorname{Ch}}
\newcommand{\SL}{\mathrm{SL}_2}

\newcommand{\Chs}{\operatorname{Ch}_{\SL}}
\newcommand{\Lie}{\operatorname{Lie}}
\newcommand{\End}{\operatorname{End}}
\newcommand{\ch}{\operatorname{ch}}

\newcommand{\Mod}[2]{{#1-\operatorname{mod}}_{#2}}
\newcommand{\Hom}{\operatorname{Hom}}

\newcommand{\Tr}{\operatorname{tr}_q}

\newcommand{\Id}{\operatorname{Id}}
\newcommand{\sk}{\operatorname{Sk}}
\newcommand{\IHom}{\underline{\operatorname{Hom}}}

\newcommand{\frsl}{\mathfrak{sl}_2}
\newcommand{\frg}{\mathfrak{g}}
\newcommand{\qgroup}[1]{\mathcal{U}_q(#1)}
\newcommand{\hgroup}[1]{\mathcal{U}_{\hbar}(#1)}

\newcommand{\mystackrel}[3][T]{\stackrel{\eqmakebox[#1]{\scriptsize#2}}{#3}}

\usepackage[inline]{enumitem}
\usepackage[super]{nth}
\newlist{stages}{enumerate}{1}
\setlist[stages]{label=\Roman*.}

\title{Kauffman Skein Algebras and Quantum Teichm\"uller Spaces via Factorisation Homology}
\author{Juliet Cooke}
\date{\today}

\begin{document}
\maketitle

\begin{abstract}
We  compute the factorisation homology of the four-punctured sphere and punctured torus over the quantum group \(\mathcal{U}_q(\mathfrak{sl}_2)\) explicitly as categories of equivariant modules using the framework developed by Ben-Zvi, Brochier, and Jordan. We identify the algebra of \(\mathcal{U}_q(\mathfrak{sl}_2)\)-invariants (quantum global sections) with the spherical double affine Hecke algebra of type $(C^{\vee}_1,C_1)$, in the four-punctured sphere case, and with the `cyclic deformation' of $U(su_2)$ in the punctured torus case. In both cases, we give an identification with the corresponding quantum Teichm\"uller space as proposed by Teschner and Vartanov as a quantization of the moduli space of flat connections.
\\ \\
\emph{Quantum Algebra:} Factorization Homology; Skein Algebras; Quantum Groups; Teichmüller Spaces, Character Varieties
\\ \\
Mathematics Subject Classification 2000: 	57K31, 16T99
\end{abstract}

\setcounter{tocdepth}{3}

\addcontentsline{toc}{section}{Introduction}
\section*{Introduction}
Factorisation homology theories of topological manifolds are generalised homology theories of manifolds whose coefficients systems are $n$-disc algebras. They may be interpreted as homology theories which are tailor-made for topological manifolds rather than general topological spaces as they satisfy a generalisation of the Eilenberg-Steenrod axioms for singular homology \cite{AyalaFrancis}. They are simultaneously an attempt to axiomatise the structure of observables in Topological Quantum Field Theories (TQFTs) which obey the strong locality principle that local observables determine global observables, and can be used to construct extended TQFTs \cite{Lurie, Scheimbauer}. Factorisation homology was first defined by Beilinson and Drinfeld in the conformal setting, and was subsequently developed in a topological setting by Lurie~\cite{Lurie}, Ayala, Francis, and Tanaka~\cite{AyalaFrancis, AyalaFrancisTanaka}. 

Given a quantum group $\qgroup{\mathfrak{g}}$ associated to a reductive algebraic group $G$, one can use the category $\Rep_q(G)$ of integrable representations of this quantum group as a coefficient system for factorisation homology. This factorisation homology $\int_{\Sigma} \Rep_q(G)$ for a surface $\Sigma$ has been studied by Ben-Zvi, Brochier and Jordan~\cite{david1, david2}. In particular, when $\Sigma$ is a punctured surface, this factorisation homology is the category of modules over an algebra $A_{\Sigma}$. The algebra $A_{\Sigma}$ is the Aleesev moduli algebra \cite{Alek94, AlekseevGrosseSchomerus96} and is determined combinatorially. Taking the subalgebra of $\hgroup{\mathfrak{sl}_2}$-invariants of this moduli algebra $A_{\Sigma}$ gives a deformation quantisation of the character variety $\Chg_G(\Sigma)$: the character variety $\Chg_G(\Sigma)$ is the moduli space of representations of the fundamental group of $\Sigma$ into $G$ and carries a canonical Atiyah--Bott--Goldman Poisson structure \cite{AtiyahBott83,Goldman84}. 

In this paper we shall concern ourselves with this moduli algebra $A_{\Sigma}$ in the cases of the four punctured sphere $\Sigma_{0,4}$ and the punctured torus $\Sigma_{1,1}$ with gauge group $G = \SL$. We shall give Poincaré–Birkhoff–Witt (PBW) bases for the algebras $A_{\Sigma_{0,4}}$ and $A_{\Sigma_{1,1}}$ before turning to our main technical result: obtaining explicit generators and relations presentations and PBW bases of the algebras  of $\mathscr{A}_{\Sigma_{0,4}}$ and $\mathscr{A}_{\Sigma_{1,1}}$ of $\qgroup{\mathfrak{sl}_2}$-invariants of these moduli algebras and hence of quantisations of the character varieties $\Chs(\Sigma_{0,4})$ and $\Chs(\Sigma_{1,1})$.

Another approach to quantising character varieties is by using skein algebras. In this paper we exhibit explicit isomorphisms between  $\mathscr{A}_{\Sigma_{0,4}}$ and $\mathscr{A}_{\Sigma_{1,1}}$ and the Kauffman bracket skein algebra of $\Sigma_{0,4}$ and $\Sigma_{1,1}$ thus relating these two approaches to quantising character varieties. During the publication process of this manuscript, there have been several papers showing the existence of an isomorphism more generally. In our followup paper~\cite{Cooke2019} we prove that there exists an isomorphism for any punctured surface $\Sigma$ and any $G$ (with generic deformation parameter). This isomorphism is lifted to an isomorphism between the internal skein algebra and the moduli algebra in~\cite{GJS2019}. For $G = \SL$ there is an alternate but related approach which proves there is an isomorphism between the moduli algebra $A_{\Sigma}$ and stated skein algebra of $\Sigma$ : see Lê and Costantino~\cite{CostantinoLe2019} or Baseilhac and Roche~\cite{BR2019} for $n$-punctured spheres and Faitig~\cite{Faitg2020} for general punctured surfaces. From this isomorphism by taking invariants, one also proves the existence of an isomorphism between the skein algebra and the invariant subalgebra of the moduli algebra. Whilst these papers prove the existence of an isomorphism in general, the results in this paper are more explicit allowing one to directly compare the presentations for the examples considered. 

The isomorphism between  $\mathscr{A}_{\Sigma_{0,4}}$ and $\mathscr{A}_{\Sigma_{1,1}}$ and the Kauffman bracket skein algebra of $\Sigma_{0,4}$ and $\Sigma_{1,1}$ also leads to an isomorphism between $\mathscr{A}_{0,4}$ and the spherical double affine Hecke algebra of type $(C^{\vee}_1, C_1)$ \cite{sahi1999nonsymmetric, Terwilliger, Berest&Samuelson}, and an isomorphism between $\mathscr{A}_{\Sigma_{1,1}}$ and the cyclic deformation of $U(\mathfrak{su}_2)$ \cite{Bullock&Przytycki, zachos1990quantum}. 

We conclude the paper by exhibiting isomorphisms between $\mathscr{A}_{\Sigma_{0,4}}$ and $\mathscr{A}_{\Sigma_{1,1}}$ and a quantisation of the $\SL$-character variety of $\Sigma_{0,4}$  and $\Sigma_{1,1}$ proposed by Teschner and Vartanov~\cite{Teschner}. This, in particular, shows that the constructions of \cite{Teschner}, which are given by generators and relations, are isomorphic to the output of a functorial construction and fit into the framework of fully extended TQFTs. This may also be useful for generalising from $\SL$ to other gauge groups.

\subsection*{Summary of Sections and Results}
\begin{description}
\item [Section 1:] In the background section we give a brief introduction to the \(\LFP{k}\) factorisation homology of oriented surfaces. We also define the algebra of \(\qgroup{\mathfrak{g}}\)-invariants. We conclude this section by recalling the definitions of reduction systems, PBW bases and the Diamond lemma.

\item[Section 2:] We give a presentation for \(\mathscr{A}_{\Sigma_{0,4}}\)  and \(\mathscr{A}_{\Sigma_{1,1}}\) thus giving an explicit presentation for the quantum \(SL_2\)-character variety of the four-punctured sphere and punctured torus.
\begin{theorem}
The algebra of \(\qgroup{\mathfrak{sl}_2}\)-invariants \(\mathscr{A}_{\Sigma_{0,4}}\) of the four-punctured sphere has a presentation given by generators \(E := \Tr(AB),\; F := \Tr(AC),\; G := \Tr(BC),\; s := \Tr(A),\; t := \Tr(B),\; u := \Tr(C)\; v := \Tr(ABC),\)
\footnote{where \(A = \left(\begin{smallmatrix} a_{11} & a_{12} \\ a_{21} & a_{22}\end{smallmatrix}\right)\) \(B = \left(\begin{smallmatrix} b_{11} & b_{12} \\ b_{21} & b_{22}\end{smallmatrix}\right)\) and \(C = \left( \begin{smallmatrix} c_{11} & c_{12} \\ c_{21} & c_{22}\end{smallmatrix} \right)\) are matrices formed out of the 12 generators of the moduli algebra \(A_{\Sigma_{0,4}}\); they correspond to loops punctures as depicted for \(E, F, G\) in \cref{figure:FPuncSphereEdges}.}
and relations 
\newcommand{\lastgen}{
\begin{cases}
\begin{aligned}
        & -E^2 - q^{-4}F^2 - G^2 - q^{-4}(s^2 + t^2 + u^2 + v^2) \\
        & + (st + uv)E + q^{-2}(su +tv)F + (sv +tu)G \\
        & -stuv + q^{-6}(q^2+1)^2
\end{aligned}
\end{cases}
}
\begin{alignat*}{8}
    &FE &{}={}&&
        q^2 &EF &{}+{}&&
        (q^2 - q^{-2})&G &{}+{}&&
        (1- q^2)&(sv + tu), \\
    &GE &{}={}&&
        q^{-2}&EG &{}+{}&&
        q^{-2}(q^2-q^{-2})&F &{}-{}&&
        (1-q^2)&(su + q^{-2}tv), \\
    &GF &{}={}&&
        q^2 & FG &{}+{}&&
        (q^2 - q^{-2})&E &{}+{}&&
        (1-q^2)&(st + uv), \\
    & EFG &= \mathclap{\hphantom{\lastgen}\lastgen}
\end{alignat*}
and \(s,t,u,v\) are central. Furthermore, the monomials 
\[
\{E^{m} F^n G^l s^a t^b u^c v^d : m, n, l, a, b, c, d \in \mathbb{N}^{0}; m\text{ or }n\text{ or }l = 0\}
\]
are a basis for the algebra.
\end{theorem}

\begin{theorem}
The algebra of \(\qgroup{\mathfrak{sl}_2}\)-invariants \(\mathscr{A}_{\Sigma_{1,1}}\) of the punctured torus with respect to \(\qgroup{\mathfrak{sl}_2}\) has a presentation given by generators \(X :=\Tr(A), Y:=\Tr(B), Z := \Tr(AB)\)and relations:
\begin{align*}
YX - q^{-1}XY &= (q-q^{-1})Z; \\
XZ - q^{-1}ZX &= - q^{-3}(q-q^{-1})Y; \\
ZY - q^{-1}YZ &= -q^{-3}(q - q^{-1})X.
\end{align*}
It has a central element 
\[L := q^5 XZY + q^3 Y^2 - q^4 Z^2 + q^3 X^2 - (q-q^{-1}).\]
and a PBW basis \[\{X^{\alpha} Y^{\beta} Z^{\gamma}: \alpha, \beta, \gamma\in \mathbb{N}^{0}\}.\] 
\end{theorem}

\item [\cref{section:DAHA}:] In this section we use the presentation for \(\mathscr{A}_{\Sigma_{0,4}}\) and \(\mathscr{A}_{1, 1}\) to prove them isomorphic to the Kauffman bracket skein algebras of \(\Sigma_{0, 4}\) and \(\Sigma_{1,1}\):
\begin{proposition}
The algebra of \(\qgroup{\mathfrak{sl}_2}\)-invariants \(\mathscr{A}_{\Sigma_{0,4}}\) is isomorphic to the Kauffman bracket skein algebra \(\sk(\Sigma_{0,4})\) with isomorphism \(\beta: \sk_q(\Sigma_{0,4})\) \(\to \mathscr{A}_{\Sigma_{0,4}}\) given by 
\begin{alignat*}{3}
    \beta(x_1) &= -q E,&\quad& \beta(p_1) &\;=& \; -q s,\\
    \beta(x_2) &= -q F,&& \beta(p_2) &=&\; - q t,\\
    \beta(x_3) &= -q G,&& \beta(p_3) &=&\; - q v,\\
   \beta(q) &= q^2, && \beta(p_4) &=&\; - q u.
\end{alignat*}
There is also an isomorphism  \(\gamma: \mathscr{A}_{\Sigma_{1,1}} \to \sk(\Sigma_{1,1})\) given by 
\begin{align*}
\gamma(q) &= q^2 \\
\gamma(X) &= i q^{-2} x_2 \\
\gamma(Y) &= i q^{-2} x_1 \\
\gamma(Z) &= -q^{-5} x_3.
\end{align*}
\end{proposition}
As a consequence we obtain explicit isomorphisms \(\mathscr{A}_{\Sigma_{0,4}} \cong S\mathscr{H}_{q,\underline{t}}\) to the spherical double affine Hecke algebra of type \((C^{\vee}_1, C_1)\) and \(\mathscr{A}_{1,1} \cong U_1(\mathfrak{su}_2)\) to the cyclical deformation of \(U(\mathfrak{su}_2)\). 

\item [\cref{section:FlatConnections}:] We summarise the paper of Teschner and Vartanov~\cite{Teschner} giving their definition of \(\mathscr{A}_b(\Sigma)\), a non-commutative deformation of the Poisson-algebra of algebraic functions on the moduli space of flat connections, and then prove
\begin{proposition}
The algebra of \(\qgroup{\mathfrak{sl}_2}\)-invariants \(\mathscr{A}_{\Sigma_{0,4}}\) is isomorphic to \(\mathscr{A}_b(S)\) with isomorphism \(\iota: \mathscr{A}_{\Sigma_{0,4}} \to \mathscr{A}_b(S)\) given by
\begin{alignat*}{3}
    \iota(q) &= e^{i \pi b^2}, && \iota(s) &= e^{-i \pi b^2} L_1, \\
    \iota(E) &= -e^{-i \pi b^2} L_u, &\quad& \iota(t) &= e^{-i \pi b^2} L_3, \\
    \iota(F) &= -e^{-i \pi b^2} L_s, &\quad& \iota(v) &= e^{-i \pi b^2} L_2, \\
    \iota(G) &= -e^{-i \pi b^2} L_t, &\quad& \iota(u) &= e^{-i \pi b^2} L_4.
\end{alignat*}
\end{proposition}
\begin{proposition}
The algebra of \(\qgroup{\mathfrak{sl}_2}\)-invariants \(\mathscr{A}_{\Sigma_{1,1}}\) is isomorphic to \(\mathscr{A}_b(\Sigma_{1,1})\) with isomorphism \(\mu: \mathscr{A}_{\Sigma_{1,1}} \to \mathscr{A}_b(\Sigma_{1,1})\) given by
\begin{align*}
\mu(Y) &= iq^{-1}s \\
\mu(X) &= iq^{-1}t \\
\mu(Z) &= -q^{- \frac{5}{2}} u \\
\mu(L) &= L_0 
\end{align*}
\end{proposition}

\end{description}

\section{Background}
\subsection{Factorisation Homology}
\label{section:FactorisationHomology}

We shall begin by defining the factorisation homology of oriented surfaces with coefficients given by a framed \(E_2\)-algebra. General introductory references for factorisation homology include Ginot \cite{Ginot} and Ayala and Francis~\cite{AyalaFrancis, AFPrimer}.

\begin{definition}
A smooth surface \(\Sigma\) is \emph{finitary} if it has a finite open cover \(\mathcal{U}\) such that if \(\{\,U_i\,\}\) is a subset of \(\mathcal{U}\) then the intersection \(\cap_i U_i\) is either empty or diffeomorphic to \(\mathbb{R}^2\). 
\end{definition}

\begin{remark}
Surfaces are assumed throughout this paper to be finitary, smooth and oriented.
\end{remark}

\begin{definition}
Let \(\Surf\) be the symmetric monoidal \((2,1)\)-category whose
\begin{enumerate}
    \item objects are oriented, finitary, smooth surfaces;
    \item \(1\)-morphisms are smooth oriented embeddings between surfaces;
    \item \(2\)-morphisms are isotopies on embeddings;
    \item symmetric monoidal product is disjoint union.
\end{enumerate}
\end{definition}

\begin{remark}
By a \((2, 1)\)-category we mean a strict \(2\)-category for which all \(2\)-morphisms are invertible.
\end{remark}

\begin{definition}
Let \(\Disc{2}\) be the full subcategory of \(\Surf\) of finite disjoint unions of \(\mathbb{R}^2\). Denote the inclusion functor by \(I: \Disc{2} \to \Surf\).
\end{definition}

\begin{definition}
Let \(\mathscr{C}^{\otimes}\) be a symmetric monoidal \((2,1)\)-category.
A \emph{framed \(E_2\)-algebra} in \(\mathscr{C}^{\otimes}\) is a symmetric monoidal functor \(F: \Disc{2} \to \mathscr{C}^{\otimes}\). As \(F\) is determined on objects by its value on a single disc, we define \(\mathscr{E} := F(\mathbb{R}^2)\), and we use \(\mathscr{E}\) to refer to the associated framed \(E_2\)-algebra.  
\end{definition}

\begin{remark}
A framed \(E_2\)-algebra is also known as a \(2\)-disk algebra. The terminology framed \(E_2\)-algebra is somewhat confusing as there is also a notion of an \(E_2\)-algebra which is a symmetric monoidal functor \(F:  \catname{Disc}^{\mathrm{fr}, \sqcup}_{2} \to \mathscr{C}^{\otimes}\) from the category of \emph{framed} discs. Using \(E_2\)-algebra as coefficients one can define factorisation homology for framed surfaces; however, we shall only consider oriented surfaces and framed \(E_2\)-algebras in this paper.
\end{remark}

\begin{remark}
Usually a framed \(E_2\)-algebra is defined where \(\mathscr{C}^{\otimes}\) and \(\Surf\) are \((\infty, 1)\)-categories rather than a \((2,1)\)-categories, but we can treat any \((2,1)\)-category as an \((\infty,1)\)-category with the only \(k\)-morphisms for \(k > 2\) being the identiy morphisms.
\end{remark}

\begin{definition}[{\cite[Definition~3.4]{AyalaFrancis}}]
A symmetric monoidal \((2,1)\)-category \(\mathscr{C}^{\otimes}\) is \emph{\(\otimes\)-presentable} if 
\begin{enumerate}
    \item \(\mathscr{C}\) is locally presentable with respect to an infinite cardinal \(\kappa\) and
    \item the monoidal structure distributes over small colimits i.e.\ the functor \(C \otimes \_ : \mathscr{C} \to \mathscr{C}\) carries colimit diagrams to colimit diagrams.
\end{enumerate}
\end{definition}

\begin{definition}
\label{defn:facthom}
Let \(\mathscr{C}^{\otimes}\) be a \(\otimes\)-presentable symmetric monoidal \((\infty,1)\)-category and let \(F: \Disc{n} \to \mathscr{C}^{\otimes}\) be a framed \(E_2\)-algebra with \(\mathscr{E} := F(\mathbb{R}^2)\). The left Kan extension of the diagram
\[
\begin{tikzcd}
\Disc{2} \arrow[r, "F"] \arrow[d, hook, "I"'] 
				&  \mathscr{C}^{\otimes} \\
\Surf \arrow[ru, dashrightarrow, "\int_{\_} \mathscr{E}"'] 
				& 
\end{tikzcd}
\]
is called the\footnote{As factorisation homology is defined via a universal construction we have uniqueness up to a contractible space of isomorphisms.} \emph{factorisation homology} with coefficients in \(\mathscr{E}\); its image on the surface \(\Sigma\) is called the factorisation homology of \(\Sigma\) over \(\mathscr{E}\) and is denoted \(\int^{\mathscr{C}^{\otimes}}_{\Sigma} \mathscr{E}\) or \(\int_{\Sigma} \mathscr{E}\) when \(\mathscr{C}^{\otimes}\) is clear from context.  
\end{definition}
\subsection{The Category \(\LFP{k}\)}
We shall now define the \((2,1)\)-category \(\LFP{k}\) which will be the ambient category \(\mathscr{C}^{\otimes}\) of the factorisation homologies considered in this paper. A general reference for this section is Borceux's `Handbook of Categorical Algebra' \cite{Handbook1, Handbook2} and we mostly follow the terminology of \cite{david1}.

\begin{definition}
Let \(k\) be a commutative ring with identity.
A \emph{\(k\)-linear category} is a category enriched over \(\kMod\), the category of left \(k\)-modules,  and a \emph{\(k\)-linear functor} is a \(\kMod\)-enriched functor.
\end{definition}

\begin{definition}
A category \(\mathscr{C}\) is \emph{locally finitely presentable} if it is locally small, cocomplete and is generated under filtered colimits by a set of finitely presentable objects.
\end{definition}

\begin{definition}
A \(k\)-linear functor is \emph{cocontinuous} if it preserves all small \(k\)-linear colimits.
\end{definition}

\begin{definition}
Let \(\LFP{k}\) denote the \((2, 1)\)-category with:
\begin{enumerate}
    \item objects: locally finitely presentable \(k\)-linear categories;
    \item \(1\)-morphisms: cocontinuous \(k\)-linear functors;
    \item \(2\)-morphisms: natural isomorphisms.
\end{enumerate}
\end{definition}

The \((2,1)\)-category \(\LFP{k}\) is a strict monoidal category with the monoidal product \(\boxtimes\) given by the Kelly--Deligne tensor product\footnote{The monoidal unit of \(\LFP{k}^{\boxtimes}\) is \(\kMod\).}:

\begin{definition}
Let \(\Cocont(\mathscr{A} \boxtimes \mathscr{B}, \mathscr{C})\) be the category of cocontinuous functors \(\mathscr{A} \boxtimes \mathscr{B} \to \mathscr{C}\) and \( \Cocont(\mathscr{A}, \mathscr{B}; \mathscr{C})\) be the category of bilinear functors \(\mathscr{A} \times \mathscr{B} \to \mathscr{C}\) which are cocontinuous in each variable separately.
The \emph{Kelly--Deligne tensor product} of \(\mathscr{A}, \mathscr{B} \in \LFP{k}\) is a category \(\mathscr{A} \boxtimes \mathscr{B} \in \LFP{k}\) together with a bilinear functor \(S: \mathscr{A} \times \mathscr{B} \to \mathscr{A} \boxtimes \mathscr{B} \in \Cocont(\mathscr{A}, \mathscr{B}; \mathscr{C})\) such that composition with \(S\) defines an equivalence of categories 
\[
\Cocont(\mathscr{A} \boxtimes \mathscr{B}, \mathscr{C}) \simeq \Cocont(\mathscr{A}, \mathscr{B}; \mathscr{C}) \cong \Cocont(\mathscr{A}, \Cocont(\mathscr{B}, \mathscr{C}))
\]
for all \(\mathscr{C} \in \LFP{k}\).
\end{definition}

\begin{remark}
Kelly \cite[Proposition~4.3]{Kelly82} proved the existence of \(\mathscr{A} \boxtimes \mathscr{B}\) for categories \(\mathscr{A}, \mathscr{B} \in \Rex\), the \((2,1)\)-category of essentially small, finitely cocomplete categories with right exact functors as \(1\)-morphisms and natural isomorphisms as \(2\)-morphisms. Franco in \cite[Theorem~18]{Franco2012} shows that for abelian categories \(\mathscr{A}, \mathscr{B}\), this tensor product \(\mathscr{A} \boxtimes \mathscr{B}\) is the Deligne tensor product of abelian categories \cite{Deligne90} when the Deligne tensor product exists; hence, the name Kelly--Deligne tensor product. For the existence of the Kelly--Deligne tensor product in \(\LFP{k}\) see \cite[Section 2.4.1]{TensorThesis} and the references therewithin. 
\end{remark}
\begin{remark}
\(\LFP{k}^{\boxtimes}\) is \(\otimes\)-presentable \cite[Section~4]{KellyLack01}~\cite[7,115]{Kelly05}~\cite[Proposition~3.5]{david1}, thus \(\LFP{k}^{\boxtimes}\) can be used as the ambient category for the factorisation homology.
\end{remark}

\subsection{\(\LFP{k}\) Factorisation Homology of Punctured Surfaces}
Let \(\mathscr{E}\) be an abelian \(k\)-linear compact-rigid balanced tensor category. The primary example of such an \(\mathscr{E}\) is \(\Rep_q(G)\).
\begin{definition}
A locally presentable monoidal category \(\mathscr{E}\) is compact-rigid if all compact objects are left and right dualisable.
\end{definition}

\begin{definition}
Let \(G\) be a connected reductive algebraic group and let \(\qgroup{\frg}\) be the quantum group of the Lie algebra \(\frg = \Lie(G)\). We assume \(q \in \mathbb{C}^*\) is generic. If \(G\) is simply connected let \(\Rep_q(G)\) be the \(\mathbb{C}\)-linear compact-rigid balanced tensor category of (possibly infinite) direct sums of finite-dimensional integrable \(\qgroup{\frg}\)-modules. If \(G\) is not simply connected let \(\Rep_q(G)\) be the subcategory of this compact-rigid balanced tensor category consisting of the \(\qgroup{\frg}\)-modules which correspond to representations of \(G\).
\end{definition}

\begin{remark}
Given the abelian \(k\)-linear compact-rigid balanced tensor category \(\mathscr{E}\) there is a canonical framed \(E_2\) algebra \(F_{\mathscr{E}}: \Disc{n} \to \LFP{k}\) such that \(F_{\mathscr{E}}(\mathbb{D}) = \mathscr{E}\). 
\end{remark}

The factorisation homology \(\int_{\Sigma} \mathscr{E}\) of the punctured surface \(\Sigma\) can be given an \(\mathscr{E}\)-module category structure as follows:
\begin{figure}[H]
    \centering
    \myincludesvg{height=2.5cm}{}{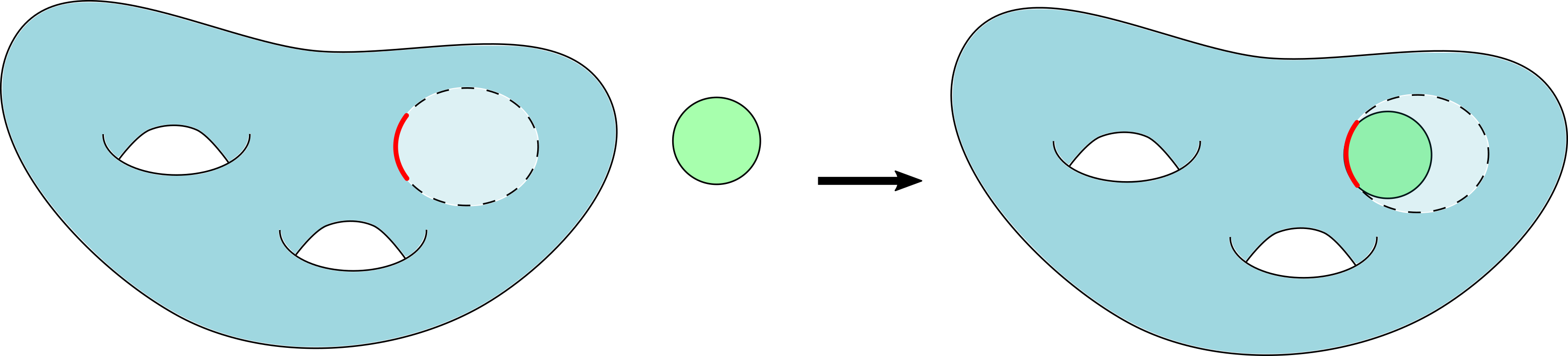}
    \caption{An illustration of the map \(\Sigma \sqcup \mathbb{D} \to \Sigma\). The surface \(\Sigma_{2,1}\) has a interval marked in red along its boundary along which the disc \(\mathbb{D}\) is attached. The resultant surface is isotopic to \(\Sigma_{2,1}\). }
\end{figure}
\noindent
Choose an interval along the boundary of \(\Sigma\)\footnote{The module structure depends on the choice of marking.}. The mapping
\(\Sigma \sqcup \mathbb{D} \to \Sigma,\)
which attaches the disc \(\mathbb{D}\) to \(\Sigma\) along the marked interval, induces a \(\int_{\mathbb{D}} \mathscr{E}\)-module structure on \(\int_{\Sigma} \mathscr{E}\). As \(\int_{\mathbb{D}} \mathscr{E} \simeq \mathscr{E}\) in \(\LFP{k}\), this means that \(\int_{\Sigma} \mathscr{E}\) is a \(\mathscr{E}\)-module. 

Not only is \(\int_{\Sigma} \mathscr{E}\) an \(\mathscr{E}\)-module category, but Ben-Zvi, Brochier and Jordan showed that it is the category of modules over an algebra \(A_{\Sigma}\) in \(\mathscr{E}\). This algebra \(A_{\Sigma}\) is an internal Hom:

\begin{definition}[{\cite[147]{EtingofTensorCategories}}]
Let \(\mathscr{M}\) be a right \(\mathscr{E}\)-module category\footnote{Note that Etignof et al.\ are assuming that \(\mathscr{M}\) is a \emph{left} \(\mathscr{E}\)-module category, whereas we are assuming that is is a \emph{right} \(\mathscr{E}\)-module category. Also note that they assume the categories are finite, but the proofs work without modification for locally finitely presentable categories.} and let \(m, n \in \mathscr{M}\). The \emph{internal Hom}\footnote{Also known as the enriched Hom.} from \(m\) to \(n\) is the object \(\IHom(m, n) \in \mathscr{E}\) which represents the functor \(x \mapsto \Hom_{\mathscr{M}}(m \cdot x, n)\) i.e.\ such that there is a natural isomorphism 
\[\eta^{m, n}: \Hom_{\mathscr{M}}(m \cdot \_\;, n) \xrightarrow{\sim} \Hom_{\mathscr{E}}(\_\;, \IHom(m, n)).\]
For any \(m \in \mathscr{M}\), the \emph{internal endomorphism algebra} \(\End(m) := \IHom(m, m)\) is an algebra object of \(\mathscr{E}\)~\cite[149]{EtingofTensorCategories}.
\end{definition}

\begin{definition}
As \(\emptyset\) is the identity for the monoidal product \(\sqcup\) in \(\Mfld{n}\),
\(
\int_{\emptyset} \mathscr{E} \simeq \kMod
\), the monoidal unit of \(\LFP{k}\).
We can embed the empty manifold into any manifold, and this embedding \(\emptyset \to \Sigma\) induces a morphism \(\kMod \to \int_{\Sigma} \mathscr{E}\). The \emph{distinguished object} \(\mathscr{O}(\mathscr{E})\) of the factorisation homology of \(\Sigma\) over \(\mathscr{E}\) is the image of \(k\) under this map.
\end{definition}

\begin{definition}
\label{defn:facthomalg}
The \emph{algebra object} \(A_{\Sigma}\) of the factorisation homology of \(\Sigma\)\footnote{The algebra object is dependent on the choice marking of \(\Sigma\)} with coefficients in \(\mathscr{E}\) is the internal endomorphism algebra of the distinguished object
\[A_{\Sigma} := \End_{\mathscr{E}}(\mathscr{O}(\mathscr{E})).\]
This is called the \emph{moduli algebra} of \(\Sigma\) in \cite{david1}.
\end{definition}

\begin{definition}[{\cite[143]{EtingofTensorCategories}}]
Let \(A\) be an algebra in \(\mathscr{E}\). 
A \emph{right module} over \(A\) in \(\mathscr{E}\) is an object \(M \in \mathscr{E}\) together with an \emph{action morphism} \(\operatorname{act}: M \otimes A \to M\) of \(\mathscr{E}\) such that certain commutative diagrams commute.
Let \(M\) and \(N\) be right modules over \(A\) in \(\mathscr{E}\). A \emph{module morphism} from \(M\) to \(N\) is a morphism \(\alpha \in \Hom_{\mathscr{E}}(M, N)\) which is compatible with the action.
The category of right modules over \(A\) in \(\mathscr{E}\) and module morphisms is denoted \(\Mod{A}{\mathscr{E}}\).
\end{definition}

\begin{proposition}{\cite[Theorem~5.14]{david1}}
\label{prop:factmodcat}
Let \(\Sigma\) be a punctured surface, and \(\mathscr{E}\) be an abelian \(\mathbb{C}\)-linear compact-rigid balanced tensor category\footnote{In Theorem 5.14, $\mathscr{E}$ is not required to be balanced. The reason we require it here is that we are working with an oriented version of factorisation homology.}.
\[
    \int_{\Sigma} \mathscr{E} \simeq \Mod{A_{\Sigma}}{\mathscr{E}},
\]
where \(A_{\Sigma}\) is the algebra object of the factorisation homology. 
\end{proposition}

\begin{remark}
Note that as the factorisation homology is equivalent to a category of modules over an algebra, it is an abelian category. 
\end{remark}

There is a combinatorial description of \(A_{\Sigma}\) in terms of the gluing pattern of the surface. 
\begin{definition}
A \emph{gluing pattern} is a bijection \[P: \{\,1, 1', \dots, n, n'\,\} \to \{\,1, 2, \dots, 2n-1, 2n\,\}\] such that \(P(i) < P(i')\) for all \(i = 1, \dots, n\). 

A gluing pattern \(P\) determines a marked surface \(\Sigma(P)\) by gluing together a disc and \(n\) handles \(H_i \cong [0,1]^2\) as follows: mark the disc with \(2n+1\) boundary intervals labelled \(1, \dots, 2n+1\); for each handle \(H_i\) mark two intervals \(i\) and \(i'\) on the boundary; glue the handles to the disc by identifying the interval \(i\) with the interval \(P(i)\) and the interval \(i'\) with the interval \(P(i')\) for all \(i = 1, \dots, n\). The final interval \(2n+1\) on the boundary of the disc gives \(\Sigma(P)\) a marking. 
\end{definition}

\begin{definition}
The handles \(H_i\) and \( H_j\), with \(i < j\) are:
\begin{enumerate}
\item \emph{positively linked} if \(P(i) < P(j) < P(i') < P(j')\),
\item \emph{positively nested} if \(P(i) < P(j) < P(j') < P(i')\),
\item \emph{positively unlinked} if \(P(i) < P(i') < P(j) < P(j')\).
\end{enumerate}
By relabelling the handles we can assume all handles are of the above forms.
\end{definition}

\begin{example}
The four-punctured sphere has the simplest possible gluing pattern with three handles
\begin{align*}
&P: \{\,1,1',2,2',3,3'\,\} \to \{\,1,2,3,4,5,6\,\}:\\
&P(1) = 1, P(1') =2, P(2) = 3, P(2') = 4, P(3) = 5, P(3') = 6.
\end{align*}
All three of its handles are positively unlinked.
\end{example}

\begin{figure}[ht!]
    \centering
    \myincludesvg{height=3cm}{}{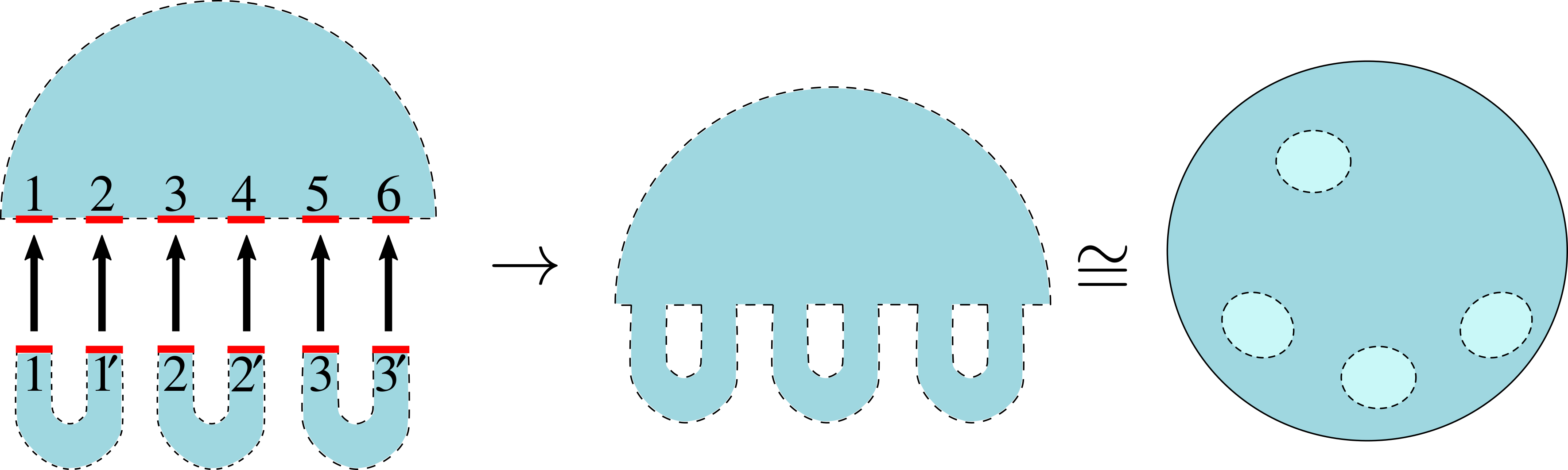}
    \caption{The gluing pattern of \(\Sigma_{0,4}\).}
\end{figure}
\begin{example}
The punctured torus has the gluing pattern 
\[P: \{\,1,1', 2, 2'\,\} \to \{\, 1, 2, 3, 4 \,\}: P(1) = 1, P(1') = 3, P(2) = 2, P(2') = 4.\]
The handles \(H_1\) and \(H_2\) are positively linked. 
\begin{figure}[ht!]
    \centering
    \myincludesvg{height=3cm}{}{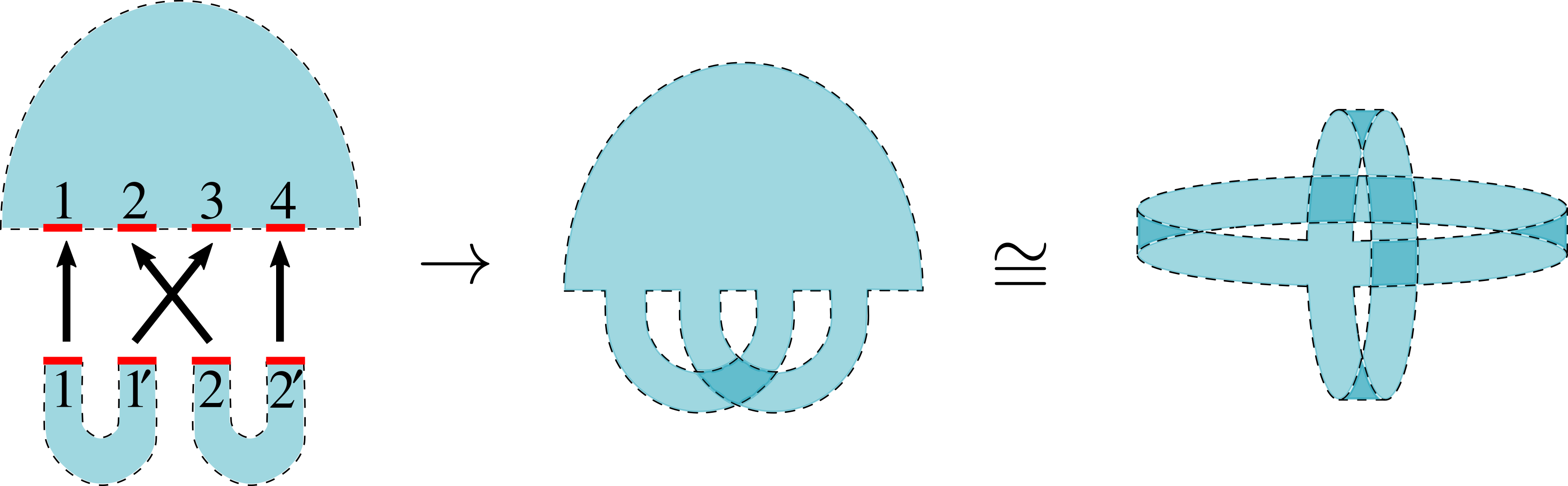}
    \caption{The gluing pattern of \(\Sigma_{1,1}\).}
\end{figure}
\end{example}

The moduli algebra \(A_{\Sigma}\) is constructed from copies of the distinguished object---one for each handle---with crossing morphisms determined from the type of handle crossing.

When \(\mathscr{E}\) is semisimple, the distinguished object 
\[
\mathscr{O}(\mathscr{E}) \cong \bigoplus_{X \text{ is simple}} X^* \otimes X.
\]
Using this we can defining the crossing morphisms.\footnote{When \(\mathscr{E}\) is not semisimple the crossing morphisms are still defined as we get that the distinguished object is a quotient of the direct sum over compact objects \cite[32]{david1}.}
\begin{definition}[{\cite[32, 36]{david1}}]
\label{defn:crossing}
Define the \emph{crossing morphism} 
\[K_{i,j}: \mathscr{O}(\mathscr{E})^{(i)} \otimes \mathscr{O}(\mathscr{E})^{(j)} \to \mathscr{O}(\mathscr{E})^{(j)} \otimes \mathscr{O}(\mathscr{E})^{(i)}\] 
as follows:
\[
    \includegraphics[height=4cm]{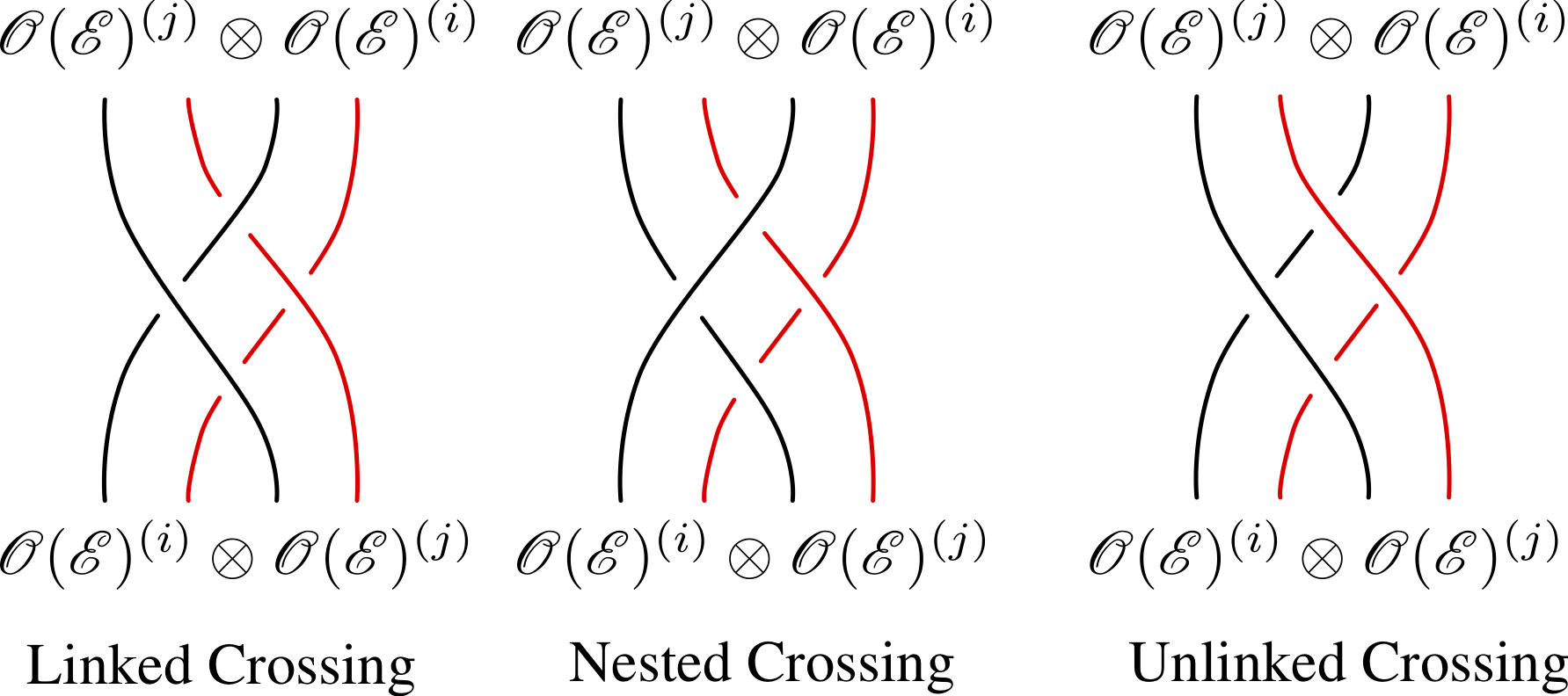}
\]
where strand crossings are determined by the braiding on \(\mathscr{E}\).
\end{definition}
As the crossing morphisms satisfy the Yang--Baxter equation, they can be used to extend the multiplication \(m: \mathscr{O}(\mathscr{E}) \otimes \mathscr{O}(\mathscr{E}) \to \mathscr{O}(\mathscr{E}) \) to a associative multiplication map \(m_n: \mathscr{O}(\mathscr{E})^{\otimes n} \otimes  \mathscr{O}(\mathscr{E})^{\otimes n} \to \mathscr{O}(\mathscr{E})^{\otimes n} \) turning \(\mathscr{O}(\mathscr{E})^{\otimes n}\) into an algebra~\cite[Theorem~3]{lebed2013}. 
\begin{figure}[ht!]
\centering
\includegraphics[height=4cm]{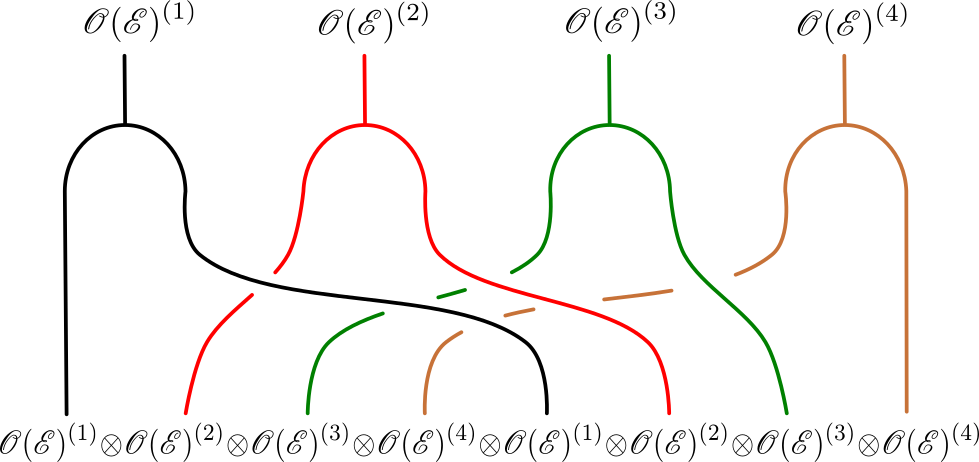}
\caption{The multiplication map for \(\mathscr{O}(\mathscr{E})^{\otimes 4}\) where the crossing of strands \(\mathscr{O}(\mathscr{E})^{(i)}\) and \(\mathscr{O}(\mathscr{E})^{(j)}\) is given by the braiding \(K_{i, j}\)}
\end{figure}

\begin{proposition}{{\cite[Theorem~5.14]{david1}}} \label{factofsurface}
Let \(\Sigma(P)\) be a surface determined by a gluing pattern \(P\) with \(n\) handles. Then \(A_{\Sigma(P)}\) is isomorphic to the algebra
\[a_P = \mathscr{O}(\mathscr{E})^{(1)} \otimes \dots \otimes \mathscr{O}(\mathscr{E})^{(n)},\]
where \(\mathscr{O}(\mathscr{E})^{(i)}\) is the distinguished object, and the crossing morphisms \(K_{i,j}: \mathscr{O}(\mathscr{E})^{(j)} \otimes \mathscr{O}(\mathscr{E})^{(i)} \to \mathscr{O}(\mathscr{E})^{(i)} \otimes \mathscr{O}(\mathscr{E})^{(j)}\) are defined in \cref{defn:crossing}.
\end{proposition}

\begin{remark}
When \(\mathscr{E} = \Rep_q(G)\) the algebra $A_{\Sigma(P)}$ is the moduli algebra of Alekseev \cite[Section~2]{Alek94}.
\end{remark}
\subsection{The Algebra of \(\qgroup{\mathfrak{g}}\)-Invariants and Character Varieties}
\label{section:InvAlgebra}
Given a surface \(\Sigma\) there are several invariants of \(\Sigma\) based on the representations of its fundamental group \(\pi_1(\Sigma)\). 

\begin{definition}
Let \(G\) be a reductive algebraic group. The \emph{representation variety} \(\mathfrak{R}_G(\Sigma)\) is the affine variety 
\[
\mathfrak{R}_G(\Sigma) = \left\{\,\rho: \pi_1(\Sigma) \to G \,\right\} 
\]
of homomorphisms from the fundamental group of \(\Sigma\) to \(G\).
\end{definition}

\begin{definition}
The \emph{character stack} \(\uChg(\Sigma)\) is the quotient \(\mathfrak{R}_G(\Sigma) / G\)  of the representation variety of the surface \(\mathfrak{R}_G(\Sigma)\) by the group \(G\) acting upon it by conjugation. 
\end{definition}

\begin{definition}
The \emph{character variety} \(\Chg_G(\Sigma)\) is the affine categorical quotient \(\mathfrak{R}_G(\Sigma) // G\) of the representation variety of the surface \(\mathfrak{R}_G(\Sigma)\) by the group \(G\) acting upon it by conjugation.  
\end{definition}

The character stack \(\uChg(\Sigma)\) is intimately related to the factorisation homology of \(\Sigma\) with coefficients in the category \(\Rep(G)\) of representations of \(G\):

\begin{theorem}{{\cite{BFN10}[Theorem~7.1]~\cite{david1}}}
If \(\Sigma\) is a surface, then there is an equivalence of categories 
\[\catname{QCoh}(\uChg(\Sigma)) \simeq \int_{\Sigma} \Rep (G)\]
between the category of quasi-coherent sheaves on the character stack \(\uChg(\Sigma)\) and the factorisation homology of the surface \(\Sigma\)  with coefficients in \(\Rep (G)\).
\end{theorem}

By replacing \(\Rep(G)\) with \(\Rep_q(G)\), one obtains a quantisation:

\begin{proposition}{{\cite[Section~7]{david1}}}
Let \(\Sigma\) be a punctured surface.
The factorisation homology \(\int_{\Sigma} \Rep_q(G)\) is a deformation quantisation of the category of sheaves on the character variety \(\Chg_G(\Sigma)\).
\end{proposition}

One can also use factorisation homology to quantise the character variety \(\Chg_G(\Sigma)\).
As the moduli algebra \(A_{\Sigma} \in \Rep_q(G)\), it is an \(\qgroup{\frg}\)-module. Hence, there is an action of the Hopf algebra \(\qgroup{\frg}\) on \(A_{\Sigma}\).

\begin{definition}
\label{defn:invalg}
We denote by \(\mathscr{A}_{\Sigma}\) the subalgebra of \(\qgroup{\mathfrak{g}}\)-invariants of the moduli algebra \(A_{\Sigma}\).
\end{definition}

\begin{proposition}[{\cite[Section 2]{Alek94}~\cite[Theorem~7.3]{david1}}]
Let \(\Sigma\) be a punctured surface. The algebra of \(\qgroup{\mathfrak{g}}\)-invariants \(\mathscr{A}_{\Sigma}\) of \(\int_{\Sigma} \Rep_q(G)\) is a quantisation of the character variety \(\Chg_G(\Sigma)\).
\end{proposition}

\begin{example} \label{Kaction}
In \cref{section:FactorisationHomology} we shall see that the algebra object \(A_{\Sigma_{0,4}}\) is generated by twelve generators 
\[\begin{pmatrix} x^1_1 & x^1_2 \\ x^2_1 & x^2_2 \end{pmatrix}\]
for \(x \in \{\,a, b, c\,\}\) and where \(x^i_j \in V^* \otimes V\). The quantum group \(\qgroup{\frsl}\) is generated by \(E, F, K^{\pm}\) whose images in the standard 2-dimensional representation are
\[E = 
\begin{pmatrix}
0 & 1 \\
0 & 0
\end{pmatrix}; \;
F = 
\begin{pmatrix}
0 & 0 \\
1 & 0
\end{pmatrix}; \;
K = 
\begin{pmatrix}
q & 0 \\
0 & q^{-1}
\end{pmatrix}.\]
It is a Hopf algebra with coproduct \(\Delta\) defined by 
\[\Delta(E) = E \otimes 1 + K^{-1} \otimes E,\; \Delta(F) = F \otimes K + 1 \otimes F,\; \Delta(K) = K \otimes K; \]
antipode \(S\) defined by
\[S(E) = KE, \; S(F) = -FK^{-1}, \; S(K) = K^{-1};\]
 and counit \(\epsilon\) defined by \(\epsilon(E) = \epsilon(F) = 0, \epsilon(K) = 1\). The vector space \(V\) with basis \(\{\,v_1, v_2\,\}\) has an \(\qgroup{\frsl}\) action on it defined by 
\begin{alignat*}{3}
    &K \cdot v_1 = q v_1; \quad &&K \cdot v_2 = q^{-1} v_2;& \\
    &E \cdot v_1 = 0; &&E \cdot v_2 = v_1;&\\
    &F \cdot v_1 = v_2; &&F \cdot v_2 = 0.&
\end{alignat*}
The action on the dual \(V^*\) is defined by \(X \cdot u^* (w) = u^* (S(X)w)\) where \(X \in \qgroup{\frsl}, u^* \in V^*, w \in V\), so on the basis \(\{\,v^1, v^2\,\}\) is given by 
\begin{alignat*}{3}
    &K \cdot v^1 = q v^1; &&K \cdot v^2 = q^{-1} v^2; & \\
    &F \cdot v^1 = -q^{-1} v^2; \quad &&F \cdot v^2 = 0; &\\
    &E \cdot v^1 = 0; &&E \cdot v^2 = -q v^1 &
\end{alignat*}
The action of \(\qgroup{\frsl}\) on \(V^* \otimes V\) is defined via the coproduct; hence, it acts on \(A_{\Sigma_{0,4}}\) as follows:
\begin{alignat*}{5}
K \cdot a^1_1 &= a^1_1; \quad &K \cdot a^1_2 &= q^2 a^1_2; \quad &K \cdot a^2_1 &= q^{-2} a^2_1; \quad &K \cdot a^2_2 &= a^2_2; \\
E \cdot a^1_1 &= q^{-1} a^1_2; \quad &E \cdot a^1_2 &= 0; \quad &E \cdot a^2_1 &= q (a^2_2 - a^1_1);\quad &E \cdot a^2_2 &= -q a^1_2; \\
F \cdot a^1_1 &= -q^{-2} a^2_1; \quad &F \cdot a^1_2 &= a^1_1 - a^2_2; &\quad F \cdot a^2_1 &=0; \quad &F \cdot a^2_2 &= a^2_1.
\end{alignat*}

An element \(x \in A_{\Sigma_{0,4}}\) is an invariant of the \(\qgroup{\frsl}\)-action if \(h \cdot v = \epsilon(h) v\) i.e.\ \(E \cdot v = F \cdot v = 0\) and \(K \cdot v = v\). 
So, the algebra of invariants quantisation of the \(SL_2\)-quantum character variety of \(\Sigma_{0,4}\) is given by the elements of \(A_{\Sigma_{0,4}}\) which are invariant under this action. We shall give a presentation for \(\mathscr{A}_{\Sigma_{0, 4}}\) in Section \cref{section:InvAlgebraSphere}.
\end{example}
\subsection{Reduction Systems and the Diamond Lemma}
\label{section:PBW}
Both the universal enveloping algebra of a Lie algebra \(\mathcal{U}(\frg)\) and its quantum group \(\qgroup{\frg}\) have a Poincare--Birkhoff--Witt basis (PBW-basis). In the case of \(\mathcal{U}(\frg)\) this means that if \(x_1, \dots, x_l\) is an ordered basis of \(\frg\) then \(\mathcal{U}(\frg)\) has a vector space basis given by the monomials 
\[y_1^{k_1} y_2^{k_2} \dots y_l^{k_l}\]
where \(k_i \in \mathbb{N}_0\) and \(x_i \mapsto y_i\) via the map \(\frg \to \mathcal{U}(\frg)\). In the case of \(\qgroup{\frg}\) this means that \(\qgroup{\frg}\) has a vector space basis given by the monomials 
\[ (X^+_1)^{a_1} \dots (X^+_n)^{a_n} K_1^{b_1} \dots K_n^{b_n} (X^-_1)^{c_1} \dots (X^-_n)^{c_n} \]
where \(a_i, c_i \geq 0\) and \(b_i \in \mathbb{Z}\).

In this section we recall the definitions and results needed to define and prove the existence of such bases. We will use these results in \cref{section:FactofSphere} and \cref{section:InvAlgebraSphere} to provide PBW-bases for the algebra objects and \(\qgroup{\mathfrak{sl}_2}\)-invariant algebras of the factorisation homology of the four-punctured sphere and punctured torus with coefficients in \(\Rep_q(\SL)\). The definitions given in this section can be found \cite[Section~1]{diamond} except those relating to the reduced degree which can be found in \cite[Section~15]{ReducedOrder}, and the main result is the diamond lemma for rings proven by Bergman \cite[Theorem~1.2]{diamond}. Let \(k\) be a commutative ring with multiplicative identity and \(X\) be an alphabet (a set of symbols from which we form words).  

\begin{definition}
A \emph{reduction system} \(S\) consists of term rewriting rules \(\sigma : W_{\sigma} \mapsto f_{\sigma}\) where \(W_{\sigma} \in \langle X \rangle\) is a word in the alphabet \(X\) and \(f_{\sigma} \in k \langle X \rangle\) is a linear combination of words. A \emph{\(\sigma\)-reduction} \(r_{\sigma}(T)\) of an expression \(T \in k\langle X \rangle\) is formed by replacing an instance of \(W_{\sigma}\) in \(T\) with \(f_{\sigma}\). For example, if \(X = \langle a,b \rangle\) and \(S = \{\,\sigma : ab \mapsto ba\,\}\) then \(r_{\sigma}(T) = aba + a\) is a \(\sigma\)-reduction of \(T= aab + a\). A \emph{reduction} is a \(\sigma\)-reduction for some \(\sigma \in S\). 
\end{definition}

\begin{definition}
The five-tuple \((\sigma, \tau, A, B, C)\) with \(\sigma, \tau \in S\) and \(A,B, C \in \langle X \rangle\) is an \emph{overlap ambiguity} if \(W_{\sigma} = AB\) and \(W_{\tau} = BC\) and an \emph{inclusion ambiguity} if \(W_{\sigma} = B\) and \(W_{\tau} = ABC\).
These ambiguities are \emph{resolvable} if reducing \(ABC\) by starting with a \(\sigma\)-reduction gives the same result as starting with a \(\tau\)-reduction. 
For example if \(S=\{\,\sigma : ab \mapsto ba, \tau : ba \mapsto a\,\}\) then \((\sigma, \tau, a,b,a)\) is an overlap ambiguity which is resolvable as \( aba \xmapsto{r_{\sigma}} ba^2 \xmapsto{r_{\tau}} a^2 \) gives the same expression as \(aba \xmapsto{r_{\tau}} a^2\).
\end{definition} 

\begin{definition}
A \emph{semigroup partial ordering} \(\leq\) on \(\langle X \rangle\) is a partial order such that \(B \leq B' \) implies that \(ABC \leq AB'C\) for all words \(A,B, B', C\); it is \emph{compatible with the reduction system} \(S\) if  for all \(\sigma \in S\) the monomials in \(f_{\sigma}\) are less than or equal to \(W_{\sigma}\).
\end{definition}

\begin{definition}
A reduction system \(S\) satisfies the \emph{descending chain condition} or is \emph{terminating} if for any expression \(T \in k \langle X \rangle\) any sequence of reductions terminates in a finite number of reductions with an irreducible expression. 
\end{definition}

\begin{lemma}[The Diamond Lemma {\cite[Theorem~1.2]{diamond}}]
Let \(S\) be a reduction system for \(k \langle X \rangle\) and let \(\leq\) be a semigroup partial ordering on \(\langle X \rangle\) compatible with the reduction system \(S\) with the descending chain condition. The following are equivalent:
\begin{enumerate}
    \item All ambiguities in \(S\) are resolvable (\(S\) is \emph{locally confluent});
    \item Every element \(a \in k \langle X \rangle\) can be reduced in a finite number of reductions to a unique expression \(r_S(a)\) (\(S\) is \emph{confluent});
    \item The algebra \(R = k \langle X \rangle / I\), where \(I\) is the two sided ideal of \(k \langle X \rangle\) generated by the elements \((W_{\sigma} - f_{\sigma})\), can be identified with the \(k\)-algebra \(k \langle X \rangle_{\mathrm{irr}}\) spanned by the \(S\)-irreducible monomials of \(\langle X \rangle\) with multiplication given by \(a \cdot b = r_S(ab)\). These \(S\)-irreducible monomials are called a Poincare--Birkhoff--Witt--basis of \(R\).
\end{enumerate}
\end{lemma}

\begin{remark}
Bergman's diamond lemma is an application to ring theory of the diamond lemma for abstract rewriting systems. An \emph{abstract rewriting system} is a set A together with a binary relation \(\to\) on A called the \emph{reduction relation} or \emph{rewrite relation}. 
\begin{enumerate}
    \item It is \emph{terminating} if there are no infinite chains \(a_0 \to a_1 \to a_2 \to \dots\). 
    \item It is \emph{locally confluent} if for all \(y \xleftarrow{} x \xrightarrow{} z\) there exists an element \(y \downarrow z \in A\) such that there are paths \(y \to \dots \to (y \downarrow z) \) and \(z \to \dots \to (y \downarrow z) \). 
    \item It is \emph{confluent} if for all \(y \xleftarrow{} \dots \xleftarrow{} x \xrightarrow{} \dots \xrightarrow{} z\) there exists an element \(y \downarrow z \in A\) such that there are paths \(y \to \dots \to (y \downarrow z) \) and \(z \to \dots \to (y \downarrow z) \). In a terminating confluent abstract rewriting system an element \(a \in A\) will always reduce to a unique reduced expression regardless of the order of the reductions used.
\end{enumerate}
The diamond lemma (or Newman's lemma) for abstract rewriting systems states that a terminating abstract rewriting system is confluent if and only if it is locally confluent.

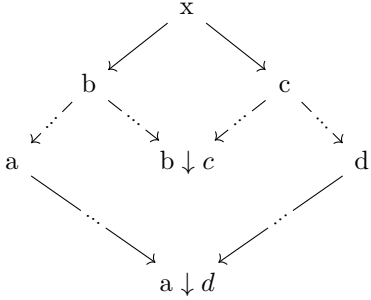
\begin{figure}[H]
\floatbox[{\capbeside\thisfloatsetup{capbesideposition={right,top},capbesidewidth=8cm}}]{figure}[\FBwidth]
{\caption{If the abstract term rewriting system is locally confluent there exists \(b \downarrow c \in A\) forming a small diamond shape. If it is confluent there exists \(a \downarrow d \in A\) forming a larger diamond shape. The diamond lemma is proven by patching together the small diamonds to give the larger diamonds and inducting on path length, hence the name.}}
{\begin{tikzcd}[row sep=1.5em, column sep=1.5em, ampersand replacement=\&]
\& \& 
\eqmakebox[a][c]{x} \ar[ld] \ar[rd] 
\& \& \\
\& \eqmakebox[a][c]{b}  \ar[rd, "\dots" description, sloped] \ar[ld, "\dots" description, sloped] 
\& \& \eqmakebox[a][c]{\vphantom{d}c} \ar[ld, "\dots" description, sloped] \ar[rd, "\dots" description, sloped] \& \\
\eqmakebox[a][c]{\vphantom{d}a} \ar[rrdd, "\dots" description, sloped]
\& \& \eqmakebox[a][c]{b} \downarrow c 
\& \& \eqmakebox[a][c]{d} \ar[lldd, "\dots" description, sloped]\\
 \& \& \& \& \\
 \& \& \eqmakebox[a][c]{\vphantom{d}a} \downarrow d \& \& 
\end{tikzcd}}
\end{figure}
\end{remark}

In this paper the semigroup partial ordering we shall use is ordering by \emph{reduced degree}:

\begin{definition}
Give the letters of the finite alphabet \(X\) an ordering \(x_1 \leq \dots \leq x_N\). Any word \(W\) of length \(n\) can be written as \(W = x_{i_1} \dots x_{i_n}\) where \(x_{i_j} \in X\). An \emph{inversion} of \(W\) is a pair \(k \leq l\) with \(x_{i_{k}} \geq x_{i_{l}}\) i.e.\ a pair with letters in the incorrect order. The number of inversions of \(W\) is denoted \(|W|\).
\end{definition}
\begin{definition}
Any expression \(T\) can be written as a linear combination of words \(T= \sum c_l W_l\). Define \(\rho_n(T):= \sum_{\operatorname{length}(W_l)=n, c_l \neq 0} |W_l|\). The \emph{reduced degree of \(T\)} is the largest \(n\) such that \(\rho_n(T) \neq 0\).
\end{definition}
\begin{definition}
Under the \emph{reduced degree ordering}, \(T \leq S\) if 
\begin{enumerate}
\item The reduced degree of  \(T\) is less than the reduced degree of \(S\), or
\item The reduced degree of \(T\) and \(S\) are equal, but  \(\rho_n(T) \leq \rho_n(S)\) for maximal nonzero \(n\). 
\end{enumerate}
\end{definition}

\section{The Algebra of \(\qgroup{\mathfrak{sl}_2}\)-Invariants of the Four-Punctured Sphere and Punctured Torus}
In this section we shall find an explicit description for the algebra of \(\qgroup{\mathfrak{sl}_2}\)-invariants \(\mathscr{A}_{\Sigma}\) of the factorisation homology \(\int_{\Sigma}^{\LFP{k}} \Rep_q(\SL)\) when \(\Sigma\) is the four-punctured surface \(\Sigma_{0,4}\) or the punctured torus \(\Sigma_{1,1}\). Throughout this section we shall always assume \(\Sigma\) is a punctured surface.

\subsection{The Factorisation Homology of the \\ Four-Punctured Sphere and Punctured Torus over \(\qgroup{\frsl}\)}
\label{section:FactofSphere}

Before considering the algebra of \(\qgroup{\mathfrak{sl}_2}\)-invariants \(\mathscr{A}_{\Sigma}\) of \(\int_{\Sigma}^{\LFP{k}} \Rep_q(\SL)\), we must first consider the factorisation homology \(\int_{\Sigma}^{\LFP{k}} \Rep_q(\SL)\). Using \cref{prop:factmodcat} we have that the factorisation homology of the four-punctured sphere and punctured torus over \(\qgroup{\frsl}\) is \(\Mod{A_{\Sigma}}{\Rep_q(\SL)}\) where \(A_{\Sigma}\) is the algebra object of the four-punctured sphere \(\Sigma_{0, 4}\) and punctured torus \(\Sigma_{1, 1}\) respectively. We shall use \cref{factofsurface} to obtain presentations of  \(A_{\Sigma_{0,4}}\) and \(A_{\Sigma_{1,1}}\). In order to do this, we need a presentation of the distinguished object \(\mathscr{O}(\Rep_q(\SL))\) and a description of \(K_{i, j}\)  in each case. Both of these depend on the choice of \(R\)-matrix for \(\Rep_q(\SL)\): we shall use the standard \(R\)-matrix. The \(R\)-matrix for \(\qgroup{\frsl}\) when evaluating on the standard representation of \(\qgroup{\frsl}\) is given by
\[R := \begin{pmatrix} R^{11}_{11} & R^{12}_{11} & R^{21}_{11} & R^{22}_{11} \\
R^{11}_{12} & R^{12}_{12} & R^{21}_{12} & R^{22}_{12} \\
R^{11}_{21} & R^{12}_{21} & R^{21}_{21} & R^{22}_{21} \\
R^{11}_{22} & R^{12}_{22} & R^{21}_{22} & R^{22}_{22} 
\end{pmatrix} := q^{\frac{1}{2}}\begin{pmatrix} q & 0& 0& 0 \\ 0& 1 & (q-q^{-1}) & 0\\ 0& 0& 1 &0 \\0 & 0& 0& q \end{pmatrix}.\]
We shall also require 
\[\tilde{R} := (\Id \otimes S)(R) = q^{-\frac{1}{2}}\begin{pmatrix} q^{-1} & 0 & 0 & 0 \\0 & 1 & q^{-2}(q^{-1}-q) & 0\\0 &0 & 1 & 0\\0 &0 &0 & q^{-1} \end{pmatrix}\]
where \(S\) is the antipode of \(\qgroup{\frsl}\).

The distinguished object \(\mathscr{O}(\Rep_q(\SL))\) is the reflection equation algebra of \(\qgroup{\SL}\) \cite[Section~6]{david1}:
\begin{definition}{\cite[Definitin~3.3]{balagovic}}
The \emph{reflection equation algebra}\footnote{This algebra also goes by other names such as the `equivariantised quantum coordinate algebra' and the `quantum
loop algebra'.} \(\mathscr{O}_q(\SL)\) is generated by the four elements \[A = \begin{pmatrix}a^1_1& a^1_2\\ a^2_1 & a^2_2\end{pmatrix}\] which satisfy the following:
\begin{enumerate}
    \item The \emph{quantum determinant} \(\det_q(A) := a^1_1 a^2_2 - q^{2}a^1_2 a^2_1 = 1\), and
    \item The \emph{reflection equation}  \(a^l_m a^p_r = \tilde{R}^{op}_{mk} (R^{-1})^{kl}_{ij} R^{sj}_{uv} R^{wu}_{or} a^i_s a^v_w\) where \(i\), \(j\), \(k\), \(l\), \(m\), \(o\), \(p\), \(r\), \(s\), \(v\), \(w\) \(\in \{\,0,1\,\}\)\footnote{The reflection equation algebra is usually given as \(R_{21} A_1 R A_2 = A_2 R_{21} A_1 R\) where \(A_1 := A \otimes I\), \(A_2:= I \otimes A\), and \(R_{21}:= \tau R \tau\), for example in \cite{REA} and \cite{REA2}. Our version is the tensor version rearranged using the relations \(\sum (R^{-1})^{ij}_{kl} R^{kl}_{mn} = \delta^i_m \delta^j_n\) and \(\sum \tilde{R}^{ij}_{kl} R^{ml}_{in} = \delta^m_k \delta^n_j\).}.
\end{enumerate}
Or more explicitly the reflection equation algebra \(\mathscr{O}_q(\SL)\) has generators \(a^1_1, a^1_2, a^2_1, a^2_2\) and relations
\begin{align}
    a^1_2 a^1_1 &= a^1_1 a^1_2 + \left(1-q^{-2}\right)a^1_2 a^2_2, \\
    a^2_1 a^1_1 &= a^1_1 a^2_1 - q^{-2}\left(1-q^{-2}\right)a^2_1 a^2_2, \\
    a^2_1 a^1_2 &= a^1_2 a^2_1 + \left(1-q^{-2}\right)\left(a^1_1 a^2_2 - a^2_2 a^2_2\right), \\
    a^2_2 a^1_1 &= a^1_1 a^2_2, \\
    a^2_2 a^1_2 &= q^2 a^1_2 a^2_2, \\
    a^2_2 a^2_1 &= q^{-2} a^2_1 a^2_2, \\
    a^1_1 a^2_2 &= 1 + q^2 a^1_2 a^2_1.
\end{align}
\end{definition}

The reflection equation algebra \(\mathscr{O}_q(\SL)\) is generated by elements of the form \(a^i_j = v^i \otimes v_j \in V^* \otimes V\) where \(V\) is the standard representation. The braiding \(\sigma\) is defined on the copies \(V,W\) of the standard representation and their duals as follows:
\begin{align*}
&\sigma_{V,W}(w\otimes v) = \tau_{V, W} \circ R(w \otimes v); \\
&\sigma_{V^*, W}(w^* \otimes v)= \tau_{V^*, W} \circ (S \otimes \Id) \circ R (w^* \otimes v) =  \tau_{V^*, W} \circ R^{-1} (w^* \otimes v); \\
&\sigma_{V, W^*}(w \otimes v^*) = \tau_{V, W^*} \circ (\Id \otimes S) \circ R (w \otimes v^*); \\
&\sigma_{V^*, W^*}(w^* \otimes v^*) = \tau_{V^*, W^*} \circ (S \otimes S)
 \circ R (w^* \otimes v^*) = \tau_{V^*, W^*} \circ R (w^* \otimes v^*).
 \end{align*}

\begin{definition}
The \emph{braiding on \(\mathscr{O}_q(\SL)\) for positively unlinked handles} \(H_i\)  and \(H_j\)  is the map 
\begin{align*}
&K_{i,j}: \mathscr{O}_q(\SL)^{(i)} \otimes \mathscr{O}_q(\SL)^{(j)} \to \mathscr{O}_q(\SL)^{(j)} \otimes \mathscr{O}_q(\SL)^{(i)}: \\ 
&K_{i,j}(y^e_f \otimes x^g_h) = \tilde{R}^{ig}_{fj} R^{ej}_{kl} R^{mn}_{ih} \left(R^{-1}\right)^{ko}_{pn} x^l_o \otimes y^p_m 
\end{align*}
where \(x^g_h\) and \(y^e_f\) are generators of \(\mathscr{O}^{(i)}_q(\SL)\) and \(\mathscr{O}^{(j)}_q(\SL)\) respectively. 
\end{definition}

So applying \cref{factofsurface} we obtain 

\begin{corollary}
\label{prop:facthom}
The factorisation homology of the four-punctured sphere with coefficients in \(\Rep_q (\SL)\) is \(\int_{\Sigma_{0,4}} \Rep_q(\SL) \simeq \Mod{A_{\Sigma_{0,4}}}{\Rep_q(\SL)}\) where \(A_{\Sigma_{0,4}}\) is an algebra with twelve generators organised into three matrices 
\[
A := \begin{pmatrix}
a^1_1 & a^1_2 \\ a^2_1 & a^2_2
\end{pmatrix},\;
B := \begin{pmatrix}
b^1_1 & b^1_2 \\ b^2_1 & b^2_2
\end{pmatrix},\;
C:= \begin{pmatrix} 
c^1_1 & c^1_2 \\ c^2_1 & c^2_2 
\end{pmatrix}
\]
subject to the relations 
\begin{align}
    x^1_1 x^2_2 &= 1 + q^2 x^1_2 x^2_1 &\text{ (determinant relation)}\\
    x^l_m x^p_r &= \tilde{R}^{op}_{mk} (R^{-1})^{kl}_{ij} R^{sj}_{uv} R^{wu}_{or} x^i_s x^v_w &\text{ (reflection equation)} \\
    y^e_f  x^g_h &= \tilde{R}^{ig}_{fj} R^{ej}_{kl} R^{mn}_{ih} (R^{-1})^{ko}_{pn} x^l_o y^p_m &\text{ (crossing relation)}
\end{align}
where \(x < y \in \{\,a,b,c\,\}\) (using ordering \(a < b < c\)), \(e,f,g,h,i,j,k,l,m,n,o,p \in \{\,0,1\,\}\),
\[
R= q^{\frac{1}{2}}\begin{pmatrix} q & 0& 0& 0 \\ 0& 1 & (q-q^{-1}) & 0\\ 0& 0& 1 & 0\\ 0& 0& 0& q \end{pmatrix}
\]
is the standard quantum \(R\)-matrix for \(\qgroup{\frsl}\) when evaluated on the standard representation of \(\qgroup{\frsl}\) and 
\[
\tilde{R} = q^{- \frac{1}{2}}\begin{pmatrix} q^{-1} & 0& 0& 0 \\ 0& 1 & q^{-2}(q^{-1}-q) & 0\\ 0& 0& 1 & 0\\ 0& 0& 0& q^{-1} \end{pmatrix}.
\]
\end{corollary}

\begin{definition}
The \emph{braiding on \(\mathscr{O}_q(\SL)\) for positively linked handles} \(H_i\)  and \(H_j\)  is the map 
\begin{align*}
&K_{i,j}: \mathscr{O}^{(i)}_q(\SL) \otimes \mathscr{O}^{(j)}_q(\SL) \to \mathscr{O}^{(j)}_q(\SL) \otimes \mathscr{O}^{(i)}_q(\SL): \\ 
&K_{i,j}(y^g_h \otimes x^e_f) = 
 \tilde{R}^{i e}_{h j} R^{g j}_{k l} R^{m n}_{i f} \left(R^{-1}\right)^{k o}_{p n} x^l_ o \otimes y^p_m
\end{align*}
where \(x^g_h\) and \(y^e_f\) are generators of \(\mathscr{O}^{(i)}_q(\SL)\) and \(\mathscr{O}^{(j)}_q(\SL)\) respectively. 
\end{definition}

So applying \cref{factofsurface} we obtain 

\begin{corollary}
\label{prop:facthomtorus}
The factorisation homology of the punctured torus with coefficients in \(\qgroup{\frsl}\) is \(\int_{\Sigma_{1,1}} \Rep_q(\SL) \simeq \Mod{A_{\Sigma_{1,1}}}{\Rep_q(\SL)}\) where \(A_{\Sigma_{1,1}}\) is an algebra with eight generators organised into two matrices 
\[
A := \begin{pmatrix}
a^1_1 & a^1_2 \\ a^2_1 & a^2_2
\end{pmatrix},\;
B := \begin{pmatrix}
b^1_1 & b^1_2 \\ b^2_1 & b^2_2
\end{pmatrix}
\]
subject to the relations 
\begin{align}
    x^1_1 x^2_2 &= 1 + q^2 x^1_2 x^2_1 &\text{ (determinant relation)}\\
    x^l_m x^p_r &= \tilde{R}^{op}_{mk} (R^{-1})^{kl}_{ij} R^{sj}_{uv} R^{wu}_{or} x^i_s y^v_w &\text{ (reflection equation)} \\
    y^g_h  x^e_f &= \tilde{R}^{i e}_{h j} R^{g j}_{k l} R^{m n}_{i f} \left(R^{-1}\right)^{k o}_{p n} x^l_ o \otimes y^p_m &\text{ (crossing relation)}
\end{align}
where \(x < y \in \{\,a,b,c\,\}\), \(e,f,g,h,i,j,k,l,m,n,o,p \in \{\,0,1\,\}\) and the R-matrices are the same as in \cref{prop:facthom}.
\end{corollary}

\subsection{PBW bases for the Algebra Objects}
We now construct a PBW basis for \(\mathscr{O}_q(\SL)\) which we shall use to construct PBW bases for \(A_{\Sigma_{0,4}}\) and \(A_{\Sigma_{1,1}}\).

\begin{proposition}
\label{pbwbasisorignial}
The set of monomials 
\[\left\{\,(a^1_1)^{\alpha} (a^1_2)^{\beta} (a^2_1)^{\gamma} (a^2_2)^{\delta} \; \middle| \; \alpha, \beta, \gamma, \delta \in \mathbb{N}_0,\; \beta \text{ or }\gamma  = 0 \,\right\}\] 
is a PBW basis for the reflection equation algebra \(\mathscr{O}_q(\SL)\) with respect to the ordering \(a^1_1 < a^1_2 < a^2_1 < a^2_2\).
\begin{proof}
The relations defining \(\mathscr{O}_q(\SL)\) can be re-expressed as the term rewriting system:
\begin{align*}
    \sigma_{1211} : a^1_2 a^1_1 &\mapsto a^1_1 a^1_2 + \left(1-q^{-2}\right)a^1_2 a^2_2,\\
    \sigma_{2111} : a^2_1 a^1_1 &\mapsto a^1_1 a^2_1 - q^{-2}\left(1-q^{-2}\right)a^2_1 a^2_2, \\
    \sigma_{2112} : a^2_1 a^1_2 &\mapsto a^1_2 a^2_1 + \left(1-q^{-2}\right)\left(a^1_1 a^2_2 - a^2_2 a^2_2\right),\\
    \sigma_{2211} : a^2_2 a^1_1 &\mapsto a^1_1 a^2_2,\\
    \sigma_{2212} : a^2_2 a^1_2 &\mapsto q^2 a^1_2 a^2_2,\\
    \sigma_{2221} : a^2_2 a^2_1 &\mapsto q^{-2} a^2_1 a^2_2,\\
    \sigma_{1221} : a^1_2 a^2_1 &\mapsto q^{-2} + q^{-2} a^1_1 a^2_2.
\end{align*}
The monomials listed in the statement of the result are the reduced monomials with respect to this term rewriting system; furthermore, there are no inclusion ambiguities, and the overlap ambiguities are 
\begin{align*}
    (\sigma_{2112}, \sigma_{1211}, a^2_1, a^1_2, a^1_1), & \quad
    (\sigma_{2212}, \sigma_{1211}, a^2_2, a^1_2, a^1_1), \\
    (\sigma_{2221}, \sigma_{2111}, a^2_2, a^2_1, a^1_1), & \quad
    (\sigma_{2221}, \sigma_{2112}, a^2_2, a^2_1, a^1_2), \\
    (\sigma_{2112}, \sigma_{1221}, a^2_1, a^1_2, a^2_1), & \quad
    (\sigma_{2212}, \sigma_{1221}, a^2_2, a^1_2, a^2_1), \\
    (\sigma_{1221}, \sigma_{2111}, a^1_2, a^2_1, a^1_1), & \quad
    (\sigma_{1221}, \sigma_{2112}, a^1_2, a^2_1, a^1_2).
\end{align*}

We shall order \(\mathscr{O}_q(\SL)\) with respect to the reduced degree where we give the generators the ordering \(a^1_1 < a^1_2 < a^2_1 < a^2_2\). This ordering is compatible with the given term rewriting systems and the rewriting will terminate, so if the ambiguities are resolvable then we can apply the diamond lemma, and we are done. It can be checked by direct calculation that the ambiguities are resolvable\footnote{We used the computer algebra system MAGMA to check this and similar computations throughout this paper.}. For example for the first ambiguity we have that both
\begin{align*}
    \left(a^2_1 a^1_2\right) a^1_1 &\mystackrel{$(\sigma_{2112})$}{=} a^1_2 \left(a^2_1 a^1_1\right) +\left(1-q^{-2} \right) \left(a^1_1 a^2_2 a^1_1 - \left(a^2_2 \right)^2 a^1_1 \right) \\
        &\mystackrel{$(\sigma_{2111},\sigma_{2211})$}{=} \left(a^1_2 a^1_1\right) a^2_1 - q^{-2}\left(1-q^{-2}\right)a^1_2 a^2_1 a^2_2 \\
            &\qquad+\left(1-q^{-2}\right)\left(\left(a^1_1\right)^2 a^2_2 - a^1_1 \left(a^2_2\right)^2\right) \\
        &\mystackrel{$(\sigma_{1211})$}{=} a^1_1 a^1_2 a^2_1 +\left(1-q^{-2}\right) a^1_2 \left(a^2_2 a^2_1\right) - q^{-2}\left(1-q^{-2}\right)a^1_2 a^2_1 a^2_2 \\
            &\qquad+ \left(1-q^{-2}\right)\left(\left(a^1_1\right)^2 a^2_2 - a^1_1 \left(a^2_2\right)^2\right) \\
        &\mystackrel{$(\sigma_{2221})$}{=} a^1_1 a^1_2 a^2_1 + q^{-2}\left(1-q^{-2}\right)a^1_2 a^2_1 a^2_2 - q^{-2}\left(1-q^{-2}\right)a^1_2 a^2_1 a^2_2 \\
            &\qquad+\left(1-q^{-2}\right) \left(\left(a^1_1\right)^2 a^2_2 - a^1_1 \left(a^2_2\right)^2\right) \\
        &\mystackrel{}{=} a^1_1 a^1_2 a^2_1 + \left(1 -  q^{-2}\right)\left(\left(a^1_1\right)^2 a^2_2 - a^1_1 \left(a^2_2\right)^2\right)
\end{align*}
and 
\begin{align*}
    a^2_1\left(a^1_2 a^1_1\right) &\mystackrel{$(\sigma_{1211})$}{=}\left(a^2_1 a^1_1\right)a^1_2 +\left(1-q^{-2}\right)a^2_1 a^1_2 a^2_2 \\
        &\mystackrel{$(\sigma_{2111})$}{=} a^1_1 a^2_1 a^1_2 - q^{-2}\left(1-q^{-2}\right)a^2_1\left(a^2_2 a^1_2\right)+\left(1-q^{-2} \right)a^2_1 a^1_2 a^2_2 \\
        &\mystackrel{$(\sigma_{2212})$}{=} a^1_1 a^2_1 a^1_2 -\left(1-q^{-2}\right) a^2_1 a^1_2 a^2_2 +\left(1-q^{-2} \right) a^2_1 a^1_2 a^2_2 \\
        &\mystackrel{}{=} a^1_1\left(a^2_1 a^1_2\right) \\
        &\mystackrel{($\sigma_{2112})$}{=} a^1_1 a^1_2 a^2_1 +\left(1-q^{-2}\right)\left(\left(a^1_1\right)^2 a^2_2 - a^1_1\left(a^2_2\right)^2\right)
\end{align*}
give the same result, so the first ambiguity is resolvable. 
\end{proof}
\end{proposition}

\begin{proposition}
\label{prop:PBWbasis}
A PBW basis for \(A_{\Sigma_{0,4}}\) is 
\begin{multline*}
\left\{ \,(a^1_1)^{\alpha_1} (a^1_2)^{\beta_1} (a^2_1)^{\gamma_1} (a^2_2)^{\delta_1} (b^1_1)^{\alpha_2} (b^1_2)^{\beta_2} (b^2_1)^{\gamma_2} (b^2_2)^{\delta_2} (c^1_1)^{\alpha_3} (c^1_2)^{\beta_3} (c^2_1)^{\gamma_3} (c^2_2)^{\delta_3} \right| \\ 
\left| \; \alpha_i, \beta_i, \gamma_i \in \mathbb{N}_0, \beta_i \text{ or } \gamma_i = 0\,\right\}. 
\end{multline*}
\end{proposition}

\begin{proof}
By \cref{pbwbasisorignial} we have a PBW basis \[\left\{\,(a^1_1)^{\alpha} (a^1_2)^{\beta} (a^2_1)^{\gamma} (a^2_2)^{\delta} \; \middle| \; \alpha, \beta, \gamma, \delta \in \mathbb{N}_0,\; \beta \text{ or }\gamma  = 0\,\right\}\]
for the reflection equation algebra \(\mathscr{O}_q(\SL)\). The algebra \(A_{\Sigma_{0,4}}\) is the tensor product of three copies of \(\mathscr{O}_q(\SL)\); hence, \begin{multline*} \left\{ \,(a^1_1)^{\alpha_1} (a^1_2)^{\beta_1} (a^2_1)^{\gamma_1} (a^2_2)^{\delta_1} (b^1_1)^{\alpha_2} (b^1_2)^{\beta_2} (b^2_1)^{\gamma_2} (b^2_2)^{\delta_2} (c^1_1)^{\alpha_3} (c^1_2)^{\beta_3} (c^2_1)^{\gamma_3} (c^2_2)^{\delta_3} \right| \\ 
\left| \; \alpha_i, \beta_i, \gamma_i \in \mathbb{N}_0,\; \beta_i \text{ or } \gamma_i = 0 \,\right\}. 
\end{multline*} is a PBW basis for it.
\end{proof}

\begin{proposition}
\label{prop:pbwbasistorus}
A PBW basis for \(A_{\Sigma_{1, 1}}\) is 
\begin{multline*}\left\{\,(a^1_1)^{\alpha_1} (a^1_2)^{\beta_1} (a^2_1)^{\gamma_1} (a^2_2)^{\delta_1} (b^1_1)^{\alpha_2} (b^1_2)^{\beta_2} (b^2_1)^{\gamma_2} (b^2_2)^{\delta_2} \; \middle| \;  \alpha_i, \beta_i, \gamma_i \in \mathbb{N}_0,\; \beta_i \text{ or } \gamma_i = 0 \,\right\}. \end{multline*}
\begin{proof}
Similar to above.
\end{proof}
\end{proposition}

We will need an alternative PBW basis for \(A_{\Sigma_{0,4}}\) in \cref{section:HilbertSeries}, so we shall now give an alternative basis for \(\mathscr{O}_q(\SL)\), and then use it to give the alternative PBW basis for \(A_{\Sigma_{0,4}}\). 

\begin{proposition}
\label{pbwbasisalt}
The set of monomials \[\left\{\, (a^2_1)^{\alpha} (a^1_1)^{\beta} (a^2_2)^{\gamma} (a^1_2)^{\delta} \; \middle| \;  \alpha, \beta, \gamma, \delta \in \mathbb{N}_0,\; \beta \text{ or }\gamma  = 0\, \right\}\] is a PBW basis for the reflection equation algebra \(\mathscr{O}_q(\SL)\) with respect to the ordering \(a^2_1 < a^1_1 < a^2_2 < a^1_2\).
\begin{proof}
A term rewriting system for \(\mathscr{O}_q(\SL)\) is
\begin{align*}
    \tau_{1211} : a^1_2 a^1_1 &\mapsto a^1_1 a^1_2 + q^{-2}(1-q^{-2})a^2_2 a^1_2, \\
    \tau_{1121} : a^1_1 a^2_1 &\mapsto a^2_1 a^1_1 - q^{-2}(1-q^{-2})a^2_1 a^2_2, \\
    \tau_{1221} : a^1_2 a^2_1 &\mapsto q^{-2}a^2_1 a^1_2 -q^{-2} (1-q^{-2})(1 - (a^2_2)^2), \\
    \tau_{2211} : a^2_2 a^1_1 &\mapsto a^1_1 a^2_2, \\
    \tau_{1222} : a^1_2 a^2_2 &\mapsto q^{-2} a^2_2 a^1_2, \\
    \tau_{2221} : a^2_2 a^2_1 &\mapsto q^{-2} a^2_1 a^2_2, \\
    \tau_{1122} : a^1_1 a^2_2 &\mapsto q^{-2} + a^2_1 a^1_2 + (1-q^{-2})(a^2_2)^2.
\end{align*}
The monomials given in the statement of the result are the reduced monomials with respect to this term rewriting system; furthermore, there are no inclusion ambiguities, and the overlap ambiguities are 
\begin{align*}
    (\tau_{1211}, \tau_{1121}, a^1_2, a^1_1, a^2_1), & \quad
    (\tau_{2211}, \tau_{1121}, a^2_2, a^1_1, a^2_1), \\
    (\tau_{1222}, \tau_{2211}, a^1_2, a^2_2, a^1_1), & \quad
    (\tau_{1222}, \tau_{2221}, a^1_2, a^2_2, a^2_1), \\
    (\tau_{2211}, \tau_{1122}, a^2_2, a^1_1, a^2_2), & \quad
    (\tau_{1211}, \tau_{1122}, a^1_2, a^1_1, a^2_2), \\
    (\tau_{1122}, \tau_{2211}, a^1_1, a^2_2, a^1_1), & \quad
    (\tau_{1122}, \tau_{2221}, a^1_1, a^2_2, a^2_1).  
\end{align*}

We shall order \(\mathscr{O}_q(\SL)\) with respect to the reduced degree where we give the generators the ordering \(a^2_1 < a^1_1 < a^2_2 < a^1_2\). This ordering is compatible with the given term rewriting systems and the rewriting will terminate, so if the ambiguities are resolvable then we can apply the diamond lemma, and we are done. It can be checked by direct calculation that the ambiguities are resolvable.
\end{proof}
\end{proposition}

\begin{corollary}
\label{altPBWbasis}
An alternative PBW basis for \(A_{\Sigma_{0,4}}\) is 
\begin{multline*}
\left\{\, (a^1_1)^{\alpha_1} (a^1_2)^{\beta_1} (a^2_1)^{\gamma_1} (a^2_2)^{\delta_1} (b^2_1)^{\alpha_2} (b^1_1)^{\beta_2} (b^2_2)^{\gamma_2} (b^1_2)^{\delta_2} (c^1_1)^{\alpha_3} (c^1_2)^{\beta_3} (c^2_1)^{\gamma_3} (c^2_2)^{\delta_3} \right| \\ 
\left| \; \alpha_i, \beta_i, \gamma_i \in \mathbb{N}_0, \beta_i \text{ or } \gamma_i = 0\,\right\}. \end{multline*}
\begin{proof}

The same as \cref{prop:PBWbasis} except we use the PBW basis
\[\left\{\,(b^2_1)^{\alpha} (b^1_1)^{\beta} (b^2_2)^{\gamma} (b^1_2)^{\delta} \; \middle| \; \alpha, \beta, \gamma, \delta \in \mathbb{N}_0, \beta \text{ or }\gamma  = 0 \,\right\}\]
\sloppy from \cref{pbwbasisalt} for the second copy of \(\mathscr{O}_q(\SL)\) in \(A_{\Sigma_{0,4}} = \mathscr{O}_q(\SL)^{\otimes 3}\).
\end{proof}
\end{corollary}
\subsection{The Algebra of \(\qgroup{\mathfrak{sl}_2}\)-Invariants of the Four-Punctured Sphere}
\label{section:InvAlgebraSphere}
We now turn to the first main result of this paper: giving a presentation of the algebra of \(\qgroup{\mathfrak{sl}_2}\)-invariants \(\mathscr{A}_{\Sigma_{0,4}}\) of \(\int_{\Sigma_{0,4}} \Repfd_q(\SL)\). As explained in \cref{section:InvAlgebra}, this algebra defines a \(\SL\)-quantum character variety of \(\Sigma_{0,4}\).

Recall from \cref{section:FactofSphere} that the generators of \(A_{\Sigma_{0,4}}\), organised into matrices, are:
\[
A := \begin{pmatrix}
a^1_1 & a^1_2 \\ a^2_1 & a^2_2
\end{pmatrix},\;
B := \begin{pmatrix}
b^1_1 & b^1_2 \\ b^2_1 & b^2_2
\end{pmatrix},\;
C:= \begin{pmatrix} 
c^1_1 & c^1_2 \\ c^2_1 & c^2_2 
\end{pmatrix}.
\]
Note that the quantum traces \(\Tr(A) = a^1_1 + q^{-2} a^2_2\), \(\Tr(B) = b^1_1 + q^{-2} b^2_2\) and \(\Tr(C) = c^1_1 + q^{-2} c^2_2\) of these matrices are invariant under the action of the quantum group on \(A_{\Sigma}\), and hence are contained in \(\mathscr{A}_{\Sigma_{0,4}}\). Furthermore, the quantum trace \(\Tr(X)\) of any matrix \(X = \sum_{i}^N A^{\alpha_i} B^{\beta_j} C^{\gamma_i}\) where \(\alpha_i, \beta_i, \gamma_i \in \mathbb{N}_0\) is also invariant under the action of the quantum group, so must also be contained in \(\mathscr{A}_{\Sigma_{0,4}}\). The quantum Cayley--Hamilton equation \(X^2 = \Tr(X)X - q^{-2} \det_q(X)\) implies that \(\Tr(X)\) is a linear combinations of the traces \(\Tr(A)\), \(\Tr(B)\), \(\Tr(C)\), \(\Tr(AB)\), \(\Tr(AC)\), \(\Tr(BC)\) and \(\Tr(ABC)\). Therefore, these seven traces generate all the invariants which are of the form \(\Tr(X)\). In this section we prove that these seven traces in fact generate the entire algebra of \(\qgroup{\mathfrak{sl}_2}\)-invariants \(\mathscr{A}_{\Sigma_{0,4}}\) and state the relations these traces satisfy.

\begin{definition}
\label{defn:Apres}
Let \(\mathscr{B}\) be the algebra with generators \(E, F, G, s, t, u, v\) subject to the relations: 
\newcommand{\lastgen}{
\begin{cases}
\begin{aligned}
        & -E^2 - q^{-4}F^2 - G^2 - q^{-4}(s^2 + t^2 + u^2 + v^2) \\
        & + (st + uv)E + q^{-2}(su +tv)F + (sv +tu)G \\
        & -stuv + q^{-6}(q^2+1)^2
\end{aligned}
\end{cases}
}

\begin{alignat}{8}
    &FE &{}={}&&
        q^2 &EF &{}+{}&&
        (q^2 - q^{-2})&G &{}+{}&&
        (1- q^2)&(sv + tu), \label{relation:FE}\\
    &GE &{}={}&&
        q^{-2}&EG &{}-{}&&
        q^{-2}(q^2-q^{-2})&F &{}+{}&&
        (1-q^{-2})&(su + tv), \label{relation:GE}\\
    &GF &{}={}&&
        q^2 & FG &{}+{}&&
        (q^2 - q^{-2})&E &{}+{}&&
        (1-q^2)&(st + uv), \label{relation:GF}\\
    & EFG &= \mathclap{\hphantom{\lastgen}\lastgen} \label{relation:EFG}
\end{alignat}
and \(s,t,u,v\) are central. 
\end{definition}

\begin{theorem} \label{thm:isom}
The map \(\Phi': \mathscr{B} \to  {\mathscr{A}_{\Sigma_{0,4}}}\) defined by:
\begin{center}
\vspace{-1em}
\begin{minipage}{0.3\linewidth}
    \begin{align*}
        E &\mapsto \Tr(AB),\\
        F &\mapsto \Tr(AC),\\
        G &\mapsto \Tr(BC),
    \end{align*}
\end{minipage}
\begin{minipage}{0.3\linewidth}
    \begin{align*}
        s &\mapsto \Tr(A),\\
        t &\mapsto \Tr(B),\\
        u &\mapsto \Tr(C),\\ 
        v &\mapsto \Tr(ABC),
    \end{align*}
\end{minipage}
\end{center}
is an isomorphism of algebras.
We denote by \(\Phi: \mathscr{B} \to \mathscr{O}_q^{\otimes 3}\) the map defined by the same formulas.
\end{theorem}

\begin{proof}[Proof of \cref{thm:isom}]
To check that \(\Phi\) is a morphism of algebras one must check that the images of relations (\ref{relation:FE})-(\ref{relation:EFG}) are satisfied in \(\mathcal{O}_q^{\otimes 3}\), which is a long but straightforward calculation, which we omit. As all quantum traces lie in \(\mathscr{A}_{\Sigma_{0,4}}\), the codomain of \(\Phi\) can be restricted to define \(\Phi'\) . So to show \(\Phi'\) is an isomorphism of algebras it remains to show \(\Phi'\) is a bijection which will be done by first proving \(\Phi\) is injective and then proving that both \(\mathscr{B}\) and \(\mathscr{A}_{\Sigma_{0,4}}\) have the same Hilbert series. 

The proof of injectivity of \(\Phi\) uses a filtration on the codomain \(\mathcal{O}_q^{\otimes 3}\).
\begin{definition}
\label{defn:filterationO}
We define a \emph{filtration on the algebra \(\mathcal{O}_q^{\otimes 3} = \bigcup_{i \in \mathbb{N}_0} F_i\)} by defining the degree of the generators as follows:
\begin{itemize}
    \item Degree 0: \(a^2_1\), \(a^2_2\), \(c^1_2\), and \(c^2_2\);
    \item Degree 1: \(a^1_1\), \(c^1_1\);
    \item Degree 2: \(a^1_2\),\(c^2_1\), \(b^1_1\), \(b^1_2\), \(b^2_1\), and \(b^2_2\).
\end{itemize}
\end{definition}

\begin{definition}
Let \(\mathcal{G}(\mathcal{O}_q^{\otimes 3}) = \bigoplus_{n \in \mathbb{N}_0} G_n\) denote the associated graded algebra of \(\mathcal{O}_q^{\otimes 3} = \cup_{i \in \mathbb{N}_0} F_i\).
\end{definition}

\begin{lemma}
The set of monomials \[\left\{ \, \Phi(E^{\epsilon} F^n G^m s^{\alpha} t^{\beta} u^{\gamma} v^{\delta}) \; \middle| \; \epsilon \text{ or } m\text{ or }n =0; \alpha, \beta, \gamma, \delta, n, m, \epsilon \in \mathbb{N}_0 \, \right\} \] is linearly independent in \(\mathcal{O}_q^{\otimes 3}\), so the homomorphism \(\Phi: \mathscr{B} \to \mathcal{O}_q^{\otimes 3}\) is injective.
\begin{proof}

Suppose to the contrary that the set 
\[
\left\{ \, \Phi\left(E^{\epsilon}F^n G^m s^{\alpha} t^{\beta} s^{\gamma} t^{\delta}\right) \; \middle| \; \epsilon \text{ or } m\text{ or }n =0;  \epsilon,m, n,\alpha, \beta, \gamma, \delta \in \mathbb{N}_{0}\, \right\} 
\]
is linearly dependent. Then for some finite indexing set \(I\) there exists scalars \(c_i\) which are not all zero such that
\begin{equation} \label{linsum0} \sum_{i \in I} c_i \Phi(E^{\epsilon_i}F^{n_i} G^{m_i} s^{\alpha_i} t^{\beta_i} u^{\gamma_i} v^{\delta_i}) = 0 \in \mathcal{O}_q^{\otimes 3}.\end{equation}
Map this sum to \(\mathcal{G}(\mathcal{O}_q^{\otimes 3})\):
\begin{equation} \label{linsum1} \sum_{i \in I} c_i \Phi(E^{\epsilon_i}F^{n_i} G^{m_i} s^{\alpha_i} t^{\beta_i} u^{\gamma_i} v^{\delta_i}) = 0 \in \mathcal{G}(\mathcal{O}_q^{\otimes 3}).\end{equation}
As \(s,t,u\) and \(v\) are central in \(\mathscr{B}\), (\ref{linsum1}) can be rearranged to give
\begin{equation} \label{linsum} \sum_{i \in I} c_i \Phi(s^{\alpha_i} E^{\epsilon_i} v^{\delta_i}  t^{\beta_i} F^{n_i}  u^{\gamma_i} G^{m_i}) = 0.\end{equation}
As \(\mathcal{G}(\mathcal{O}_q^{\otimes 3}) \) is graded, we can assume that all the terms in expression (\ref{linsum}) are in the maximal degree; we also know that
\begin{alignat*}{3}
    \Phi(X) &= \Tr(AB) &&= a^1_2 b^2_1 &\in \mathcal{G}_4, \\
    \Phi(F) &= \Tr(AC) &&= a^1_2 c^2_1 &\in \mathcal{G}_4, \\
    \Phi(G) &= \Tr(BC) &&= b^1_2 c^2_1 &\in \mathcal{G}_4, \\
    \Phi(s) &= \Tr(A) &&= a^1_1 &\in \mathcal{G}_1, \\
    \Phi(t) &= \Tr(B) &&= b^1_1 + q^{-1} b^2_2 \;&\in \mathcal{G}_2, \\
    \Phi(u) &= \Tr(C) &&= c^1_1 &\in \mathcal{G}_1, \\
    \Phi(v) &= \Tr(ABC) &&= a^1_2 b^2_2 c^2_1 &\in \mathcal{G}_6,
\end{alignat*}
so expression (\ref{linsum}) implies that:
\begin{equation} \label{linsum2}\sum_{i \in I, S(i) = N} c_i (a^1_1)^{\alpha_i} (a^1_2 b^2_1)^{\epsilon_i} (a^1_2 b^2_2 c^2_1)^{\delta_i} (b^1_1 +b^2_2)^{\beta_i} (a^1_2 c^2_1)^{n_i} (c^1_1)^{\gamma_i} (b^1_2 c^2_1)^{m_i}  = 0,\end{equation}
where \(S(i):= \alpha_i + \gamma_i + 4(\epsilon_i + n_i + m_i +\beta_i) + 6 \delta_i\) and \(N \in \mathbb{N}_{0}\). The crossing relations (\cref{prop:facthom}):

\begin{alignat*}{8}
    b^1_1  a^1_2 &= a^1_2  b^1_1 &\in \mathcal{G}_4, &\quad&
    b^2_1  a^1_2 &= q^{-2}a^1_2  b^2_1 \;&\in \mathcal{G}_4, \\
    b^2_2  a^1_2 &= a^1_2  b^2_2 &\in \mathcal{G}_4, &&
    b^2_2  b^1_1 &= b^1_1 b^2_2 &\in \mathcal{G}_4, \\
    c^1_1  b^1_2 &= b^1_2  c^1_1 &\in \mathcal{G}_3, &&
    c^1_2  b^2_2 &= b^2_2  c^1_2 &\in \mathcal{G}_2, \\
    c^2_1  a^1_2 &= q^{-2}a^1_2  c^2_1 \;&\in \mathcal{G}_2, &&
    c^2_1  b^1_1 &= b^1_1  c^2_1 &\in \mathcal{G}_4, \\
    c^2_1  b^1_2 &= q^{-2}b^1_2  c^2_1 &\in \mathcal{G}_4 , &&
    c^2_1  b^2_2 &= b^2_2  c^2_1 &\in \mathcal{G}_4, \\
    b^2_2 b^1_1 &= b^1_1 b^2_2 &\in \mathcal{G}_4, &&
    b^2_2 b^1_2 &= q^2 b^1_2 b^2_2 &\in \mathcal{G}_4, \\
    c^2_1 c^1_1 &= c^1_1 c^2_1 &\in \mathcal{G}_3,
\end{alignat*}
can be used to reorder the term in expression (\ref{linsum2}) to give
\begin{equation} \label{linsum3}\sum_{\substack{i \in I, \\S(i) = N}} \sum_{k=0}^{\beta_i} c_i q^{A_{i,k}} (a^1_1)^{\alpha_i} (a^1_2)^{\delta_i + \epsilon_i + \gamma_i} (b^2_1)^{\epsilon_i}  (b^1_1)^k (b^2_2)^{\beta_{i}-k + \delta_i } (b^1_2)^{m_i} (c^1_1)^{\gamma_i} (c^2_1)^{\delta_i + n_i + m_i} = 0,\end{equation}
for some constants \(A_{i,k} \in \mathbb{Z}\). 

Using the basis for \(A_{\Sigma_{0,4}}\) given in \cref{altPBWbasis}, the expression \((\ref{linsum3})\) is linear combination of distinct monomials which are in the basis of \(\mathcal{G}(\mathcal{O}^{\otimes 3})\), so all the coefficients must be zero. This is a contradiction as we assumed that not all the \(c_i\) were zero.
\end{proof}
\end{lemma}

In order to compute the Hilbert series of \(\mathscr{B}\), \(\mathscr{B}\) must be filtered.

\begin{definition}
\label{def:gradingB}
We define a \emph{filtration on the algebra \(\mathscr{B}\)} by defining the degree of the generators as follows:
\begin{itemize}
    \item Degree 1: \(s,t, u\);
    \item Degree 2: \(E, F, G\);
    \item Degree 3: \(v\).
\end{itemize}
\end{definition}

\begin{lemma}
The algebras \(\mathscr{B}\) and \(\mathscr{A}_{\Sigma_{0,4}}\) have the same Hilbert series when \(\mathscr{B}\) is given the filtration defined directly above and \(\mathscr{A}_{\Sigma_{0,4}}\) the filtration by degree.
\begin{proof}
The Hilbert series of \(\mathscr{A}_{\Sigma_{0,4}}\) is computed in \cref{section:HilbertSeries} and is \(\frac{1-t+t^2}{(1-t)^6(1+t)^2}\). By \cref{prop:PBWbasissphere}, a basis of \(\mathscr{G}(\mathscr{B})\) over \(\mathbb{C}[s, t, u, v]\) is \(\left\{ \, E^n F^m G^l \; \middle| \; n \text{ or } m\text{ or }l =0 \, \right\}\), so a basis of \(\mathscr{B}\) over \(\mathbb{C}\) is 
\[\left\{ \, E^n F^m G^l s^a t^b u^c v^d\; \middle| \; n \text{ or } m\text{ or }l =0;\; a, b, c, d, n, m, l \in \mathbb{N}_0 \, \right\}.\]
Therefore, there is a grading preserving vector space isomorphism
\begin{align*}
    \mathscr{G}(\mathscr{A}) &\to \langle E, F, G \rangle \otimes \mathbb{C}[s] \otimes \mathbb{C}[t] \otimes \mathbb{C}[u] \otimes \mathbb{C}[v]: \\
E^a F^b G^c s^d t^e u^f v^g &\mapsto (E^a F^b G^c) \otimes s^d \otimes t^e \otimes u^f \otimes v^g
\end{align*}
where \(\langle E, F, G \rangle\) is the subalgebra of \(\mathscr{A}\) generated by \(E,F,G\); hence,
\[
h_{\mathscr{A}}(t) = h_{\langle E, F, G \rangle}(t) \cdot h_{\mathbb{C}[s]}(t) \cdot h_{\mathbb{C}[t]}(t) \cdot h_{\mathbb{C}[u]}(t) \cdot h_{\mathbb{C}[v]}(t).
\]
If \(x=s,t,u\) the algebra \(\mathbb{C}[x]\) is the polynomial algebra graded by degree, so \((\mathbb{C}[x])[n]\) has basis \(\{\,x^n\, \}\), and 
\[
h_{\mathbb{C}[x]}(t) = \sum_{n=0}^{\infty} (\dim \left(\mathbb{C}[x]\right)[n])t^n = \sum_{n=0}^{\infty} t^n = \frac{1}{1-t}.
\]
The algebra \(\mathbb{C}[v]\) is the polynomial algebra graded by \(3\) times the degree, so \((\mathbb{C}[x])[n]\) has basis \(\left\{ \,x^{\frac{n}{3}}\, \right\}\) if \(n \equiv 0 \mod 3\) and \(\emptyset\) otherwise, and  
\[
h_{\mathbb{C}[v]}(t) = \sum_{n=0}^{\infty} (\dim \left(\mathbb{C}[x]\right)[n])t^n = \sum_{n=0}^{\infty} t^{3n} = \frac{1}{1-t^3}.
\]
The algebra \(\langle E, F, G \rangle[k]\) has basis 
\[
\left\{ \, E^a F^b G^c \; \middle| \; a+b+c=n; a\text{ or }b\text{ or }c\text{ is }0 \, \right\}
\]
if \(k=2n\) is even and the basis is \(\emptyset\) otherwise. Assume \(k\) is even so \(k=2n\). If \(n=0\) then the basis has one element \(\left\{ \,0\, \right\}\). If \(n\neq0\) then the basis is
\begin{align*}
&\left\{ \, E^a F^b G^c\; \middle| \; a+b+c=n; a\text{ or }b\text{ or }c\text{ is }0 \,\right\} \\
&= \left\{ \, E^a F^b G^c \; \middle| \; a+b+c=n; \text{ one of } a,b,c \text{ is }0 \, \right\} \\
&\quad \sqcup \left\{ \, E^a F^b G^c \; \middle| \; a+b+c=n; \text{ two of } a,b,c \text{ is }0 \, \right\} \\
&= \left\{ \, E^a F^b \; \middle| \; a+b=n; a,b \neq 0 \, \right\} \sqcup \left\{ \, F^b G^c \; \middle| \; b+c=n; b,c \neq 0 \, \right\} \\
& \quad \sqcup \left\{ \, E^a G^c \; \middle| \; a+c=n; a,c \neq 0 \, \right\}\sqcup \left\{ \, E^n, F^n, G^n \, \right\}
\end{align*}
which has \(3n\) elements. Hence, the Hilbert series of \(\langle E, F, G \rangle\) is 
\[h_{\langle E, F, G \rangle}(t) = \sum_{n=0}^{\infty} (\dim \left(\langle E, F, G \rangle \right)[n])t^n = 1+ \sum_{n=1}^{\infty} 3n t^{2n} = 1+\frac{3t^2}{(1-t^2)^2}.\]
Thus 
\begin{align*}
h_{\mathscr{A}_{\Sigma_{0,4}}}(t) 
&= h_{\langle E, F, G \rangle}(t) \cdot h_{\mathbb{C}[s]}(t) \cdot h_{\mathbb{C}[t]}(t) \cdot h_{\mathbb{C}[u]}(t) \cdot h_{\mathbb{C}[v]}(t) \\
&= \left( 1+\frac{3t^2}{(1-t^2)^2} \right) \frac{1}{(1-t)^3(1-t^3)} \\
&= \frac{1-t+t^2}{(1-t)^6(1+t)^2},
\end{align*}
which means that \(\mathscr{B}\) and \(\mathscr{A}_{\Sigma_{0,4}}\) have the same Hilbert series. 
\end{proof}
\end{lemma}
The homomorphism \(\Phi'\) is filtered if we give \(\mathscr{B}\) the filtration defined in \cref{def:gradingB} and \(\mathscr{A}_{\Sigma_{0,4}}\) the filtration by degree. It is injective and the Hilbert series of \(\mathscr{B}\) and \(\mathscr{A}_{\Sigma_{0,4}}\) are equal, so \(\Phi'\) is an isomorphism. This concludes the proof of \cref{thm:isom}.
\end{proof}
\subsection{The Algebra of \(\qgroup{\mathfrak{sl}_2}\)-Invariants of the Punctured Torus}
We now obtain a presentation of the algebra of \(\qgroup{\mathfrak{sl}_2}\)-invariants for our second surface, the punctured torus. This is simpler than the four-punctured torus case, and the proofs follow in a similar manner. 

\begin{definition}
\label{defn:Tpres}
Let \(\mathscr{T}\) be the algebra with generators \(X, Y, Z\) and relations: 
\begin{align*}
YX - q^{-1}XY &= (q-q^{-1})Z; \\
XZ - q^{-1}ZX &= - q^{-3}(q-q^{-1})Y; \\
ZY - q^{-1}YZ &= -q^{-3}(q - q^{-1})X.
\end{align*}
It has a central element 
\[L := q^5 XZY + q^3 Y^2 - q^4 Z^2 + q^3 X^2 - (q-q^{-1}).\]
\end{definition}

\begin{proposition}
\label{prop:PBWbasisT}
The set of monomials \[\left\{ \, X^{\alpha} Y^{\beta} Z^{\gamma} \; \middle| \; \alpha, \beta, \gamma\in \mathbb{N}_0\, \right\}\] 
is a PBW basis for the algebra \(\mathscr{T}\).
\end{proposition}

\begin{proof}
We use the reduced degree with the generators ordered by \(X<Y<Z\) as our ordering. From the relations of \(\mathscr{T}\) we obtain the term rewriting system
\begin{align*}
\sigma_{YX} &: YX \mapsto q^{-1}XY + (q-q^{-1})Z; \\
\sigma_{ZX} &: ZX \mapsto qXZ +  q^{-2}(q-q^{-1})Y; \\
\sigma_{ZY} &: ZY \mapsto q^{-1}YZ -q^{-3}(q - q^{-1})X.
\end{align*}
this term rewriting system is compatible with the ordering, and its only ambiguity \((\sigma_{ZY}, \sigma_{YX}, Z, X, Y)\) is resolvable, so by the diamond lemma the reduced monomials \(\left\{ \,X^{\alpha} Y^{\beta} Z^{\gamma} \; \middle| \; \alpha, \beta, \gamma\in \mathbb{N}_0\, \right\}\) form a PBW basis for the algebra. 
\end{proof}

Organise the generators of \(A_{\Sigma_{1,1}}\) into matrices as follows:
\[
A := \begin{pmatrix}
a^1_1 & a^1_2 \\ a^2_1 & a^2_2
\end{pmatrix},
B := \begin{pmatrix}
b^1_1 & b^1_2 \\ b^2_1 & b^2_2
\end{pmatrix}.
\]

\begin{theorem}
\label{thm:torusiso}
Define the map $\Psi: \mathscr{T} \to \mathcal{O}_q^{\otimes 2}$ by
\begin{align*}
X &\mapsto \Tr(A), \\
Y &\mapsto \Tr(B), \\
Z &\mapsto \Tr(AB).
\end{align*}
The restricted map $\Psi': \mathscr{T} \to \mathscr{A}_{\Sigma_{1,1}}$ is an algebra isomorphism.
\end{theorem}

\begin{proof}
To check that \(\Psi\) is a morphism of algebras one must check that the images of the three relations are satisfied in \(\mathcal{O}_q^{\otimes 2}\), which is a long but straightforward calculation. As all quantum traces lie in \(\mathscr{A}_{\Sigma_{1,1}}\), the codomain of \(\Psi\) can be restricted to define \(\Psi'\) . So to show \(\Psi'\) is an isomorphism of algebras it remains to show \(\Psi'\) is a bijection which will be done by proving \(\Psi\) is injective and that both \(\mathscr{T}\) and \(\mathscr{A}_{1,1}\) have the same Hilbert series. 

\begin{lemma}
The set of monomials \[\left\{ \,\Psi\left(X^{\alpha} Y^{\beta} Z^{\gamma}\right) \; \middle| \; \alpha, \beta, \gamma \in \mathbb{N}_{0}\, \right\} \] is linearly independent in \(\mathcal{O}_q^{\otimes 2}\), so the homomorphism \(\Psi: \mathscr{T} \to \mathcal{O}_q^{\otimes 2}\) is injective.
\begin{proof}
In this proof we use the filtration in defined in \cref{defn:filterationO} restricted to \(\mathcal{O}_q^{\otimes 2}\). Suppose to the contrary that the set 
\[\left\{ \,\Psi\left(X^{\alpha} Y^{\beta} Z^{\gamma}\right) \; \middle| \; \alpha, \beta, \gamma \in \mathbb{N}_{0}\, \right\} \]
is linearly dependent then for some finite indexing set \(I\) there exists scalars \(c_i\) which are not all zero such that
\begin{equation} \label{linsumtorus0} \sum_{i \in I} c_i \Psi(X^{\alpha_i} Y^{\beta_i} Z^{\gamma_i})  = 0 \in \mathcal{O}_q^{\otimes 2}.\end{equation}
Map this to \(\mathcal{G}(\mathcal{O}_q^{\otimes 2})\):
\begin{equation} \label{linsumtorus1} \sum_{i \in I} c_i \Psi(X^{\alpha_i} Y^{\beta_i} Z^{\gamma_i})  = 0 \in \mathscr{G}(\mathcal{O}_q^{\otimes 2}).\end{equation}
As \(\mathcal{G}(\mathcal{O}_q^{\otimes 2}) \) is graded, we can assume that all the terms in expression (\ref{linsumtorus1}) are in the maximal degree; we also know that

\begin{alignat*}{4}
	\Phi(X) &= \Tr(A) &&= a^1_1 &\in \mathcal{G}_1, \\
    \Phi(Y) &= \Tr(B) &&= b^1_1 + q^{-1} b^2_2 \;&\in \mathcal{G}_2, \\
    \Phi(Z) &= \Tr(AB) &&= a^1_2 b^2_1 &\in \mathcal{G}_4,
\end{alignat*}

so expression (\ref{linsumtorus1}) implies that:
\begin{equation} \label{linsumtorus2}\sum_{i \in I, S(i) = N} c_i (a^1_1)^{\alpha_i} (b^1_1 + q^{-1}b^2_2)^{\beta_i}  (a^1_2 b^2_1)^{\gamma_i}  = 0,\end{equation}
where \(S(i):= \alpha_i + 4(\beta_i + \gamma_i) \) and \(N \in \mathbb{N}_{0}\). The crossing relations
\begin{alignat*}{6}
    b^1_1  a^1_2 &= a^1_2  b^1_1 &&\in \mathcal{G}_4, \quad&
    b^2_1  a^1_2 &= q^{-2}a^1_2  b^2_1 &&\in \mathcal{G}_4, \\
    b^2_2  a^1_2 &= a^1_2  b^2_2 &&\in \mathcal{G}_4, &
    b^2_2  b^1_1 &= b^1_1 b^2_2 &&\in \mathcal{G}_4, \\
    b^2_2 b^1_2 &= q^2 b^1_2 b^2_2 &&\in \mathcal{G}_4,
\end{alignat*}
can be used to reorder the term in expression (\ref{linsumtorus2}) to give
\begin{equation} \label{linsumtorus3}\sum_{\substack{i \in I, \\S(i) = N}} \sum_{k=0}^{\beta_i} c_i q^{A_{i,k}} (a^1_1)^{\alpha_i} (a^1_2)^{\gamma_i} (b^1_1)^k (b^2_1)^{\gamma_i} (b^2_2)^{\beta_i - k} = 0,\end{equation}
for some constants \(A_{i,k} \in \mathbb{Z}\). 

Using the basis for \(A_{\Sigma_{1,1}}\) given in \cref{prop:pbwbasistorus}, the expression (\ref{linsumtorus3}) is linear combination of distinct monomials which are in the basis of \(\mathcal{G}(\mathcal{O}^{\otimes 2})\), so all the coefficients must be zero. This is a contradiction as we assumed that not all the \(c_i\) were zero.
\end{proof}
\end{lemma}

In order to compute the Hilbert series of \(\mathscr{T}\), \(\mathscr{T}\) must be filtered.

\begin{definition}
\label{defn:gradingofT}
We define a \emph{filtration on the algebra \(\mathscr{T}\)} by defining the degree of the generators as follows:
\begin{itemize}
    \item Degree 1: \(X, Y\);
    \item Degree 2: \(Z\).
\end{itemize}
\end{definition}

\begin{lemma} \label{gradingofT} 
The associated graded algebra \(\mathscr{G}(\mathscr{T})\) has a PBW basis 
\[\left\{ \, X^{\alpha} Y^{\beta} Z^{\gamma} \; \middle| \; \alpha, \beta, \gamma \in \mathbb{N}_0 \, \right\}.\]
\end{lemma}
\begin{proof}
The associated graded algebra \(\mathscr{G}(\mathscr{T})\) is the algebra with generators \(X, Y, Z\) subject to the relations: 
\[YX =  q^{-1}XY + (q-q^{-1})Z; \quad XZ = q^{-1}ZX; \quad ZY = q^{-1}YZ;\]
We can apply the diamond lemma with the above relations as the term rewriting system. 
\end{proof}

\begin{lemma}
The algebras \(\mathscr{T}\) and \(\mathscr{A}_{\Sigma_{1,1}}\) have the same Hilbert series when \(\mathscr{T}\) is given the filtration in \cref{defn:gradingofT} and \(\mathscr{A}_{\Sigma_{1,1}}\) the filtration by degree.
\begin{proof}
The Hilbert series of \(\mathscr{A}_{\Sigma_{1,1}}\) is computed in \cref{section:HilbertSeries} and is \(\frac{1}{(1-t)^2(1-t^2)}\). We note from \cref{gradingofT} that 
\[\left\{ \, X^{\alpha} Y^{\beta} Z^{\gamma} \; \middle| \; \alpha, \beta, \gamma \in \mathbb{N}_0 \, \right\}.\]
is a basis of \(\mathscr{G}(\mathscr{T})\), so there is a grading preserving vector space isomorphism
\begin{align*}
    \mathscr{G}(\mathscr{T}) &\to \mathbb{C}[X] \otimes \mathbb{C}[Y] \otimes \mathbb{C}[X]: \\
X^{\alpha} Y^{\beta} Z^{\gamma} &\mapsto X^{\alpha} \otimes Y^{\beta} \otimes Z^{\gamma};
\end{align*}
hence,
\[
h_{\mathscr{T}}(t) = h_{\mathbb{C}[X]}(t) \cdot h_{\mathbb{C}[Y]}(t) \cdot h_{\mathbb{C}[Z]}(t).
\]
If \(x=X, Y\) the algebra \(\mathbb{C}[x]\) is the polynomial algebra graded by degree, so \((\mathbb{C}[x])[n]\) has basis \(\left\{ \,x^n\, \right\}\), and 
\[
h_{\mathbb{C}[x]}(t) = \sum_{n=0}^{\infty} (\dim \left(\mathbb{C}[x]\right)[n])t^n = \sum_{n=0}^{\infty} t^n = \frac{1}{1-t}.
\]
The algebra \(\mathbb{C}[Z]\) is the polynomial algebra graded by two times the degree, so \((\mathbb{C}[Z])[n]\) has basis \(\left\{ \,Z^{\frac{n}{2}}\, \right\}\) if \(n \equiv 0 \mod 2\) and \(\emptyset\) otherwise, and  
\[
h_{\mathbb{C}[Z]}(t) = \sum_{n=0}^{\infty} (\dim \left(\mathbb{C}[Z]\right)[n])t^n = \sum_{n=0}^{\infty} t^{2n} = \frac{1}{1-t^2}.
\]Thus 
\begin{align*}
h_{\mathscr{T}}(t) 
&=  h_{\mathbb{C}[X]}(t) \cdot h_{\mathbb{C}[Y]}(t) \cdot h_{\mathbb{C}[Z]}(t) \\
&= \frac{1}{(1-t)^2(1-t^2)},
\end{align*}
which means that \(\mathscr{T}\) and \(\mathscr{A}_{\Sigma_{1,1}}\) have the same Hilbert series. 
\end{proof}
\end{lemma}
The homomorphism \(\Psi'\) is filtered if we give \(\mathscr{T}\) the filtration in \cref{gradingofT} and \(\mathscr{A}_{\Sigma_{1,1}}\) the filtration by degree. It is injective and the Hilbert series of \(\mathscr{T}\) and \(\mathscr{A}_{\Sigma_{1,1}}\) are equal, so \(\Psi'\) is an isomorphism. This concludes the proof of \cref{thm:torusiso}.
\end{proof}

\section{Isomorphisms}
\subsection{Isomorphisms with Skein Algebras, Spherical Double Affine Hecke Algebras and Cyclic Deformations}
\label{section:DAHA}

In this section we use the presentation of the algebras of \(\qgroup{\mathfrak{sl}_2}\)-invariants \(\mathscr{A}_{0,4}\) and  \(\mathscr{A}_{1,1}\) of the four-punctured sphere \(\Sigma_{0,4}\) and punctured torus \(\Sigma_{1,1}\) over \(\qgroup{\frsl}\) obtained in the previous section. We state isomorphisms between \(\mathscr{A}_{0,4}\) and two isomorphic algebras: \(S\mathscr{H}_{q,\underbar{t}}\), the spherical double affine Hecke algebra of type \((C^{\vee}_1, C_1)\), and \(\sk(\Sigma_{0,4})\), the Kauffman bracket skein algebra of the four-punctured sphere. We also state isomorphisms between \(\mathscr{A}_{1,1}\) and two isomorphic algebras: \(U_q(\mathfrak{su}_2)\), a cyclic deformation of \(U(\mathfrak{su}_2)\), and \(\sk(\Sigma_{1,1})\), the Kauffman bracket skein algebra of the punctured torus.

\subsubsection*{The Kauffman Bracket Skein Algebra}

\begin{definition}
The \emph{Kauffman bracket skein module \(\sk_q
(M)\)} of an oriented \(3\)-manifold \(M\) (possibly with boundary) is the vector space of formal linear sums of isotopy classes of framed links without contractible components in \(M\) (but including the empty link) on which we impose the Kauffman bracket skein relations: 
\begin{align*}
\smalldiagram{g4547} &= q^{-1} \smalldiagram{g4563} + q\; \smalldiagram{g4571}, \\ 
\smalldiagram{g4617} &= -q^{2} -q^{-2}.
\end{align*}
\end{definition}

Whilst in general it is difficult to find explicit presentations for skein algebras, presentations for the Kauffman bracket skein algebras of our surfaces, \(\Sigma_{0,4}\) and \(\Sigma_{1,1}\), are known.

\begin{definition} \label{defn:skeinalgebra}
The \emph{Kauffman bracket skein algebra \(\sk(\Sigma)\)} of the surface \(\Sigma\) is the Kauffman bracket skein module \(\sk(\Sigma \times [0,1])\). It is an algebra with multiplication given by stacking copies of \(\Sigma \times [0,1]\) on top of each other and retracting.  
\end{definition}

\begin{figure}[H]
\centering
\includegraphics[height=3cm]{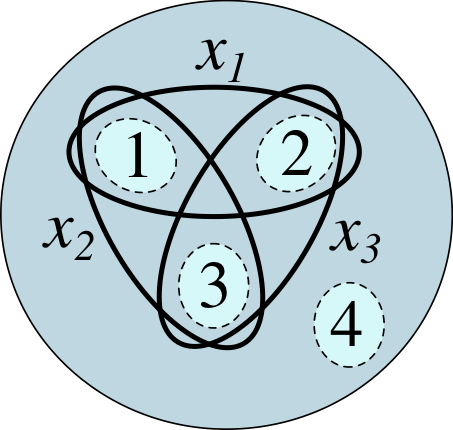}
\caption{The loops \(x_1,  x_2\) and \(x_3\)}
\label{figure:FPuncSphereEdges}
\end{figure}

\begin{theorem}{{\cite[Theorem~3.1]{Bullock&Przytycki}~\cite[Theorem~2.4]{Berest&Samuelson}}}
\label{thm:presofskeinspehre}
Let \(p_i\) denote the loops around the four punctures of \(\Sigma_{0,4}\) and let \(x_i\) denote the loops around punctures 1 and 2, 2 and 3, 1 and 3 respectively (see \cref{figure:FPuncSphereEdges}). The Kauffman bracket skein algebra \(\sk(\Sigma_{0,4})\) has a presentation where the generators are \(x_i\) and \(p_i\), and the relations are 
\begin{align*}
    [x_i, x_{i+1}]_{q^2} &= (q^4 - q^{-4})x_{i+2} - (q^2 - q^{-2})p_i \text{ (indices taken modulo 3)}; \\
    \Omega_K &= (q^2 +q^{-2})^2 - (p_1 p_2 p_3 p_4 + p_1^2 + p_2^2 + p_3^2 + p_4^2);
\end{align*}
where \([a, b]_q := qab - q^{-1}ba\) is the quantum Lie bracket and 
\[\Omega_K := -q^2 x_1 x_2 x_3 + q^4 x_1^2 + q^{-4} x_2^2 + q^4 x_3^2 + q^2 p_1 x_1 + q^{-2}p_2 x_2 + q^2 p_3 x_3.\]
\end{theorem}

\begin{theorem}[{\cite[Theorem~2.1]{Bullock&Przytycki}}]
\label{thm:skeintoruspres}
The Kauffman bracket skein algebra \(\sk(\Sigma_{1,1})\) has a presentation with generators \(x_1, x_2, x_3\) and relations 
\[[x_i, x_{i+1}]_{q} = (q^2-q^{-2}) x_{i+2} \text{ (indices taken modulo 3)}.\]
\end{theorem}

\subsubsection*{The Spherical Double Affine Hecke Algebras \(S\mathscr{H}_{q,\underline{t}}\) and \(SH_{q,t}\), and the Cyclic Deformation of \(U(\mathfrak{su}_2)\)}
Double Affine Hecke Algebras (DAHAs) were introduced by Cherednik \cite{cherednik1992double}, who used them to prove Macdonald's constant term conjecture for Macdonald polynomials, but have since found wider ranging applications particularly in representation theory \cite{CerednikIntro, CerednikJones}. DAHAs can be associated to different root systems with Cherednik's original DAHA being associated to the \(A^1\) root system. 

\begin{definition}
The \emph{\(A^1\) double affine Hecke algebra (DAHA)} \(H_{q,t}\) is the algebra with generators \(X^{\pm 1}\), \(Y^{\pm 1}\) and \(T\), and relations 
\[T X T = X^{-1}, \quad T Y^{-1} T = Y, \quad XY = q^2 Y X T^2, \quad (T - t)(T + t^{-1}) = 0.\]
\end{definition}

The element \(e = (T + t^{-1})/(t + t^{-1})\) is an idempotent of \(H_{q,t}\), and is used to define the spherical subalgebra \(SH_{q, t} := e H_{q, t} e\). 

\begin{theorem}[{\cite[11]{samuelson2014}~\cite[Section~2]{Terwilliger}}]
The spherical double affine Hecke algebra \(SH_{q, t}\) has a presentation with generators \(x,y,z\) and relations 
\[[x, y]_q = (q^2 - q^{-2}) z, \quad  [z, x]_q = (q^2 - q^{-2})y, \quad [y, z]_q = (q^2 - q^{-2})x\]
\[q^2 x^2 + q^{-2}y^2 + q^2 z^2 - q x y z = \left( \frac{t}{q} - \frac{q}{t} \right)^2 + \left(q + \frac{1}{q} \right)^2\]
where \([a,b]_q := qab - q^{-1} ba\) is the quantum Lie bracket.
\end{theorem}

The double affine Hecke algebra \(\mathscr{H}_{q,\underbar{t}}\) of type \((C^{\vee}_1, C_1)\) is a \(5\)-parameter deformation of the affine Weyl group \(\mathbb{C}[X^{\pm}, Y^{\pm}] \rtimes \mathbb{Z}_2\) with deformation parameters \(q \in \mathbb{C}^{\times}\) and \(\underbar{t} = (t_1, t_2, t_3, t_4) \in (\mathbb{C}^*)^4\). It can be given an abstract presentation with generators are \(T_0, T_1, T_0^{\vee}, T_1^{\vee}\) and relations:
\begin{align*}
     (T_0-t_1)(T_0 + t^{-1}_1) &= 0, \\
     (T_0^{\vee}-t_2)(T_0^{\vee} + t^{-1}_2) &= 0, \\
     (T_1-t_3)(T_1 + t^{-1}_3) &= 0, \\
     (T_1^{\vee}-t_4)(T_1^{\vee} + t^{-1}_4) &= 0, \\
     T_1^{\vee} T_1 T_0 T_0^{\vee} &=q.
\end{align*}
It generalises Cherednik's double affine Hecke algebras of rank 1 as \(H_{q;t}:= \mathscr{H}_{q, (1,1,t^{-1},1)}\). The element \(e = (T_1 + t^{-1}_3)/(t_3 + t_3^{-1})\) is an idempotent of \(\mathscr{H}_{q,\underbar{t}}\), and is used to define the spherical subalgebra \(S\mathscr{H}_{q, \underbar{t}} := e \mathscr{H}_{q, \underbar{t}} e\).

\begin{theorem}{{\cite[Theorem~2.20]{Berest&Samuelson}}} \label{thm:sDAHA pres}
The spherical double affine Hecke algebra \(S\mathscr{H}_{q, \underbar{t}}\) of type \((C^{\vee}_1, C_1)\) has a presentation with generators \(x,y,z\) and relations
\begin{align*}
    [x,y]_q &= (q^2 - q^{-2})z - (q-q^{-1})\gamma \\
    [y,z]_q &= (q^2 - q^{-2})x - (q-q^{-1})\alpha \\
    [z,x]_q &= (q^2 - q^{-2})y - (q-q^{-1})\beta \\ 
    \Omega &= \overline{t_1}^2 + \overline{t_2}^2 + \overline{q t_3}^2 + \overline{t_4}^2 - \overline{t_1} \overline{t_2} (\overline{q t_3}) \overline{t_4} + (q + q^{-1})^2
\end{align*}
where 
\begin{align*}
    \alpha &:= \overline{t_1} \overline{t_2} + \overline{q t_3} \overline{t_4}, \\ 
    \beta &:= \overline{t_1} \overline{t_4} + \overline{q t_3} \overline{t_2}, \\
    \gamma &:= \overline{t_2} \overline{t_4} + \overline{q t_3} \overline{t_1}, \\
    \Omega &:= -qxyz + q^2 x^2 + q^{-2} y^2 + q^2 z^2 - q \alpha x - q^{-1} \beta y - q \gamma z, \\
    [a, b]_q &:= qab - q^{-1}ba  \text{ is the quantum Lie bracket. } \\
\end{align*}
\end{theorem}

Using the presentation for the Kauffman bracket skein algebra \(\sk(\Sigma_{0,4})\) (\cref{thm:presofskeinspehre}) and the type \((C^{\vee}_1, C_1)\) spherical DAHA above, it is easy to see:

\begin{corollary}[{\cite[Corollary 2.10]{Berest&Samuelson}}]
\label{prop:Bereset&SamuelsonIso}
There is an isomorphism \(\delta: \sk(\Sigma_{0,4}) \to S\mathscr{H}_{q,\underline{t}}\) given by 
\begin{alignat*}{4}
    \beta(x_1) &= x,&\quad& \beta(p_1) &\;=& \;i \overline{t_1},\\
    \beta(x_2) &= y,&& \beta(p_2) &=&\;i \overline{t_2},\\
    \beta(x_3) &= z,&& \beta(p_3) &=&\;i \overline{q t_3},\\
   \beta(q) &= q^2, && \beta(p_4) &=&\;i \overline{t_4}.
\end{alignat*}
\end{corollary}

We now define the cyclic deformation of \(U(\mathfrak{su_2})\) and relate it to \(\sk(\Sigma_{1,1})\).

\begin{definition}[{\cite[3]{Bullock&Przytycki}~\cite[5]{zachos1990quantum}}]
The \emph{cyclic deformation of \(U(\mathfrak{su}_2)\)} is given by 
\[U_q(\mathfrak{su}_2) := \mathbb{C} \langle y_1, y_2, y_3 | [y_i, y_{i+1}]_q = y_{i+2} \rangle.\]where indices are taken modulo 3. 
\end{definition}

\begin{proposition}[{\cite[Corollary~2.2]{Bullock&Przytycki}}]
\label{prop:cyclicskein}
When \( (q^2 - q^{-2})\) is non-invertible there is an isomorphism 
\[\nu: \sk(\Sigma_{1,1}) \to U_q(\mathfrak{su}_2): x_i \mapsto (q^2 - q^{-2}) y_i.\]
\end{proposition}

Note that the element \(q^2 x_1^2 + q^{-2}x_2^2 + q^2 x_3^2 - q x_1 x_2 x_3\) is central in \(U_q(\mathfrak{su}_2)\) and setting it equal to \(\left( \frac{t}{q} - \frac{q}{t} \right)^2 + \left(q + \frac{1}{q} \right)^2 \) recovers the spherical DAHA \(SH_{q, t}\). 

\subsubsection*{Relation to Algebra of \(\qgroup{\mathfrak{sl}_2}\)-Invariants}
\begin{proposition}
\label{prop:sphereinviso}
There is an isomorphism \(\alpha: S\mathscr{H}_{q,\underline{t}} \to \mathscr{A}_{\Sigma_{0,4}}\) defined by
\begin{alignat*}{3}
    \alpha(x) &= -q E,&\quad& \alpha(\overline{t_1}) &= iqs,\\
    \alpha(y) &= -q F,&& \alpha(\overline{t_2}) &= iqt,\\
    \alpha(z) &= -q G,&&  \alpha(\overline{q t_3}) &= iqv,\\
    & && \alpha(\overline{t_4}) &= iqu.
\end{alignat*}
\end{proposition}
\begin{proof}
By rewriting the relations in the presentation of \(\mathscr{A}_{\Sigma}\) given in \cref{defn:Apres} in terms of the quantum Lie bracket \([\cdot , \cdot]_q\), we see that the algebra of \(\qgroup{\mathfrak{sl}_2}\)-invariants \(\mathscr{A}_{\Sigma}\)  has generators \(E, F, G, u, v, s, t\) and relations:
\begin{align*}
    [E, F]_q &= - q^{-1}(q^2 - q^{-2}) G + (q-q^{-1})(sv+tu) \\
    [F, G]_q &= - q^{-1}(q^2 - q^{-2}) E + (q-q^{-1})(st+uv) \\
    [G, E]_q &= - q^{-1}(q^2 - q^{-2}) F + (q -q^{-1})(su + tv) \\
    \tilde{\Omega} &= -q^2 s^2 + -q^{2}t^2 - q^{2}u^2 - q^{2}v^2 - q^4 stuv + q^{-2}(q^2 + 1)^2
\end{align*}
where
\begin{align*}
\tilde{\Omega} &= q^4 E F G - q^4 (st + uv) E - q^{2}(su +tv) F - q^4(sv + tu) G \\
&+ q^4 E^2 + F^2 + q^4 G^2.
\end{align*}
Also note that
\begin{align*}
    \alpha(\Omega) &=  \alpha(-qxyz + q^2 x^2 + q^{-2} y^2 + q^2 z^2 - q \alpha x - q^{-1} \beta y - q \gamma z) \\
    &= q^4 EFG + q^4 E^2 + F^2 + q^4 G^2 - q^4 (st+uv) E -q^2 (su+tv) F - q^4 (sv+tu) G \\
    &= \tilde{\Omega}.
\end{align*}

The map \(\alpha\) is clearly bijective, so it remains to show it is a algebra homomorphism:
\begin{align*}
    &\alpha \left( [x,y]_q - (q^2 - q^{-2})z + (q-q^{-1})\gamma \right) \\
    &\quad = q^2 [E,F]_q + (q^2 - q^{-2}) q^2 G - (q-q^{-1}) q^2 (sv +tu) \\
    &\quad= q^2 \left([E,F]_q + (q^2 - q^{-2}) G - (q-q^{-1}) (sv +tu) \right) \\ &\quad= 0 \\
    &\text{and similarly for the next two relations. For the final relation:} \\
    & \alpha(\overline{t_1}^2 + \overline{t_2}^2 + \overline{q t_3}^2 + \overline{t_4}^2 - \overline{t_1} \overline{t_2} \overline{q t_3} \overline{t_4} + (q + q^{-1})^2) - \Omega))\\
    &\quad= -q^2 s^2 - q^2 t^2 - q^2 v^2 -q^2 u^2 - q^4 stuv + (q + q^{-1})^2) - \tilde{\Omega} \\
    &\quad=0.
\end{align*}
\end{proof}

\begin{corollary}
\label{cor:skeinalginsphere}
There is an isomorphism  \(\beta: \sk_q(\Sigma_{0,4}) \to \mathscr{A}_{\Sigma_{0,4}}\) defined by 
\begin{alignat*}{3}
    \beta(x_1) &= -q E,&\quad& \beta(p_1) &\;=& \; -q s,\\
    \beta(x_2) &= -q F,&& \beta(p_2) &=&\; - q t,\\
    \beta(x_3) &= -q G,&& \beta(p_3) &=&\; - q v,\\
   \beta(q) &= q^2, && \beta(p_4) &=&\; - q u.
\end{alignat*}
\end{corollary}
\begin{proof}
Immediate from \cref{prop:Bereset&SamuelsonIso} and \cref{prop:sphereinviso}.
\end{proof}

\begin{proposition}
\label{prop:skeinalgintorus}
There is an isomorphism  \(\gamma: \mathscr{A}_{\Sigma_{1,1}} \to \sk(\Sigma_{1,1})\) defined by 
\begin{align*}
\gamma(q) &= q^2, \\
\gamma(X) &= i q^{-2} x_2, \\
\gamma(Y) &= i q^{-2} x_1, \\
\gamma(Z) &= -q^{-5} x_3.
\end{align*}
\begin{proof}
Immediate from \cref{thm:torusiso} and \cref{thm:skeintoruspres}.
\end{proof}
\end{proposition}
Hence by \cref{prop:cyclicskein}, \(\mathscr{A}_{1,1}\) is also isomorphic to \(U_q(\mathfrak{su}_2)\). 
\subsection{Isomorphism with a Quantisation of the Moduli Space of Flat Connections}
\label{section:FlatConnections}

Teschner and Vartanov proposed a quantisation for the \(\SL\)-character varieties of surfaces \cite{Teschner} by stating generators and relations for the quantisation of \(\Chg_{\SL}(\Sigma_{0,4})\) and \(\Chg_{\SL}(\Sigma_{1,1})\). The quantisation for other surfaces is then given by decomposing the surface into such surfaces. In this section we shall briefly outline this decomposition before showing that Teschner and Vartanov's quantisation of \(\Chg_{\SL}(\Sigma)\) coincides with the algebra of \(\qgroup{\mathfrak{sl}_2}\)-invariants \(\mathscr{A}_{\Sigma}\) quantisation for the base cases \(\Sigma=\Sigma_{0,4}\) and \(\Sigma_{1,1}\).

\begin{definition}
The Poisson algebra of algebraic functions on \(\Chg_G(\Sigma)\) is denoted \(\mathcal{A}(\Sigma)\).
\end{definition}

\begin{definition}
We can associate to the Riemann surface \(\Sigma\) a \emph{pants decomposition} \(\sigma = (C_{\sigma}, \Gamma_{\sigma})\) where:
\begin{enumerate}
\item The \emph{cut system} \(C_{\sigma} = \left\{\,\gamma_1, \dots, \gamma_n \,\right\}\) is a set of homotopy classes of simple closed curves on \(\Sigma\) such that cutting along these curves produces a pants decomposition 
\[
\Sigma \backslash C_{\sigma} \simeq \sqcup_{\nu} \Sigma_{0,3}^{\nu} \sqcup_{\mu} \Sigma_{0,1}^{\mu}
\] 
where the \(\Sigma^{\nu}_{0,3}\) are the `pairs of pants' and the \(\Sigma^{\mu}_{0,1}\) are discs which are used to fill any unwanted punctures; 
\item The \emph{Moore--Seiberg graph} \(\Gamma_{\sigma}\) is a \(3\)-valent graph specifying branch cuts, and is needed to distinguish when a Dehn twist has been applied to \(\Sigma\). 
\end{enumerate}
\end{definition}

We shall now describe a presentation for \(\mathcal{A}(\Sigma)\) which is dependent to a choice of pants decomposition. By Dehn's theorem, a curve \(\gamma\) can be classified uniquely up to homotopy by the \emph{Dehn parameters}
\[
\left\{\, (p_i, q_i) \;\middle| \; i = 1 \dots n \,\right\},
\]
where \(p_i\) is the intersection number between \(\gamma\) and \(\gamma_i \in C_{\sigma}\), and \(q_i\) is the twisting number between \(\gamma\) and \(\gamma_i \in C_{\sigma}\). 

Each curve \(e \in \Gamma_{\sigma}\) which does not end in the boundary of \(\Sigma\) lies in a subspace \(\Sigma_e\) which is homotopic to either \(\Sigma_{0,1}\) or \(\Sigma_{1,1}\): if \(e\) is a loop then \(\Sigma_{e} \simeq \Sigma_{1,1}\), and if it is not then \(\Sigma_{e} \simeq \Sigma_{0,4}\). To \(e\) we assign the curves:
\begin{enumerate}
    \item \(\gamma^e_s := \gamma_e\) is the unique curve \(\gamma_e \in C_{\sigma}\) which lies in the interior of \(\Sigma_e\); it is the curve in the cut system for \(\Sigma\) which also defines a cut system for \(\Sigma_{e}\);
    \item \(\gamma^e_t\) has Dehn parameters \(\left\{\,(p^e_i, 0) \;\middle| \; i = 1, \dots, n\,\right\}\);
    \item \(\gamma^e_u\) has Dehn parameters \(\left\{\,(p^e_i, \delta_{i,e}) \;\middle| \; i = 1, \dots, n\,\right\}\)    
\end{enumerate}
where \(p^e_i := \begin{cases}
2 \delta_{i,e} & \text{ if }  \Sigma_{e} \simeq \Sigma_{0,4} \\
\delta_{i,e} & \text{ if }  \Sigma_{e} \simeq \Sigma_{1,1}.
\end{cases}\)

\begin{definition}
Let \(\gamma\) be a closed curve on \(\Sigma\). Its \emph{geodesic length function} is \(L_{\gamma}:= \nu_{\gamma} \Tr(\rho(\gamma))\) where \(\nu\) is a sign and \(\rho: \pi_1(\Sigma) \to \SL\) is the uniformisation representation.
\end{definition}

\begin{remark}
The geodesic length functions depend only on the homotopy class of the curve, and they satisfy the `skein' relation 
\[
L_{S(\gamma_1, \gamma_2)} = L_{\gamma_1} L_{\gamma_2}
\] 
where \(S(\gamma_1, \gamma_2)\) is a curve with a crossing point and \(\gamma_1, \gamma_2\) are the curves which result from the symmetric smoothing operation:
\[
\smalldiagram{g8014}  \xmapsto{\text{S}} \smalldiagram{g4563} + \smalldiagram{g4571}
\]
\end{remark}

\begin{proposition}{{\cite[Section~2.5.4]{Teschner}}}
The generators of \(\mathcal{A}(\Sigma)\) are 
\[\left\{\, L^e_s, L^e_t, L^e_u \;\middle| \; e \in \Gamma_{\sigma} \text{ is an interior edge}\,\right\}\]
where \(L^e_k = |L_{\gamma^e_k}|\). There is a single relation \(\mathcal{P}_e(L^e_s, L^e_t, L^e_u)\) on \(\mathcal{A}(\Sigma)\) for each internal edge \(e\):
\begin{align*}
\mathcal{P}_e(L^e_s, L^e_t, L^e_u) &= - L^e_s L^e_t L^e_u +(L^e_s)^2 +(L^e_t)^2 +(L^e_u)^2 \\
&+L^e_s (L_3 L_4 + L_1 L_2) + L^e_t (L_2 L_3 + L_1 L_4) + L^e_u(L_1 L_3 + L_2 L_4) \\
&-4 + L^2_1 + L^2_2 + L^2_3 + L^2_4 + L_1 L_2 L_3 L_4 \text{ when } \Sigma_e \simeq \Sigma_{0, 4}, \text{and} \\
\mathcal{P}_e(L^e_s, L^e_t, L^e_u) &= - L^e_s L^e_t L^e_u + (L^e_s)^2 +(L^e_t)^2 +(L^e_u)^2 +L_0 -2 \text{ when } \Sigma_e \simeq \Sigma_{1, 1},
\end{align*}
where \(L_1, L_2, L_2, L_4\) are loops around the four punctures of \(\Sigma_{0,4}\), and \(L_0\) is a loop around the single puncture of \(\Sigma_{1,1}\). The Poisson bracket on \(\mathcal{A}(\Sigma)\) is given by
\[\left\{\, L_{\gamma_1}, L_{\gamma_2} \,\right\} = L_{A(\gamma_1,\gamma_2)},\]
where \(A\) is the antisymmetric smoothing operation:
\[
\smalldiagram{g8014}  \xmapsto{\text{A}} \smalldiagram{g4563} - \smalldiagram{g4571}
\]
\end{proposition}

\begin{figure}[H]
\centering
\includegraphics[height=4cm]{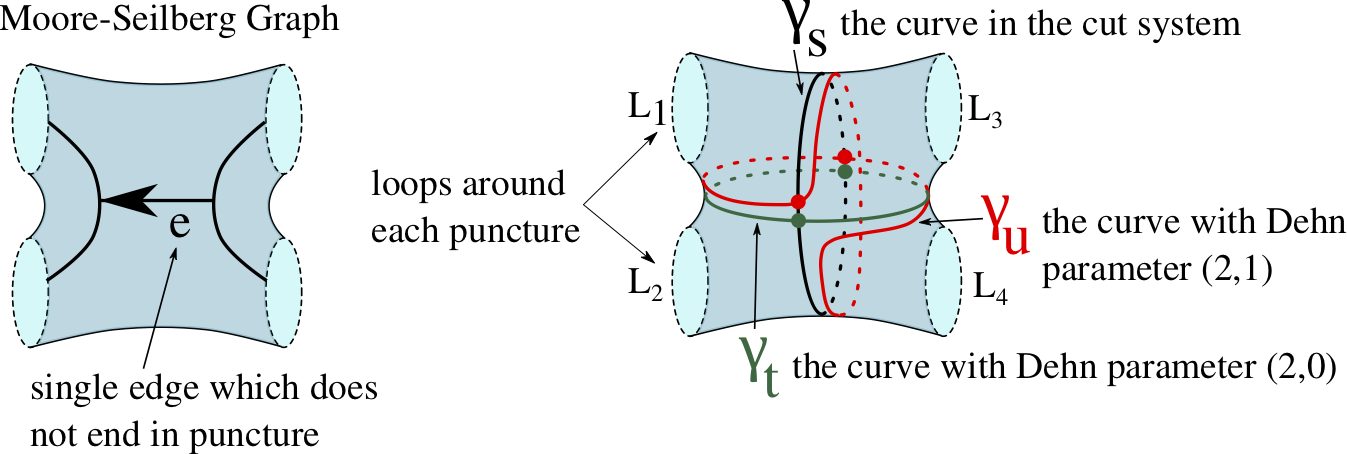}
\caption{Applied to the four-punctured sphere.}
\end{figure}

As \(\mathcal{A}(\Sigma)\) is given by local data on copies of \(\Sigma_{0,4}\) and \(\Sigma_{1,1}\), Teschner and Vartanov state the deformation for these basic surfaces. 

\begin{definition}[{\cite[39--40]{Teschner}}]
The deformation \(\mathcal{A}_b (\Sigma_{0,4})\) of \(\mathscr{A} (\Sigma_{0,4})\) is generated by \(L_s, L_t, L_u, L_1, L_2, L_3, L_4\) with relations \begin{align*}
\mathcal{Q}_e (L_s, L_t, L_u) &=e^{\pi i b^2} L_s L_t - e^{-\pi i b^2} L_t L_s \\
&- (e^{2 \pi i b^2} - e^{-2 \pi i b^2})L_u - (e^{\pi i b^2}-e^{-\pi i b^2})(L_1 L_3 + L_2 L_4) \\
\mathcal{P}_e(L_s, L_t, L_u) &= -e^{\pi i b^2} L_s L_t L_u + e^{2 \pi i b^2} L_u^2 + e^{2 \pi i b^2} L_s^2 + e^{-2 \pi i b^2} L_t^2 \\
    &+ e^{\pi i b^2} (L_1 L_3 + L_2 L_4) L_u + e^{\pi i b^2} (L_3 L_4 + L_2 L_1 ) L_s \\
    &+ e^{-\pi i b^2} (L_1 L_4 + L_2 L_3 ) L_t + L_1^2 + L_3^2 + L_2^2 + L_4^2  + L_1 L_3 L_2 L_4 \\
    &- (2\cos(\pi b^2))^2
\end{align*}
where the quadratic relations \(\mathcal{Q}_e\) arise from deforming the Poisson bracket.
\end{definition}

\begin{definition}[{\cite[40]{Teschner}}]
The deformation \(\mathcal{A}_b (\Sigma_{1,1})\) of \(\mathscr{A} (\Sigma_{1,1})\) is generated by \(L_s, L_t, L_u, L_0\) with relations 
\begin{align*}
\mathcal{Q}_e (L_s, L_t, L_u) &= e^{\frac{\pi i }{2}} L_s L_t -  e^{-\frac{\pi i }{2}} L_t L_s -(e^{\pi i b^2} - e^{-\pi i b^2}) L_u\\
\mathcal{P}_e(L_s, L_t, L_u) &= e^{\pi i b^2} L_s^2 +e^{-\pi i b^2} L_t^2 + e^{\pi i b^2} L_u^2 - e^{\frac{\pi i }{2}} L_s L_t L_u + L_0 - 2\cos(\pi b^2)
\end{align*}
\end{definition}

Using the presentation for the algebras of \(\qgroup{\mathfrak{sl}_2}\)-invariants \(\mathscr{A}_{\Sigma_{0,4}}\) and \(\mathscr{A}_{\Sigma_{1,1}}\) from \cref{section:InvAlgebraSphere}, we see that we have the following isomorphisms:

\begin{proposition}
The algebra of \(\qgroup{\mathfrak{sl}_2}\)-invariants \(\mathscr{A}_{\Sigma_{0,4}}\) is isomorphic to \(\mathcal{A}_b(\Sigma_{0,4})\) with isomorphism \(\iota: \mathscr{A}_{\Sigma_{0,4}} \to \mathcal{A}_b(\Sigma_{0,4})\) defined by
\begin{alignat*}{3}
    \iota(q) &= e^{i \pi b^2}, && \iota(s) &= e^{-i \pi b^2} L_1, \\
    \iota(E) &= -e^{-i \pi b^2} L_u, &\quad& \iota(t) &= e^{-i \pi b^2} L_3, \\
    \iota(F) &= -e^{-i \pi b^2} L_s, &\quad& \iota(v) &= e^{-i \pi b^2} L_2, \\
    \iota(G) &= -e^{-i \pi b^2} L_t, &\quad& \iota(u) &= e^{-i \pi b^2} L_4.
\end{alignat*}
\end{proposition}
\begin{proof}
The map \(\kappa: S \mathscr{H}_{q, \underline{t}} \to \mathcal{A}_b(\Sigma_{0,4})\) defined by
\begin{align*}
    q &\mapsto e^{i \pi b^2}, & \overline{t_1} &\mapsto i L_1, \\
    x &\mapsto L_u, & \overline{t_2} &\mapsto i L_3, \\
    y &\mapsto L_s, & \overline{q t_3} &\mapsto i L_2, \\
    z &\mapsto L_t, & \overline{t_4} & \mapsto i L_4, \\
\end{align*}
maps \(S \mathscr{H}_{q, \underline{t}}\) to an algebra generated by \(L_s, L_t, L_u\) with relations 
\begin{align*}
    0 &=e^{\pi i b^2} L_u L_s -e^{-\pi i b^2} L_s L_u - (e^{2 \pi i b^2} - e^{-2 \pi i b^2})L_t - (e^{\pi i b^2}-e^{-\pi i b^2})(L_1 L_4 + L_2 L_3 ) \\
    0 &=e^{\pi i b^2} L_s L_t - e^{-\pi i b^2} L_t L_s - (e^{2 \pi i b^2} - e^{-2 \pi i b^2})L_u - (e^{\pi i b^2}-e^{-\pi i b^2})(L_1 L_3 + L_2 L_4) \\
    0 &= e^{\pi i b^2} L_t L_u - e^{-\pi i b^2} L_u L_t - (e^{2 \pi i b^2} - e^{-2 \pi i b^2})L_s - (e^{\pi i b^2}-e^{-\pi i b^2})(L_3 L_4 + L_2 L_1 ) \\
    0 &= -e^{\pi i b^2} L_s L_t L_u + e^{2 \pi i b^2} L_u^2 + e^{2 \pi i b^2} L_s^2 + e^{-2 \pi i b^2} L_t^2 \\
    &+ e^{\pi i b^2} (L_1 L_3 + L_2 L_4) L_u + e^{\pi i b^2} (L_3 L_4 + L_2 L_1 ) L_s + e^{-\pi i b^2} (L_1 L_4 + L_2 L_3 ) L_t \\
    &+ L_1^2 + L_3^2 + L_2^2 + L_4^2 + L_1 L_3 L_2 L_4 - (2\cos(\pi b^2))^2
\end{align*}
which is just the algebra \(\mathcal{A}_b(\Sigma_{0,4})\). Hence the algebra \(\mathscr{A}_{\Sigma_{0,4}}\) is isomorphic to both \(S \mathscr{H}_{q, \underline{t}}\) and \(\mathcal{A}_b(\Sigma_{0,4})\) and isomorphism \(\iota: \mathscr{A}_{\Sigma_{0,4}} \to \mathcal{A}_b(\Sigma_{0,4})\) is given by \(\kappa \circ \alpha^{-1}\).
\end{proof}

\begin{proposition}
The algebra of \(\qgroup{\mathfrak{sl}_2}\)-invariants \(\mathscr{A}_{\Sigma_{1,1}}\) is isomorphic to \(\mathcal{A}_b(\Sigma_{1,1})\) with isomorphism \(\mu: \mathscr{A}_{\Sigma_{1,1}} \to \mathcal{A}_b(\Sigma_{1,1})\) defined by
\begin{align*}
\mu(Y) &= i e^{-i \pi b^2}L_s \\
\mu(X) &= ie^{-i \pi b^2}L_t \\
\mu(Z) &= -e^{- \frac{5}{2}i \pi b^2} L_u \\
\mu(L) &= L_0 
\end{align*}
\begin{proof}
Follows from \cref{thm:torusiso}.
\end{proof}
\end{proposition}

\appendix
\addtocontents{toc}{\protect\setcounter{tocdepth}{1}} 
\section{Hilbert Series of the Algebras of Invariants}
\label{section:HilbertSeries}

In this section we shall compute the graded character of the algebra objects \(A_{\Sigma_{0,4}}\) and \(A_{\Sigma_{1,1}}\), and then use these to compute the Hilbert series of the algebras of invariants \(\mathscr{A}_{\Sigma_{0,4}}\) and \(\mathscr{A}_{\Sigma_{1,1}}\) which we will need in the proof of presentation of \(\mathscr{A}_{\Sigma_{0,4}}\) and \(\mathscr{A}_{\Sigma_{1,1}}\). A Hilbert series encodes the dimensions of the graded parts of an algebra. 

\begin{definition}
The \emph{associated graded algebra} of the \(\mathbb{Z}_+\) filtered algebra \(A = \bigcup_{n \in \mathbb{Z}_{+}} A(n)\) is 
\vspace{-1em}
\[
\mathscr{G}(A) = \bigoplus_{n \in \mathbb{Z}_{+}} A[n] \text{ where } A[n] = \begin{cases}
A(0) & \text{ for }n = 0 \\
\faktor{A(n)}{A(n-1)} &\text{ for } n > 0.
\end{cases}
\]
\end{definition}

\begin{definition}
The \emph{Hilbert series} of the \(\mathbb{Z}_+\) graded vector space \(A = \bigoplus_{n \in \mathbb{Z}_{+}} A[n]\) is the formal power series\[h_A(t) = \sum \dim(A[n]) t^n.\]The  Hilbert series of a \(\mathbb{Z}_+\) graded algebra \(A\) is the Hilbert series of its underlying \(\mathbb{Z}_+\) graded vector space, and
the Hilbert series of the \(\mathbb{Z}_+\) filtered algebra \(A = \bigcup_{n \in \mathbb{Z}_{+}} A(n)\) is the Hilbert series of the associated graded algebra \(\mathscr{G}(A)\).
\end{definition}

A graded character of a filtered/graded representation encodes the dimensions of graded parts and weight spaces simultaneously.

\begin{definition}
Let \(V\) be a vector space acted on by \(\qgroup{\frsl}\) and let \(V^k\) denote the \(q^k\)-weight space of \(V\) where \(k \in \mathbb{Z}\). The \emph{character} of \(V\) is the formal power series
\[
\ch_V(u) = \sum_{k \in \Lambda} \dim \left(V^{k} \right) u^{k}.
\]
\end{definition}

\begin{definition}
Let \(V = \bigoplus_{n} V[n]\) be a graded vector space acted on by \(\qgroup{\frsl}\). The \emph{graded character} of \(V\) is  
\[
h_V(u,t) := \sum_n \ch_{V[n]}(u) t^n = \sum_{n,k} \dim \left(V[n]^k\right) u^k t^n,
\]
where \(V[n]^k\) is the \(q^k\)-weight space of \(V[n]\). If \(V\) is filtered rather than graded the graded character of \(V\) \(h_V(u,t)\) is \(h_{\mathscr{G}(V)}(u,t)\), the graded character of associated graded vector space \(\mathscr{G}(V)\).
\end{definition}

Let \(\Sigma = \Sigma_{0,4}\) or \(\Sigma_{1,1}\). Both \(A_{\Sigma}\) and its subalgebra \(\mathscr{A}_{\Sigma}\) have filtrations by degree:
\[
A_{\Sigma} = \bigcup_{n \in \mathbb{Z}_+} A(n); \; \mathscr{A}_{\Sigma} = \bigcup_{n \in \mathbb{Z}_+} \mathscr{A}(n)
\]
where \(A(n)\) and \(\mathscr{A}(n)\) are the span of monomials in \(A_{\Sigma}\) and \(\mathscr{A}_{\Sigma}\) respectively with at most \(n\) generators.

\begin{remark}
Unless otherwise stated, Hilbert series will always assume grading by degree, and the action of \(\qgroup{\frsl}\) will always be that stated in Example \ref{Kaction}.
\end{remark}

As \(\mathscr{A}_{\Sigma}\) is the part of \(A_{\Sigma}\) with weight \(1=q^0\) under the action of \(\qgroup{\frsl}\), the terms of the graded character \(h_{A_{\Sigma}}(u,v)\) where \(k=0\) give the Hilbert series \(h_{\mathscr{A}_{\Sigma}}(t)\); hence, we shall:
\begin{stages}
\item Compute the graded character of \(\mathscr{O}_q(\SL)\) which we use to
\item Compute the graded character of \(A_{\Sigma}\), and then
\item Extract the terms of the graded character which give the Hilbert series of \(\mathscr{A}_{\Sigma}\).
\end{stages}

\subsubsection{The Graded Character of the Algebra Objects \(A_{\Sigma_{0,4}}\) and \(A_{\Sigma_{1,1}}\)}
\begin{proposition}
The graded character of \(\mathscr{O}_q(\SL)\) is 
\begin{align*}
    h_{\mathscr{O}_q}(u,t) = \frac{(1 + t)}{(1 - t) (1-u^2t) (1-u^{-2}t)}.
\end{align*}
\end{proposition}
\begin{proof}
Recall from Proposition \ref{pbwbasisorignial} that \(\mathscr{O}_q(\Rep_q(\SL))\) has basis
\[
\left\{\,(a^1_1)^{\alpha} (a^1_2)^{\beta} (a^2_1)^{\gamma} (a^2_2)^{\delta} \; \middle| \; \alpha, \beta, \gamma, \delta \in \mathbb{N}_0;\; \beta \text{ or }\gamma  = 0\,\right\}.
\]
We shall denote \(X_{\alpha, \beta, \gamma, \delta}:= (a^1_1)^{\alpha} (a^1_2)^{\beta} (a^2_1)^{\gamma} (a^2_2)^{\delta}\). The \(n^{th}\) graded part \(\mathscr{O}_q[n]:=\left(\mathscr{O}_q(\Rep_q(\SL))\right)[n]\) has basis 
\[
\left\{\, X_{\alpha, \beta, \gamma, \delta}  \; \middle| \; \alpha, \beta, \gamma, \delta \in \mathbb{N}_0;\; \beta \text{ or }\gamma  = 0;\; \alpha + \beta + \gamma + \delta = n\,\right\}.
\]
We can see from Example \ref{Kaction} that \(a^1_1, a^1_2, a^2_1, a^2_2\) have weights \(1, q^2, q^{-2}, 1\) respectively, so
\[
K \cdot X_{\alpha, \beta, \gamma, \delta}= K \cdot \left((a^1_1)^{\alpha} (a^1_2)^{\beta} (a^2_1)^{\gamma} (a^2_2)^{\delta} \right) = q^{2 \beta - 2\gamma} (a^1_1)^{\alpha} (a^1_2)^{\beta} (a^2_1)^{\gamma} (a^2_2)^{\delta} = q^{2(\beta-\gamma)} X_{\alpha, \beta, \gamma, \delta},
\]
and \(X_{\alpha,\beta,\gamma,\delta}\) has weight \(q^{2(\beta-\gamma)}\). This means that \(\mathscr{O}_q[n]^k\), the \(q^k\) weight space of \(\mathscr{O}_q[n]\), has basis
\[
\left\{\,X_{\alpha,\beta,\gamma,\delta}  \; \middle| \; \alpha, \beta, \gamma, \delta \in \mathbb{N}_0;\; \beta \text{ or }\gamma  = 0;\; \alpha + \beta + \gamma + \delta = n ;\; 2(\beta-\gamma) = k\,\right\}.
\]
If \(k\) is odd the final condition is never satisfied, and thus \(\mathscr{O}_q[n]^k = \emptyset\). If \(k=2m\) for \(m \geq 0\) then we get the basis 
\begin{align*}
&\left\{\,X_{\alpha,\beta,\gamma,\delta}  \; \middle| \; \alpha, \beta, \gamma, \delta \in \mathbb{N}_0;\; \beta \text{ or }\gamma  = 0;\; \alpha + \beta + \gamma + \delta = n;\; 2(\beta-\gamma) = 2m\,\right\}\\
&= \left\{\,X_{\alpha,\beta, 0, \delta}  \; \middle| \; \alpha, \beta, \gamma, \delta \in \mathbb{N}_0;\; \alpha + \beta + \delta = n;\; \beta = m \,\right\}\\
&\text{ as } \beta - \gamma \geq 0 \text{ and } \beta \text{ or } \gamma = 0 \text{ implies } \gamma = 0 \\
&= \left\{\,X_{\alpha,m,0,\delta}  \; \middle| \; \alpha, \delta \in \mathbb{N}_0;\; \alpha + \delta = n-m\,\right\}.
\end{align*}
which is empty if \(m>n\) and has \(n-m+1\) elements otherwise. Finally, if \(k=-2m\) for \(m>0\) then we get the basis 
\begin{align*}
&\left\{\,X_{\alpha,\beta,\gamma,\delta}  \; \middle| \; \alpha, \beta, \gamma, \delta \in \mathbb{N}_0;\; \beta \text{ or }\gamma  = 0;\; \alpha + \beta + \gamma + \delta = n;\; 2(\beta-\gamma) = -2m\,\right\}\\
&= \left\{\,X_{\alpha,0, \gamma,\delta}  \; \middle| \; \alpha, \beta, \gamma, \delta \in \mathbb{N}_0;\; \alpha + \gamma + \delta = n;\; \gamma = m\,\right\}\\
&\text{ as } \beta - \gamma \leq 0 \text{ and } \beta \text{ or } \gamma = 0 \text{ implies } \beta = 0 \\
&= \left\{\,X_{\alpha,0,m,\delta}  \; \middle| \; \alpha, \delta \in \mathbb{N}_0;\; \alpha + \delta = n-m\,\right\}.
\end{align*}
which is empty if \(m>n\) and has \(n-m+1\) elements otherwise. Hence,
\[
\dim \mathscr{O}_q[n]^k = \begin{cases}
n-m+1 &\text{ if } k=2m \text{ for some } m\geq 0 \\
n-m+1 &\text{ if } k=-2m \text{ for some } m\geq 1 \\
0 &\text{ otherwise},
\end{cases}
\]
so the character of \(\mathscr{O}_q[n]\) is 
\begin{align*}
h_{\mathscr{O}_q[n]}(u) &= \left(\sum_{m=0}^n (n-m+1)u^{2m}\right) + \left(\sum_{m=1}^n (n-m+1)u^{-2m}\right) \\
&= \frac{u^{-2 n} (u^{2 + 2 n}-1)^2}{(u^2-1)^2},
\end{align*}
and the graded character of \(\mathscr{O}_q\) is
\begin{align*}
    h_{\mathscr{O}_q}(u,t) = \sum_{n=0}^{\infty} \frac{u^{-2 n} (u^{2 + 2 n}-1)^2}{(u^2-1)^2} t^n = \frac{(1 + t)}{(1 - t) (1-u^2t) (1-u^{-2}t)}.
\end{align*}
\end{proof}

We note that if \(V = \bigoplus_{n} V(n)\) and \(W = \bigoplus_{n} W(n)\) are two graded vector spaces acted on by \(\qgroup{\frsl}\) then \(h_{V \otimes W}(u,t) = h_V(u,t) \cdot h_W (u,t)\).

\begin{corollary}
\label{cor:charactecalc}
The graded character of \(A_{\Sigma_{0,4}}\) is 
\[
h_{A_{\Sigma_{0,4}}}(u,t) = \left(\frac{(1 + t)}{(1 - t) (1-u^2t) (1-u^{-2}t)}\right)^3.
\]
\end{corollary}
\begin{proof}
We have from Proposition \ref{factofsurface} that \(A_{\Sigma_{0,4}} \cong \mathscr{O}_q \otimes \mathscr{O}_q \otimes \mathscr{O}_q\); hence,
\begin{align*}
    h_{A_{\Sigma_{0,4}}}(u,t) = h_{\mathscr{O}_q}(u,t) \cdot h_{\mathscr{O}_q}(u,t) \cdot h_{\mathscr{O}_q}(u,t)  = \left(\frac{(1 + t)}{(1 - t) (1-u^2t) (1-u^{-2}t)}\right)^3.
\end{align*}

\end{proof}

\begin{corollary}
\label{cor:charactecalctorus}
The graded character of \(A_{\Sigma_{1,1}}\) is 
\[
h_{A_{\Sigma_{1,1}}}(u,t) = \left(\frac{(1 + t)}{(1 - t) (1-u^2t) (1-u^{-2}t)}\right)^2.
\]
\end{corollary}
\begin{proof}
We have from Proposition \ref{factofsurface} that \(A_{\Sigma_{1,1}} \cong \mathscr{O}_q \otimes \mathscr{O}_q\); hence,
\begin{align*}
    h_{A_{\Sigma_{1,1}}}(u,t) = h_{\mathscr{O}_q}(u,t) \cdot h_{\mathscr{O}_q}(u,t) = \left(\frac{(1 + t)}{(1 - t) (1-u^2t) (1-u^{-2}t)}\right)^2.
\end{align*}

\end{proof}

\subsubsection{The Hilbert Series of \(\mathscr{A}_{\Sigma_{0,4}}\) and \(\mathscr{A}_{\Sigma_{1,1}}\)}

\begin{proposition}
\label{prop:characterfromweyl}
Let \(\Sigma\) be any punctured surface and \(A_{\Sigma}\) be the algebra object of \(\int_{\Sigma} \Rep_q(\SL)\). The graded character of \(A_{\Sigma}\) is
\[h_{A_{\Sigma}}(u,t) = \sum_{n,k} m_{n,k} \frac{u^{k+1} - u^{-k-1}}{u - u^{-1}} t^n\]
for \(m_{n,k} \in \mathbb{Z}_+\).
\begin{proof}
As integrable representations of \(\qgroup{\frsl}\) are semisimple, any finite-dimensional representation \(V\) of \(\qgroup{\frsl}\) when \(q\) is generic can be decomposed into \(V = \bigoplus_{k \in \mathbb{Z}_+} V[k]^{m_k}\) where \(m_k \in \mathbb{Z}_+\) and \(V[k]\) is an irreducible representation with character given by the Weyl character formula:
\[\ch_{V(k)} = u^k + u^{k-2} + \dots + u^{-k+2} + u^{-k} = \frac{u^{k+1} - u^{-k-1}}{u - u^{-1}}.\]
Applying this to \(V = A_{\Sigma}[n]\) the degree \(n\) part of \(\mathscr{G}(A_{\Sigma})\) gives 
\begin{align*}
    h_{A_{\Sigma}}(u,t) &= h_{\mathscr{G}(A_{\Sigma})}(u,t) \\
        &= \sum_n \ch_{V[n]}(u) t^n \\
        &= \sum_n \ch_{\bigoplus_k V[n](k)^{m_{n,k}}}(u) t^n  \\
        &=\sum_{n,k} m_{n,k} \ch_{V[n](k)}(u) t^n \\
        &= \sum_{n,k} m_{n,k} \frac{u^{k+1} - u^{-k-1}}{u - u^{-1}}t^n.
\end{align*}
\end{proof}
\end{proposition}

\begin{corollary}
\label{cor:ucoeff}
Let \(A_{\Sigma}\) be the algebra object and \(\mathscr{A}_{\Sigma}\) be the algebra of \(\qgroup{\mathfrak{sl}_2}\)-invariants of the factorisation homology of \(\int_{\Sigma} \Rep_q(\SL)\) for a punctured surface \(\Sigma\). The Hilbert series \(h_{\mathscr{A}_{\Sigma}}(t)\) is given by the \(u\) coefficient of \((u-u^{-1})\cdot h_{A_{\Sigma}}(u,t)\).
\begin{proof}
From Proposition \ref{prop:characterfromweyl} we have that
\begin{align*}
h_{A_{\Sigma}}(u,t) &= \sum_{n,k} m_{n,k} \frac{u^{k+1} - u^{k-1}}{u - u^{-1}} t^n \\
\implies (u-u^{-1}) h_{A_{\Sigma}}(u,t) &= \sum_{n,k} m_{n,k} (u^{k+1} - u^{k-1}) t^n
\end{align*}
where 
\[h_{\mathscr{A}_{\Sigma}}(t) = \sum_{n} m_{n,0} t^n,\]
so \(h_{\mathscr{A}_{\Sigma}}(t)\) is given by the \(u\) coefficient of \((u-u^{-1})\cdot h_{A_{\Sigma}}(u,t)\).
\end{proof}
\end{corollary}

\begin{proposition}
The Hilbert series of \(\mathscr{A}_{\Sigma_{0,4}}\) is 
\[h_{\mathscr{A}_{\Sigma_{0,4}}}(t) =  \frac{t^2 - t+ 1}{(1-t)^6(1+t)^2}.\]
\end{proposition}
\begin{proof}
From Corollary \ref{cor:charactecalc} we have that
\begin{align*}
h_{A_{\Sigma_{0,4}}}(u,t) &= \left(\frac{(1 + t)}{(1 - t) (1-u^2t) (1-u^{-2}t)}\right)^3 \\
&= \frac{1}{(1-t)^6}\Bigg(\frac{ t^3}{(u^2 - t)^3} +\frac{3 t^2}{(1-t^2) (u^2 - t)^2} \\
&+\frac{3 (t^2 + 1) t}{(1-t^2)^2 (u^2 - t)} + \frac{1}{(1- t u^2)^3} \\
&+ \frac{3 t^2}{(1-t^2) (1-t u^2)^2}
+ \frac{3 t^2 (t^2 + 1)}{(1-t^2)^2 (1-t u^2)} \Bigg)
\end{align*}
where
\begin{align*}
\frac{1}{(1-u^2t)} &= \sum_{i=0}^{\infty} (u^2 t)^i = 1 +u^2t +u^4t^2 +\dots \\
\frac{1}{(u^2 -t)} &= u^{-2} \sum_{i=0}^{\infty} (u^{-2}t)^i = u^{-2} + u^{-4}t + \dots \\
\end{align*}
so the \(u\) coefficient of \((u-u^{-1})\cdot h_{A_{\Sigma_{0,4}}}(u,t)\) is
\[
    \frac{1}{(1-t)^6} \Bigg( (1 - 3t) +\frac{3 t^2 (1 - 2t)}{(1-t^2)}  + \frac{3t^2 (1 - t) (t^2+1)}{(1-t^2)^2} \Bigg) = \frac{t^2 - t+ 1}{(1-t)^6(1+t)^2}
\]
which by Corollary \ref{cor:ucoeff} is the Hilbert series of \(\mathscr{A}_{\Sigma_{0,4}}\). 
\end{proof}

\begin{proposition}
The Hilbert series of \(\mathscr{A}_{1,1}\) is 
\[h_{\mathscr{A}_{\Sigma_{1,1}}} = \frac{1}{(1-t)^3 (1+t)}.\]
\begin{proof}
From Corollary \ref{cor:charactecalctorus} we have that
\begin{align*}
h_{A_{\Sigma_{1,1}}}(u,t) &= \left(\frac{(1 + t)}{(1 - t) (1-u^2t) (1-u^{-2}t)}\right)^2 \\
&= \frac{(1+t)^2}{(1-t)^2(1-t^2)^2} \bigg( \frac{2t^2}{(1- t^2) (1-tu^2)} + \frac{t^2}{(u^2 - t)^2} \\
&+ \frac{2t}{(1-t^2) (u^2 - t)}  + \frac{1}{(1-tu^2)^2 }\bigg),
\end{align*}
so the \(u\) coefficient of \((u - u^{-1})h_{A_{\Sigma_{1,1}}}(u,t)\) is 
\[ \frac{(1+t)^2}{(1-t)^2(1-t^2)^2} \left( \frac{2t^2 (1-t)}{(1-t^2)} + (1-2t) \right) = \frac{1}{(1-t)^3 (1+t)}\]
which by Corollary \ref{cor:ucoeff} is the Hilbert series of \(\mathscr{A}_{\Sigma_{1,1}}\). 
\end{proof}
\end{proposition}
uii\section{PBW Basis of \(\mathscr{G}(\mathscr{B})\)}
Recall from \cref{defn:Apres} the definition of \(\mathscr{B}\).
As the elements \(u, v, s\) and \(t\) are central, instead of considering \(\mathscr{B}\) as an algebra over \(\mathbb{C}\) with seven generators, we can consider \(\mathscr{B}\) as an algebra over the polynomial ring \(\mathbb{C}[s, t, u, v]\) with generators \(E, F, G\), i.e.\ \(\mathscr{B} = \mathbb{C}[s, t, u, v] \langle E, F, G \rangle\)\footnote{The algebra \(\langle E, F, G \rangle\) denotes the subalgebra of \(\mathscr{B}\) generated by \(E\), \(F\) and \(G\) not the free algebra.}.

\begin{proposition}
\label{prop:PBWbasissphere}
A PBW-basis for the associated graded algebra \(\mathscr{G}(\mathscr{B})\) over \(\mathbb{C}[s, t, u, v]\) is  \[\left\{ \, E^n F^m G^l \; \middle| \; n \text{ or } m\text{ or }l =0 \, \, \right\}.\]
\begin{proof}

A term rewriting system for \(\mathscr{G}(B)\) is given by 
\begin{align*}
\sigma_{FE}&: FE \mapsto q^2EF + dG + ea\\
\sigma_{GF}&: GF \mapsto q^2 FG + dE + ec \\
\sigma_{GE}&: GE \mapsto q^{-2}EG - q^{-2}dF + fb \\
\sigma_{EF^nG}&: EF^nG \mapsto f(n)
\end{align*}
where 
\[a:= sv + tu,\; b:= su + tv,\; c:= st + uv,\; d := (q^2-q^{-2}),\; e:= (1-q^2),\; f:= (1-q^{-2})\]
and \(f(n)\) is defined recursively as follows\footnote{This recursion relation arises from applying \(\sigma_{FE}^{-1}\) to \(EF^nG\); one could equally apply \(\sigma_{GF}^{-1}\) which would give an alternate term rewriting system.}:
\begin{align*}
    f(1) &:=  -E^2 -q^{-4}F^2 - G^2 +cE + q^{-2}bF + aG \\
    &\quad+ \left(- q^{-4}(s^2 + t^2 + u^2 + v^2) - stuv + q^{-6}(q^2 +1)^2\right) \\
    f(n) &:=  q^{-2} Ff(n-1) + (q^{-4}-1)GF^{n-1}G + (1-q^{-2})aF^{n-1}G.
\end{align*}

We shall use the above term rewriting system for \(\mathscr{G}(B)\) and apply the diamond lemma. In order to do this we must first show that all the ambiguities of the term rewriting system are resolvable. The ambiguities are 
\begin{enumerate}
    \item \((\sigma_{GF}, \sigma_{FE}, G, F, E)\),
    \item \( (\sigma_{FE}, \sigma_{EF^nG}, F, E, F^nG) \),
    \item \((\sigma_{GE}, \sigma_{EF^nG}, G, E, F^nG) \),
    \item \( (\sigma_{EF^nG},\sigma_{GE}, EF^n, G, E) \),
    \item \( (\sigma_{EF^nG}, \sigma_{GF}, EF^n, G, F)\).
\end{enumerate}

The first ambiguity \((\sigma_{GF}, \sigma_{FE}, G, F, E)\) is resolvable by direct calculation:
\begin{align*}
    GFE &\xmapsto{\sigma_{GF}} q^2 FGE + dE^2 + ecE \\
    &\xmapsto{\sigma_{GE}} FEG - dF^2 + q^2fbF+ dE^2 + ecE \\
    &\xmapsto{\sigma_{FE}} q^2EFG + dG^2 + eaG - dF^2 + q^2fbF+ dE^2 + ecE
\end{align*}
is equal to
\begin{align*}
    GFE &\xmapsto{\sigma_{FE}} q^2GEF + dG^2 + eaG \\
    &\xmapsto{\sigma_{GE}} EGF - dF^2 + q^2fbF + dG^2 + eaG \\
    &\xmapsto{\sigma_{GF}} q^2 EFG + dE^2 + ecE - dF^2 + q^2fbF + dG^2 + eaG.
\end{align*}

The second ambiguity \((\sigma_{FE}, \sigma_{EF^nG}, F, E, F^nG)\) also follows directly:
\begin{align*}
    FEF^nG &\xmapsto{\sigma_{FE}} q^2EF^{n+1}G + dGF^nG + eaF^nG \\
    &\xmapsto{\sigma_{EF^{n+1}G}} Ff(n) -dGF^{n}G + (q^2-1)aF^{n}G + dGF^nG + eaF^nG\\
    &= Ff(n)
\end{align*}
is equal to 
\begin{align*}   
FEF^nG &\xmapsto{\sigma_{EF^nG}} Ff(n).
\end{align*}

For the remainder of the ambiguities we proceed by induction on \(n\). For the third ambiguity
\((\sigma_{GE}, \sigma_{EF^nG}, G, E, F^nG)\) one direction is given by:
\begin{align}
GEF^nG &\xmapsto{\sigma_{GE}} q^{-2}EGF^nG - q^{-2}dF^{n+1}G + fbF^nG \nonumber\\
&\xmapsto{\sigma_{GF}} EFGF^{n-1}G + q^{-2}dE^2F^{n-1}G + q^{-2}ecEF^{n-1}G \nonumber\\
&\quad- q^{-2}dF^{n+1}G + fbF^nG \nonumber\\
&\xmapsto{\sigma_{EFG}} \Big(-E^2 -q^{-4}F^2 - G^2 +cE + q^{-2}bF + aG \nonumber\\
&\quad- q^{-4}(s^2 + t^2 + u^2 + v^2) - stuv + q^{-6}(q^2 +1)^2\Big)F^{n-1}G \nonumber\\
&\quad+ (1-q^{-4})E^2F^{n-1}G + (q^{-2}-1)cEF^{n-1}G - q^{-2}dF^{n+1}G + fbF^nG \nonumber\\
&= \Big(-q^{-4}E^2 -F^2 - G^2 +q^{-2}cE + bF + aG \nonumber\\
&\quad- q^{-4}(s^2 + t^2 + u^2 + v^2) - stuv + q^{-6}(q^2 +1)^2\Big)F^{n-1}G \; \text{ for all } n \geq 1 \tag{\(\dagger\)}\\
&\xmapsto{\sigma_{EF^{n-1}G}^2} \Big(-F^2 - G^2 + bF + aG - q^{-4}(s^2 + t^2 + u^2 + v^2) - stuv  \nonumber\\
&\quad+ q^{-6}(q^2 +1)^2\Big)F^{n-1}G - q^{-4}Ef(n-1) + q^{-2}cf(n-1) \text{ when } n \neq 1. \tag{\(\ddagger\)} 
\end{align}
This equals the other direction when \(n=1\):
\begin{align*}
GEFG &\xmapsto{\sigma_{EFG}}  -GE^2 -q^{-4}GF^2 - G^3 +cGE + q^{-2}bGF + aG^2 \\
&\quad- q^{-4}(s^2 + t^2 + u^2 + v^2)G - stuvG + q^{-6}(q^2 +1)^2G \\
&\xmapsto{\sigma_{GE}^3} -q^{-4}E^2G + q^{-4}dEF - q^{-2}fbE + q^{-2}dFE - fbE -q^{-4}GF^2 - G^3 \\
&\quad+q^{-2}cEG - q^{-2}dcF + fbc + q^{-2}bGF + aG^2 \\
&\quad+\left(- q^{-4}(s^2 + t^2 + u^2 + v^2) - stuvG + q^{-6}(q^2 +1)^2\right)G \\
&\xmapsto{\sigma_{GF}^3} -q^{-4}E^2G + q^{-4}dEF - q^{-2}fbE + q^{-2}dFE - fbE \\
&\quad- F^2G -q^{-2} dFE -q^{-2} ecF -q^{-4} dEF -q^{-4} ecF - G^3 \\
&\quad+q^{-2}cEG - q^{-2}dcF + fbc + bFG + q^{-2}dbE + q^{-2}ebc + aG^2 \\
&\quad+\left(- q^{-4}(s^2 + t^2 + u^2 + v^2) - stuvG + q^{-6}(q^2 +1)^2\right)G \\
&=\Big(-q^{-4}E^2 - F^2 - G^2 + q^{-2}cE + bF   + aG \\
&\quad- q^{-4}(s^2 + t^2 + u^2 + v^2) - stuvG + q^{-6}(q^2 +1)^2\Big)G \\
&= (\dagger)
\end{align*}
And in the general case:
\begin{align*}
GEF^nG &\xmapsto{EF^nG} q^{-2} GFf(n-1) + (q^{-4}-1)G^2F^{n-1}G + (1-q^{-2})aGF^{n-1}G \\
        &\xmapsto{\sigma_{GF}} FGf(n-1) + q^{-2}dEf(n-1) + q^{-2}ecf(n-1) \\
        &\quad+ (q^{-4}-1)G^2F^{n-1}G + (1-q^{-2})aGF^{n-1}G \\
        &\xmapsto{} q^{-2}FEGF^{n-1}G - q^{-2}dF^{n+1}G + fbF^{n}G \\
        &\quad+ q^{-2}dEf(n-1) + q^{-2}ecf(n-1) + (q^{-4}-1)G^2F^{n-1}G \\
        &\quad+ (1-q^{-2})aGF^{n-1}G \text{ by the induction assumption}\\
        &\xmapsto{\sigma_{FE}} EFGF^{n-1}G + q^{-2}dG^2F^{n-1}G + q^{-2}eaGF^{n-1}G - q^{-2}dF^{n+1}G \\
        &\quad+ fbF^{n}G + q^{-2}dEf(n-1) + q^{-2}ecf(n-1) + (q^{-4}-1)G^2F^{n-1}G \\
        &\quad+ (1-q^{-2})aGF^{n-1}G \\
        &\xmapsto{\sigma_{EFG}} \Big( -E^2 -q^{-4}F^2 - G^2 +cE + q^{-2}bF + aG \\
    &\quad+ \left(- q^{-4}(s^2 + t^2 + u^2 + v^2) - stuv + q^{-6}(q^2 +1)^2\right) \Big)F^{n-1}G \\
    &\quad+ q^{-2}dG^2F^{n-1}G + q^{-2}eaGF^{n-1}G - q^{-2}dF^{n+1}G + fbF^{n}G \\
    &\quad+ q^{-2}dEf(n-1) + q^{-2}ecf(n-1) + (q^{-4}-1)G^2F^{n-1}G\\
    &\quad+ (1-q^{-2})aGF^{n-1}G \\
    &= \Big( -E^2 -F^2 - G^2 +cE + bF + aG + \Big(- q^{-4}(s^2 + t^2 + u^2 + v^2) \\
    &\quad- stuv + q^{-6}(q^2 +1)^2\Big) \Big)F^{n-1}G + q^{-2}dEf(n-1) + q^{-2}ecf(n-1) \\
    &\xmapsto{\sigma_{EF^{n-1}G}^2} \Big(-F^2 - G^2 + bF + aG + \Big(- q^{-4}(s^2 + t^2 + u^2 + v^2) \\
    &\quad- stuv + q^{-6}(q^2 +1)^2\Big) \Big)F^{n-1}G - q^{-4}Ef(n-1) +q^{-2}cf(n-1) \\
    &= (\ddagger)
\end{align*}

For the fourth ambiguity \((\sigma_{EF^nG},\sigma_{GE}, EF^n, G, E)\), one direction is:
\begin{align*}
    EF^nGE &\xmapsto{\sigma_{GE}}  q^{-2}EF^nEG - q^{-2}dEF^{n+1} + fbEF^n \\
    &\xmapsto{\sigma_{FE}} EF^{n-1}(EFG + q^{-2}dG^2 + q^{-2}eaG - q^{-2}dF^2 + fbF) \\
    &\xmapsto{\sigma_{EFG}}EF^{n-1}\Big(-E^2 -F^2 - q^{-4}G^2 +cE + bF + q^{-2}aG \\
    &\quad- q^{-4}(s^2 + t^2 + u^2 + v^2) - stuv + q^{-6}(q^2 +1)^2\Big).
\end{align*}
This equals the other direction when \(n=1\):
\begin{align*}
    EFGE &\xmapsto{\sigma_{EFG}} -E^3 -q^{-4}F^2E - G^2E +cE^2 + q^{-2}bFE + aGE \\
    &\quad+\left(- q^{-4}(s^2 + t^2 + u^2 + v^2) - stuv + q^{-6}(q^2 +1)^2\right)E \\
    &= E\left(-E^2 +cE - q^{-4}(s^2 + t^2 + u^2 + v^2) - stuv + q^{-6}(q^2 +1)^2\right)E \\
    &\quad-q^{-4}F^2E - G^2E + q^{-2}bFE + aGE \\
    &\xmapsto{\sigma_{FE}^3 \circ \sigma_{GE}^3}E\left(-E^2 +cE - q^{-4}(s^2 + t^2 + u^2 + v^2) - stuv + q^{-6}(q^2 +1)^2\right) \\
    &\quad-EF^2 -q^{-2} dGF -q^{-2}eaF - q^{-4}dFG -q^{-4} eaF -q^{-4}EG^2 + q^{-4}dFG \\
    &\quad-q^{-2} fbG + q^{-2}dGF - fbG + bEF + q^{-2}dbG + q^{-2}eab + q^{-2}aEG \\
    &\quad- q^{-2}daF + fab \\
    &= E\Big(-E^2 -F^2 -q^{-4}G^2 + cE + bF + q^{-2}aG - q^{-4}(s^2 + t^2 + u^2 + v^2) \\
    &\quad- stuv + q^{-6}(q^2 +1)^2\Big).
\end{align*}
And in the general case:
\begin{align*}
    EF^nGE &\xmapsto{\sigma_{EF^nG}} q^{-2} Ff(n-1)E + (q^{-4}-1)GF^{n-1}GE + (1-q^{-2})aF^{n-1}GE \\
    &\mapsto q^{-4}FEF^{n-1}EG - q^{-4}dFEF^{n} + q^{-2}fbFEF^{n-1} + (q^{-4}-1)GF^{n-1}GE\\
    &\quad+ (1-q^{-2})aF^{n-1}GE \text{ by the induction assumption} \\
    &\xmapsto{\sigma_{FE}^2} EF^{n-1}EFG + q^{-2}dEF^{n-1}G^2 + q^{-2}eaEF^{n-1}G + q^{-4}dGF^{n-1}EG\\
    &\quad+ q^{-4}eaF^{n-1}EG - q^{-4}dFEF^{n} + q^{-2}fbFEF^{n-1} + (q^{-4}-1)GF^{n-1}GE \\
    &\quad+ (1-q^{-2})aF^{n-1}GE\\
    &\xmapsto{\sigma_{EFG}}EF^{n-1}\Big(  -E^2 -q^{-4}F^2 - q^{-4}G^2 +cE + q^{-2}bF + q^{-2}aG \\
    &\quad+ \left(- q^{-4}(s^2 + t^2 + u^2 + v^2) - stuv + q^{-6}(q^2 +1)^2\right) \Big) \\
    &\quad+ q^{-4}dGF^{n-1}EG + q^{-4}eaF^{n-1}EG - q^{-4}dFEF^{n} + q^{-2}fbFEF^{n-1}\\
    &\quad+ (q^{-4}-1)GF^{n-1}GE + (1-q^{-2})aF^{n-1}GE \\
    &\xmapsto{\sigma_{GE}^2 \circ \sigma_{FE}^2} EF^{n-1}\Big(  -E^2 -F^2 - q^{-4}G^2 +cE + bF + q^{-2}aG \\
    &\quad+ \left(- q^{-4}(s^2 + t^2 + u^2 + v^2) - stuv + q^{-6}(q^2 +1)^2\right) \Big). \\
\end{align*}

For the final ambiguity \((\sigma_{EF^nG}, \sigma_{GF}, EF^n, G, F)\), one direction is:
\begin{align*}
    EF^nGF &\xmapsto{\sigma_{GF}} q^2 EF^{n+1}G + dEF^nE + ecEF^n \\
    &\xmapsto{EF^{n+1}G} Ff(n) + q^2(q^{-4}-1)GF^{n}G + q^2(1-q^{-2})aF^{n}G + dEF^nE\\
    &\quad+ ecEF^n.
\end{align*}
When \(n=1\) this gives 
\begin{align*}
    EFGF &\mapsto -FE^2 -q^{-4}F^3 - FG^2 +cFE + q^{-2}bF^2 + aFG \\
    &\quad+\left(- q^{-4}(s^2 + t^2 + u^2 + v^2) - stuv + q^{-6}(q^2 +1)^2\right)F \\
    &\quad+ q^2(q^{-4}-1)GFG + q^2(1-q^{-2})aFG + dEFE + ecEF \\
    &\xmapsto{\sigma_{FE}^3} -E^2F -q^{-2} dEG -q^{-2} eaE - dGE - eaE -q^{-4}F^3 - FG^2 + cEF\\
    &\quad+ dcG + eac + q^{-2}bF^2 + q^2aFG \\
    &\quad+\left(- q^{-4}(s^2 + t^2 + u^2 + v^2) - stuv + q^{-6}(q^2 +1)^2\right)F + q^2(q^{-4}-1)GFG \\
    &\xmapsto{\sigma_{GE} \circ \sigma_{GF}} -E^2F -q^{-2} dEG -q^{-2} eaE -q^{-2}dEG + q^{-2}d^2F - dfb \\
    &\quad- eaE -q^{-4}F^3 - FG^2 + cEF + dcG + eac + q^{-2}bF^2 + q^2aFG \\
    &\quad+\left(- q^{-4}(s^2 + t^2 + u^2 + v^2) - stuv + q^{-6}(q^2 +1)^2\right)F -q^2 dFG^2 \\
    &\quad- d^2EG - decG \\
    &= -E^2F + cEF -d^2EG +daE -q^{-4}F^3 + q^{-2}bF^2- q^4FG^2 + q^2aFG \\
    &\quad+ q^{-2}d^2F -q^2 d cG - dfb  + eac \\
    &\quad+\left(- q^{-4}(s^2 + t^2 + u^2 + v^2) - stuv + q^{-6}(q^2 +1)^2\right)F.  
\end{align*}

This equals the other direction when \(n=1\):
\begin{align*}
    EFGF &\xmapsto{\sigma_{EFG}} -E^2F -q^{-4}F^3 - G^2F +cEF + q^{-2}bF^2 + aGF \\
    &\quad+ \left(- q^{-4}(s^2 + t^2 + u^2 + v^2) - stuv + q^{-6}(q^2 +1)^2\right)F \\
    &\xmapsto{\sigma_{GE} \circ \sigma_{GF}^3} -E^2F -q^{-4}F^3 -q^4FG^2 -q^2 dEG -q^2 ecG -q^{-2}dEG \\
    &\quad+ q^{-2}d^2F - dfb - ecG +cEF + q^{-2}bF^2 + q^2 aFG + daE + eac \\
    &\quad+ \left(- q^{-4}(s^2 + t^2 + u^2 + v^2) - stuv + q^{-6}(q^2 +1)^2\right)F \\
    &= -E^2F +cEF -d^2EG + daE -q^{-4}F^3 + q^{-2}bF^2 -q^4FG^2 + q^2 aFG \\
    &\quad+ q^{-2}d^2F - (1+q^2) ecG  - dfb + eac \\
    &\quad+ \left(- q^{-4}(s^2 + t^2 + u^2 + v^2) - stuv + q^{-6}(q^2 +1)^2\right)F. \\  
\end{align*}
And in the general case:
\begin{align*}
    EF^nGF &\xmapsto{\sigma_{EF^nG}} q^{-2} Ff(n-1)F + (q^{-4}-1)GF^{n-1}GF + (1-q^{-2})aF^{n-1}GF \\
    &\mapsto  FEF^{n}G + q^{-2}dFEF^{n-1}E + q^{-2}ecFEF^{n-1} + (q^{-4}-1)GF^{n-1}GF \\
    &\quad+ (1-q^{-2})aF^{n-1}GF \text{ by the induction assumption}\\
    &\xmapsto{\sigma_{GF}^2 \circ \sigma_{FE}^2} FEF^{n}G + dEF^nE + q^{-2}d^2GF^{n-1}E + q^{-2}deaF^{n-1}E \\
    &\quad+ ecEF^{n} + q^{-2}decGF^{n-1} +q^{-2} e^2acF^{n-1} + q^2 (q^{-4}-1)GF^{n}G \\
    &\quad+ (q^{-4}-1)dGF^{n-1}E + (q^{-4}-1)ecGF^{n-1} + q^2 (1-q^{-2})aF^{n}G \\
    &\quad+ (1-q^{-2})daF^{n-1}E + (1-q^{-2})eacF^{n-1} \\
    &= FEF^{n}G + dEF^nE + ecEF^{n} -dGF^{n}G + q^2 (1-q^{-2})aF^{n}G \\
    &\xmapsto{\sigma_{EF^nG}} Ff(n) + dEF^nE + ecEF^{n} -dGF^{n}G + q^2 (1-q^{-2})aF^{n}G.
\end{align*}

Hence, all ambiguities in the reduction system are resolvable.  It remains to show that the reduction algorithm eventually terminates. We proceed by induction on the degree of the expression. As no rules apply to expressions of degree one, the reduction algorithm trivially terminates. Consider an expression \(T \in \mathbb{C} \langle E, F, G \rangle\) of degree \(n\); it is a finite linear combination of words in \(\langle E, F, G \rangle\) and can be reduced in a finite number of steps using the reduction rules \(\sigma_{FE}\), \(\sigma_{GE}\) and \(\sigma_{GF}\) to a finite linear combination of words of the form \(E^{\alpha_i} F^{\beta_i} G^{\gamma_i}\) for some \(\alpha_i, \beta_i, \gamma_i \in \mathbb{N}_0\) such that \(\alpha_i + \beta_i + \gamma_i \leq n\): if each of these monomials is reducible in a finite number of reductions so is \(T\). Either \(\beta_{i}=0\) and the monomial \(E^{\alpha_i} F^{\beta_i} G^{\gamma_i}\) is reduced, or the only reduction we can apply is \(\sigma_{EF^{\beta_i}G}\) which reduces the degree, so the result follows by induction. As every expression can be reduced fully in a finite number of reductions and all ambiguities are resolvable, the diamond lemma applies giving the result. 
 \end{proof}
\end{proposition}

\section*{Acknowledgements}
The author would like to thank David Jordan for his guidance in writing this paper and for his many comments and corrections on the drafts. The author would also like to thank Thomas Wright for his discussions and help with programming, Tim Weelinck for proofreading, Peter Samuelson for his assistance navigating the literature around DAHAs and skein theory, and Yuri Berest for encouraging her to include the puncture torus example. The research was funded through a EPSRC studentship, ERC grant STG-637618 and the F.R.S-FNRS.

\printbibliography[heading=bibintoc,title={References}]
\end{document}